\documentclass[aop,MSNbibl,seceqn,citesort,dvips]{arximspdf}

\doi{10.1214/11-AOP682}
\volume{40}
\issue{6}
\pubyear{2012}
\firstpage{2483}
\lastpage{2538}

\makeatletter

\newtheorem{theorem}{Theorem}[section]
\newtheorem{cor}[theorem]{Corollary}
\newtheorem{lemma}[theorem]{Lemma}
\newtheorem{prop}[theorem]{Proposition}
\newproclaim{defn}[theorem]{Definition}

\newproclaim{remark}[theorem]{Remark}

\newtheorem{conj}[theorem]{Conjecture}

\newcommand{\sL}{\mathcal{L}}

\newcommand{\bK}{{\mathbb K}}
\newcommand{\bN}{{\mathbb N}}
\newcommand{\bR}{{\mathbb R}}

\newcommand{\wt}{\widetilde}
\newcommand{\wtt}{\tilde}
\newcommand{\wh}{\widehat}
\newcommand{\E}{{\mathbb E}}
\newcommand{\PP}{{\mathbb P}}

\newcommand{\eps}{\varepsilon}

\makeatother

\begin{document}
\begin{frontmatter}

\title{Dirichlet heat kernel estimates for fractional Laplacian
with gradient perturbation}
\runtitle{Dirichlet heat kernel estimates
for $\Delta^{\alpha/2}+b\cdot\nabla$}

\begin{aug}
\author[A]{\fnms{Zhen-Qing} \snm{Chen}\thanksref{t1}\ead[label=e1]{zchen@math.washington.edu}},
\author[B]{\fnms{Panki} \snm{Kim}\corref{}\thanksref{t2}\ead[label=e2]{pkim@snu.ac.kr}} and
\author[C]{\fnms{Renming} \snm{Song}\ead[label=e3]{rsong@math.uiuc.edu}}
\runauthor{Z.-Q. Chen, P. Kim and R. Song}
\affiliation{University of Washington,
Seoul National University and University~of~Illinois}
\address[A]{Z.-Q. Chen\\
Department of Mathematics\\
University of Washington\\
Seattle, Washington 98195\\
USA\\
\printead{e1}}
\address[B]{P. Kim\\
Department of Mathematical Sciences \\
Seoul National University\\
San56-1 Shinrim-dong Kwanak-gu\\
Seoul 151-747\\
Republic of Korea\\
\printead{e2}}
\address[C]{R. Song\\
Department of Mathematics\\
University of Illinois\\
Urbana, Illinois 61801\\
USA\\
\printead{e3}}
\end{aug}

\thankstext{t1}{Supported in part by NSF Grants DMS-09-06743 and DMR-10-35196.}
\thankstext{t2}{Supported by Basic Science Research Program
through the National Research Foundation of Korea (NRF) grant funded
by the Korea government (MEST) (2010-0001984).}

\received{\smonth{11} \syear{2010}}
\revised{\smonth{5} \syear{2011}}

%
\begin{abstract}
Suppose that $d\geq2$ and $\alpha\in(1, 2)$. Let $D$ be a bounded
$C^{1,1}$ open set in $\bR^d$ and $b$ an $\bR^d$-valued function
on $\bR^d$ whose components are in a certain Kato class of the
rotationally symmetric $\alpha$-stable process. In this paper, we
derive sharp two-sided heat kernel estimates for
$\sL^b=\Delta^{\alpha/2}+b\cdot\nabla$ in $D$ with zero exterior
condition. We also obtain the boundary Harnack principle for $\sL^b$
in $D$ with explicit decay rate.
\end{abstract}

%

\pdfkeywords{60J35, 47G20, 60J75, 47D07,
Symmetric $\alpha$-stable process,
gradient operator,
heat kernel, transition density,
Green function, exit time,
Levy system, boundary Harnack inequality,
Kato class}

\end{frontmatter}

\section{Introduction}

Throughout this paper we assume $d\ge2$, $\alpha\in(1, 2)$ and
that $X$ is a (rotationally) symmetric $\alpha$-stable process on
$\bR^d$. The infinitesimal generator of $X$ is
$\Delta^{\alpha/2}:=-(-\Delta)^{\alpha/2}$.
We will use $B(x, r)$ to denote the open ball centered at $x\in
\bR^d$ with radius $r>0$.
\begin{defn}\label{dkc}
For a function $f$ on $\bR^d$, we define for $r>0$,
\[
M^\alpha_f(r)=\sup_{x\in\bR^d}\int_{B(x, r)}\frac{|f|(y)}
{|x-y|^{d+1-\alpha}} \,dy.
\]
A function $f$ on $\bR^d$ is said to belong to the Kato class
$\bK_{d, \alpha-1}$ if
\mbox{$ \lim_{r\downarrow0}M^\alpha_{f}(r)\!=\!0$}.
\end{defn}

Since $1<\alpha<2$, using H\"older's inequality, it is easy to see
that for every $p> d/(\alpha-1)$, $ L^\infty(\bR^d; dx)+L^p(\bR^d;
dx) \subset\bK_{d, \alpha-1}$. Throughout this paper we will assume
that $b=(b^1,\ldots, b^d)$ is an $\bR^d$-valued function on $\bR^d$
such that $|b|\in\bK_{d, \alpha-1}$. Define
$\sL^b=\Delta^{\alpha/2}+b\cdot\nabla$. Intuitively, the
fundamental solution $p^b(t, x, y)$ of $\sL^b$ and the fundamental
solution $p(t, x, y)$ of $\Delta^{\alpha/2}$, which is also the
transition density of $X$, should be related by the following
Duhamel formula:
%
%
\begin{equation}\label{e11c}
p^b (t, x, y)
=p (t, x, y)+
\int^t_0\int_{\bR^d}p^b (s, x, z) b(z)\cdot\nabla_z p (t-s, z,
y)\,dz\,ds.\hspace*{-25pt}
\end{equation}
Applying the above formula repeatedly, one expects that $p^b(t, x,
y)$ can be expressed as an infinite series in terms of $p$ and its
derivatives. This motivates the following definition. Define
$p^b_0(t, x, y)= p(t, x, y)$ and, for $k\geq1$,
%
%
\begin{equation}\label{e12}
p^b_k (t, x, y):= \int_0^t \int_{\bR^d} p^b_{k-1} (s, x, z) b(z)
\cdot
\nabla_z p(t-s, z, y) \,dz\,ds.
\end{equation}

The following results are shown in~\cite{BJ},
Theorem 1, Lemmas 15 and 23, and their proofs. Here and in the sequel, we use $:=$ as a
way of definition. For $a, b\in\bR$, $a\wedge b:=\min\{a, b\}$ and
$a\vee b:=\max\{a, b\}$.
\begin{theorem}\label{T11}
\textup{(i)} There exist\vspace*{1pt} $T_0>0$ and $c_1>1$ depending on $b$
only through the rate at which $ M^\alpha_{|b|}(r)$ goes to zero
such that $\sum_{k=0}^\infty p^b_k(t, x, y)$ converges locally
uniformly on $(0, T_0]\times\bR^d \times\bR^d$ to a positive
jointly continuous function $p^b(t, x, y)$ and that on $(0,
T_0]\times\bR^d \times\bR^d$,
%
%
\begin{equation}\label{e10}\quad
c_1^{-1}\biggl( t^{-d/\alpha} \wedge
\frac{t}{|x-y|^{d+\alpha}}\biggr) \le p^b (t, x, y) \le c_1 \biggl(
t^{-d/\alpha} \wedge\frac{t}{|x-y|^{d+\alpha}}\biggr).
\end{equation}
Moreover, $\int_{\bR^d} p^b(t, x, y) \,dy=1$ for every $t\in(0, T_0]$ and
$x\in\bR^d$.

\mbox{}\textup{\hphantom{i}(ii)} The function $p^b(t, x, y)$ defined in \textup{(i)} can be extended
uniquely to a positive jointly continuous function on $(0,
\infty) \times\bR^d \times\bR^d$
so that for all $s, t\in(0, \infty)$
and $(x,y) \in\bR^d \times\bR^d$,
$\int_{\bR^d} p^b(t, x, z) \,dz=1$ and
%
%
\begin{equation}\label{esemi}
p^b(s+t,x,y)=\int_{\bR^d} p^b(s,x,z)p^b(t,z,y) \,dz.
\end{equation}

\textup{(iii)} If we define
%
%
\begin{equation}\label{e11}
P^b_tf(x):= \int_{\bR^d} p^b(t,x,y) f(y)\,dy,
\end{equation}
then for any $f,g \in C^\infty_c(\bR^d)$, the space of smooth
functions with compact supports,
\[
\lim_{t \downarrow0}
\int_{\bR^d} t^{-1}
\bigl(P^b_tf(x)-f(x)\bigr) g(x)\,dx= \int_{\bR^d}
(\sL^b f)(x) g(x)\,dx.
\]
Thus $p^b(t,x,y)$ is the fundamental solution of $\sL^b$ in the
distributional sense.
\end{theorem}

Here and in the rest of this paper, the meaning of the phrase
``depending on $b$ only via the rate at which $M^\alpha_{|b|}(r)$
goes to zero'' is\vadjust{\goodbreak} that the statement is true for any $\bR^d$-valued
function $\wt b$ on $\bR^d$ with
\[
M^\alpha_{|\wtt b|} (r)\le M^\alpha_{|b|}(r)
\qquad\mbox{for all } r>0.
\]
Note that the Green function $G(x, y)$ of $X$ is $c/|x-y|^{d-\alpha}$
and so $|\nabla_x G(x, y)|\leq c/|x-y|^{d-\alpha+1}$.
This indicates that $\bK_{d, \alpha-1}$ is the right class
of functions for gradient perturbations of fractional Laplacian.
The same phenomenon happens for $\Delta+ b\cdot\nabla$; see~\cite{CrZ}.

It is easy to show (see Proposition~\ref{P23} below) that
the operators $\{P^b_t; t\geq0\}$ defined by~(\ref{e11}) form a
Feller semigroup and so there exists a conservative Feller process
$X^b=\{X^b_t, t\geq0, \PP_x, x\in\bR^d\}$ in $\bR^d$ such that
$P^b_t f(x)=\E_x[f(X^b_t)]$. The process $X^b$ is, in general,
nonsymmetric. We call $X^b$ an $\alpha$-stable process with
drift~$b$, since its infinitesimal generator is $\sL^b$.

For any open subset $D\subset\bR^d$, we define $\tau^b_D=\inf\{t>0\dvtx
X^b_t\notin D\}$. We will use $X^{b, D}$ to denote the subprocess of
$X^b$ in $D$; that is, $X^{b,
D}_t(\omega)=X^b_t(\omega)$ if $t< \tau^b_D(\omega)$ and $X^{b,
D}_t(\omega)=\partial$ if $t\geq\tau^b_D(\omega)$, where $\partial$
is a cemetery state. The subprocess\vspace*{1pt} of $X$ in $D$ will be denoted by
$X^D$. Throughout this paper, we use the convention that for every
function $f$, we extend its definition to $\partial$ by setting
$f(\partial)=0$. The infinitesimal generator of $X^{b, D}$ is
$\sL^b|_D$, that is, $\sL^b$ on $D$ with zero exterior condition.
The process $X^{b, D}$ has a transition density $p^b_D(t, x, y)$
with respect to the Lebesgue measure; see~(\ref{econti1}) below.
The transition density $p^b_D(t, x, y)$ of $X^{b, D}$ is the
fundamental solution of $\sL^b |_D$. The transition density of $X^D$
is denoted by $p_D(t, x, y)$, and it is the fundamental solution of
$\sL|_D$.

The purpose of this paper is to establish the following sharp
two-sided estimates on $p^{b}_D(t, x, y)$ in Theorem~\ref{Tmain}.
To state this theorem, we first recall that an open set $D$ in
$\bR^d$ is said to be a $C^{1,1}$ open set if there exist a
localization radius $R_0>0$ and a constant $\Lambda_0>0$ such that
for every $z\in\partial D$, there exist a $C^{1,1}$-function
$\phi=\phi_z\dvtx\bR^{d-1}\to\bR$ satisfying $\phi(0)=0$, $\nabla
\phi
(0)=(0,\ldots, 0)$, $\| \nabla\phi\|_\infty\leq\Lambda_0$, $|
\nabla\phi(x)-\nabla\phi(z)| \leq\Lambda_0 |x-z|$ and an
orthonormal coordinate system $CS_z$: $y=(y_1,\ldots, y_{d-1},
y_d):=(\wt y, y_d)$ with its origin at $z$ such that
\[
B(z,R_0)\cap D=\{y\in B(0, R_0) \mbox{ in } CS_z\dvtx y_d >\phi(\wt y)
\}.
\]
The pair $(R_0, \Lambda_0)$ is called the characteristics of the
$C^{1,1}$ open set $D$. We remark that in some literatures, the
$C^{1,1}$ open set defined above is called a \textit{uniform} $C^{1,1}$
open set as $(R_0, \Lambda_0)$ is universal for every $z\in\partial
D$. For $x\in D$, let $\delta_{D}(x)$ denote the Euclidean distance
between $x$ and $\partial D$. Note that a bounded $C^{1,1}$ open
set may be disconnected.
\begin{theorem}\label{Tmain} Let $D$ be a bounded $C^{1,1}$ open
subset of $ \bR^d$ with $C^{1,1}$ characteristics $(R_0,
\Lambda_0)$. Define
\[
f_D(t, x, y)=\biggl( 1\wedge
\frac{\delta_D(x)^{\alpha/2}}{\sqrt{t}}\biggr) \biggl( 1\wedge
\frac{\delta_D(y)^{\alpha/2}}{\sqrt{t}}\biggr) \biggl(
t^{-d/\alpha}
\wedge\frac{t}{|x-y|^{d+\alpha}}\biggr).
\]
For each $T>0$, there are constants $c_1=c_1(T, R_0, \Lambda_0, d,
\alpha, \operatorname{diam}(D), b) \geq1$ and $c_2=c_2(T, d, \alpha
, D, b)
\geq1$ with the dependence on $b$ only through the rate at
which $ M^\alpha_{|b|}(r)$ goes to zero such that:
\begin{longlist}
\item on $(0, T]\times D\times D$,
\[
c_1^{-1} f_D(t, x, y) \leq p^{b}_D(t, x, y)
\leq c_1 f_D(t, x, y);
\]
\item on $[T,
\infty)\times D\times D$,
\[
c_2^{-1} e^{-t \lambda_0^{b, D}} \delta_D (x)^{\alpha/2} \delta_D
(y)^{\alpha/2}\leq p^{b}_D(t, x, y) \leq c_2
e^{-t \lambda_0^{b, D}} \delta_D (x)^{\alpha/2}\delta_D
(y)^{\alpha/2},
\]
where $-\lambda_0^{b, D}:=\sup\operatorname{Re} (\sigma(\sL^b|_D)) <0$.
\end{longlist}
\end{theorem}

Here $\operatorname{diam}(D)$ denotes the diameter of $D$. At first glance,
one might think that the estimates in Theorem~\ref{Tmain} can be
obtained from the estimates for $p_D(t, x, y)$ by using a Duhamel
formula similar to~(\ref{e11c}) with $p^b$, $p$ and $\bR^d$
replaced by $p^b_D$, $p_D$ and $D$, respectively. Unfortunately
such an approach does not work for $p^b_D(t, x, y)$. This is because
unlike the whole space case, we do not have a good control on
$\nabla_z p_D(s, z, y)$ when $z$ is near the boundary of $D$.
When $D=\bR^d$, $p (t, x, y)$ is the transition density of the
symmetric $\alpha$-stable process, and there is a nice bound for
$\nabla_z p(t, z, y)$. This is the key reason why the result in
Theorem~\ref{T11}(i) can be established by using Duhamel's
formula. Instead, we establish Theorem~\ref{Tmain} by using
probabilistic means through the Feller process $X^b$. More
specifically, we adapt the road map outlined in our paper~\cite{CKS}
that establishes sharp\vspace*{1pt} two-sided Dirichlet heat kernel estimates for
symmetric $\alpha$-stable processes in $C^{1,1}$ open sets. Clearly,
many new and major difficulties arise when adapting the strategy
outlined in~\cite{CKS} to $X^b$. Symmetric stable processes are
L\'evy processes that are rotationally symmetric and self-similar.
The Feller process\vspace*{1pt} $X^b$ here is typically nonsymmetric, which is
the main difficulty that we have to overcome. In addition, $X^b$ is
neither self-similar nor rotationally symmetric. Specifically, our
approach consists of the following four ingredients:

\begin{longlist}
\item determine
the L\'evy system of $X^b$ that describes how the process jumps;

\item derive an approximate stable-scaling property of $X^b$
in bounded $C^{1,1}$ open sets, which will be used to derive heat
kernel estimates in bounded $C^{1,1}$ open sets for small time
$t\in(0, T]$ from that at time $t=1$;

\item
establish sharp two-sided estimates with explicit boundary decay
rate on the Green functions of $X^b$ and its suitable dual process
in $C^{1,1}$ open sets with sufficiently small diameter;

\item
prove the intrinsic ultracontractivity of (the nonsymmetric
process) $X^b$ in bounded open sets, which will give sharp two-sided
Dirichlet heat kernel estimates for large time.
\end{longlist}

In step (ii), we choose a large ball $E$ centered at the origin so
that our bounded $C^{1,1}$ open set $D$ is contained in $\frac14E$.
Then we derive heat kernel estimates in $D$ at time $t=1$ carefully
so that the constants depend on the quantity $M$ defined in
(\ref{MMM}), not on the diameter of $D$ directly. Note that the
constant $M$ has the correct scaling property, while the diameter of
$D$ does not. In fact, the constant $c_1$ in Theorem~\ref{Tmain}
depends on the diameter of $D$ only through $M$.

We also establish the boundary Harnack principle for $X^b$ and its
suitable dual process in $C^{1,1}$ open sets with explicit boundary
decay rate (Theorem~\ref{ubhp}). However, we like to point out that
Theorem~\ref{ubhp} is not used in the proof of Theorem~\ref{Tmain}.

By integrating the two-sided heat kernel estimates in Theorem
\ref{Tmain} with respect to $t$, one can easily get the
following estimates on the Green function $ G^b_D(x,
y)=\int_0^\infty p^b_D(t, x, y)\,dt$.
\begin{cor}\label{C12} Let $D$ be a bounded $C^{1,1}$ open set in
$\bR^d$. Then there is a constant $c=c(D, d, \alpha, b)\ge1$
with the dependence on $b$ only through the rate at which $
M^\alpha_{|b|}(r)$ goes to zero such that on $D\times D$,
\begin{eqnarray*}
&&c^{-1} \frac{1} { |x-y|^{d-\alpha}} \biggl(1\wedge\frac{
\delta_D(x)^{\alpha/2} \delta_D(y)^{\alpha/2}}{ |x-y|^{\alpha}}
\biggr)
\\
&&\qquad \leq G^b_D(x, y) \leq\frac{c} {|x-y|^{d-\alpha}}
\biggl(1\wedge\frac{ \delta_D(x)^{\alpha/2}
\delta_D(y)^{\alpha/2}}{ |x-y|^{\alpha}} \biggr).
\end{eqnarray*}
\end{cor}

The above result was obtained independently as the main result in
\cite{BJ2}. Clearly the heat kernel $p^{b}_D(t, x, y)$ contains much
more information than the Green function $G^b_D(x, y)$. The
estimates on $p^b_D(t, x, y)$ are not studied in~\cite{BJ2}.

The sharp two-sided estimates for $p_D(t, x, y)$, corresponding to
the case $b=0$ in Theorem~\ref{Tmain}, were first established in
\cite{CKS}. Theorem~\ref{Tmain} indicates that short time Dirichlet
heat kernel estimates for the fractional Laplacian in bounded
$C^{1,1}$ open sets are stable under gradient perturbations. Such
stability should hold for much more general open sets.

We say that an open set $D$ is $\kappa$-fat if there exists an
$R_1>0$ such that for every $x\in D$ and $r\in(0, R_1]$, there is
some $y$ such that $B(y, \kappa r)\subset B(x, r)\cap D$. The pair
$(R_1, \kappa)$ is called the characteristics of the $\kappa$-fat
open set $D$.
\begin{conj}
Let $T>0$ and $D$ be a bounded $\kappa$-fat open subset of $\bR^d$.
Then there is a constant $c_1\geq1$ depending only on $T$, $D$,
$\alpha$ and $b$ with the dependence on $b$ only through the
rate at which $M^\alpha_{|b|}(r)$ goes to zero such that
\[
c_1^{-1} p_D(t, x, y)
\leq p^b(t, x, y) \leq c_1 p_D(t, x, y)\qquad
\mbox{for } t\in(0, T] \mbox{ and } x, y\in D
\]
and
\[
c_1^{-1} G_D( x, y)
\leq G_D^b( x, y) \leq c_1 G_D( x, y)\qquad
\mbox{for } x, y\in D.
\]
\end{conj}

In the remainder of this paper, the constants $C_1, C_2, C_3, C_4$
will be fixed throughout this paper. The lower case constants $c_0,
c_1, c_2, \ldots$ can change from one appearance to another. The
dependence of the constants on the dimension $d$ and the stability
index $\alpha$ will not be always mentioned explicitly. We will use
$dx$ to denote the Lebesgue measure in $\bR^d$. For a Borel set
$A\subset\bR^d$, we also use $|A|$ to denote its Lebesgue measure.
The space of continuous functions on $\bR^d$ will be denoted as
$C(\bR^d)$, while $C_b(\bR^d)$ and $C_\infty(\bR^d)$ denote the
space of bounded continuous functions on $\bR^d$ and the space of
continuous functions on $\bR^d$ that vanish at infinity,
respectively.
For two nonnegative functions $f$ and $g$, the notation $f\asymp g$
means that there are positive constants $c_1$ and $c_2$ so that
$c_1g( x)\leq f (x)\leq c_2 g(x)$ in the common domain of definition
for $f$ and $g$.

\section{\texorpdfstring{Feller property and L\'evy system}{Feller property and Levy system}}\label{SKP}

Recall that $d\ge2$ and $\alpha\in(1, 2)$. A~(rotationally) symmetric
$\alpha$-stable process $X=\{X_t, t\geq0, \PP_x, x\in\bR^d\}$ in
$\bR^d$ is a L\'evy process such that
\[
\E_x \bigl[ e^{i\xi\cdot(X_t-X_0)} \bigr] = e^{-t|\xi|^{\alpha}}
\qquad\mbox{for every } x\in\bR^d \mbox{ and } \xi\in\bR^d.
\]
The infinitesimal generator of this process
$X$
is the fractional Laplacian $ \Delta^{\alpha/2}$,
which is a prototype of nonlocal operators. The fractional Laplacian
can be written in the form
%
%
\begin{equation}\label{efraclap}
\Delta^{\alpha/2} u(x) = \lim_{\eps\downarrow0}\int_{\{y\in
\bR^d\dvtx |y-x|>\eps\}} \bigl(u(y)-u(x)\bigr) \frac{\mathcal{A}(d, - \alpha)
}{|x-y|^{d+\alpha}} \,dy,
\end{equation}
where $ \mathcal{A}(d, - \alpha):=
\alpha2^{\alpha-1}\pi^{-d/2} \Gamma(\frac{d+\alpha}2)
\Gamma(1-\frac{\alpha}2)^{-1}. $

We will use $p(t, x, y)$ to denote the transition density of $X$ (or
equivalently the heat kernel of the fractional Laplacian $
\Delta^{\alpha/2}$). It is well known (see, e.g.,~\cite{BG,CK})
that
\[
p (t, x, y) \asymp t^{-d/\alpha} \wedge
\frac{t}{|x-y|^{d+\alpha}} \qquad\mbox{on } (0,
\infty)\times\bR^d \times\bR^d.
\]

The next two lemmas will be used
later.
\begin{lemma}\label{lkato1}
If $f$ is a function belonging to
$\bK_{d, \alpha-1}$, then for any compact subset $K$ of $\bR^d$,
\[
\sup_{x\in\bR^d}\int_{K}\frac{|f|(y)} {|x-y|^{d-\alpha}}
\,dy<\infty.
\]
\end{lemma}
\begin{pf}
This follows immediately from the fact that
$d-\alpha<d+1-\alpha$. We omit the details.\vadjust{\goodbreak}
\end{pf}

Recall that we are assuming that $b$ is
an $\bR^d$-valued function on $\bR^d$ such that $|b|\in\bK_{d,
\alpha-1}$.
\begin{lemma}\label{lkato2}
If $f$ is a function belonging to
$\bK_{d, \alpha-1}$, then
\[
\lim_{t\to0}\sup_{x\in\bR^d}\int^t_0P^b_s|f|(x)\,ds=0.
\]
\end{lemma}
\begin{pf}
By~(\ref{e10}),
\begin{eqnarray*}
&&\lim_{t\to0}\sup_{x\in\bR^d}\int^t_0P^b_s|f|(x)\,ds \\
&&\qquad\leq c_1
\lim_{t\to0}\sup_{x\in\bR^d}\int^t_0 \biggl(s\int_{B(x,
s^{1/\alpha})^c}\frac{|f(y)|}{|y-x|^{d+\alpha}}\,dy\\
&&\qquad\quad\hspace*{68.6pt}{} +
s^{-d/\alpha}\int_{B(x, s^{1/\alpha})}|f(y)|\,dy\biggr) \,ds.
\end{eqnarray*}
So it suffices to show that the right-hand side is zero.
Clearly, for any $s\le1$, we have
%
%
\begin{equation}\label{ekato3}
\int_{B(x, s^{1/\alpha})}|f(y)|\,dy\le
(s^{1/\alpha})^{d+1-\alpha}\sup_{x\in\bR^d}\int_{B(x,
1)}\frac{|f(y)|}{|y-x|^{d+1-\alpha}}\,dy.
\end{equation}
Applying~\cite{So1}, Lemma 1.1, we get
%
%
\begin{equation}\label{ekato4}
\sup_{x\in\bR^d}\int_{B(x,
s^{1/\alpha})^c}\frac{|f(y)|}{|y-x|^{d+\alpha}}\,dy \le
c_2(s^{1/\alpha})^{d+1-\alpha} (s^{1/\alpha})^{-(d+\alpha)} =
c_2s^{1/\alpha-2}.\hspace*{-35pt}
\end{equation}
Now the conclusion follows immediately from
(\ref{ekato3}) and~(\ref{ekato4}).
\end{pf}

By the\vspace*{1pt} semigroup property of $p^b (t, x, y)$
and~(\ref{e10}), there are constants
$c_1$, \mbox{$c_2\geq1$}, such that on $(0, \infty)\times\bR^d \times\bR^d$,
%
%
\begin{eqnarray}\label{e11a}
&&c^{-1}_1 e^{-c_2t} \biggl( t^{-d/\alpha} \wedge
\frac{t}{|x-y|^{d+\alpha}}\biggr) \nonumber\\[-8pt]\\[-8pt]
&&\qquad\le p^b (t, x, y) \le c_1
e^{c_2t} \biggl( t^{-d/\alpha} \wedge
\frac{t}{|x-y|^{d+\alpha}}\biggr) .\nonumber
\end{eqnarray}

\begin{prop}\label{P23}
The family of operators $\{P^b_t; t\geq0\}$ defined by
(\ref{e11}) forms a Feller semigroup. Moreover, it satisfies the
strong Feller property; that is, for each $t>0$, $P^b_t f$ maps
bounded measurable functions to continuous functions.
\end{prop}
\begin{pf}
Since $p^b(t, x, y)$ is continuous, by the bounded convergence
theorem, $P^b_t$ enjoys the strong Feller property.\vadjust{\goodbreak} Moreover, for
every $f\in C_\infty(\bR^d)$ and $t>0$,
\[
\lim_{x\to\infty} |P^b_t f(x)|
\leq\lim_{x\to\infty} c_1
e^{c_2t} \int_{\bR^d} \biggl( t^{-d/\alpha} \wedge
\frac{t}{|y|^{d+\alpha}}\biggr) |f(x+y)| \,dy =0
\]
and so $P^b_t f \in C_\infty(\bR^d)$. By~(\ref{e11a}), we have
\begin{eqnarray*}
&&\sup_{t\le t_0} \sup_{ x \in\bR^d}
\PP_x(|X^b_t -X^b_0| \ge\delta)\\
&&\qquad\le c_1 e^{c_2 t_0} \sup_{t\le t_0} \sup_{ x \in\bR^d}
\int_{ \{y\in\bR^d\dvtx |x-y| \ge\delta\}} \biggl( t^{-d/\alpha}
\wedge
\frac{t}{|x-y|^{d+\alpha}}\biggr) \,dy\\
&&\qquad= c_3 e^{c_2 t_0} \sup_{t\le t_0} \int^\infty_{\delta}
r^{d-1}
\biggl( t^{-d/\alpha} \wedge
\frac{t}{r^{d+\alpha}}\biggr) \,dr\\
&&\qquad\le c_4 e^{c_2 t_0} \int^\infty_{\delta t_0^{-1/\alpha}}
u^{d-1} \biggl( 1 \wedge
\frac{1}{u^{d+\alpha}}\biggr) \,du
\end{eqnarray*}
for some $c_3=c_3(d)>0$ and $c_4=c_4(d)>0$. Thus
%
%
\begin{equation}\label{econ}
\lim_{t_0 \downarrow0} \sup_{t\le t_0}\sup_{x \in\bR^d}
\PP_x(|X^b_t -X^b_0| \ge\delta)
= 0.
\end{equation}
For every $f\in C_b(\bR^d)$, $x\in\bR^d$ and $\eps>0$, there is a
$\delta>0$ so that $|f(y)-f(x)|<\eps$ for every $y\in B(x,
\delta)$. Therefore we have by~(\ref{econ}),
\begin{eqnarray*}
&&\lim_{t\downarrow0} |P^b_t f(x) -f(x)|\\
&&\qquad= \lim_{t\downarrow0} \biggl|\int_{\bR^d} p^b(t, x, y) \bigl(f(y)-f(x)\bigr)
\,dy \biggr| \\
&&\qquad\leq \lim_{t\downarrow0} \int_{\{y\in\bR^d\dvtx |y-x|<\delta\}}
p^b(t, x, y)
|f(y)-f(x)| \,dy \\
&&\qquad\quad{} + \lim_{t\downarrow0} 2\| f\|_\infty
\PP_x(|X^b_t -x| \ge\delta) < \eps.
\end{eqnarray*}
Therefore for every $f\in C_b(\bR^d)$ and $x\in\bR^d$,
$\lim_{t\downarrow0} P^b_t f(x)=f(x)$. This completes the proof of
the proposition.
\end{pf}

We will need the next result, which is an extension of Theorem
\ref{T11}(iii).
\begin{prop}\label{peffee}
For any $f\in C^\infty_c(\bR^d)$ and $g\in C_\infty(\bR^d)$, we have
\[
\lim_{t \downarrow0} \int_{\bR^d} t^{-1} \bigl(P^b_tf(x)-f(x)\bigr) g(x)\,dx=
\int_{\bR^d} (\sL^b f)(x) g(x)\,dx.
\]
\end{prop}
\begin{pf}
This proposition can be proved by following the proof of
\cite{BJ}, Theorem 1, with some obvious modifications. Indeed, one
can follow the same argument of the proof of~\cite{BJ}, Theorem 1,
until the second display in~\cite{BJ},\vadjust{\goodbreak} page 195, with $f\in
C^\infty_c(\bR^d)$ and $g\in C_\infty(\bR^d)$. Let $\eps>0$, and use
the same notation as in~\cite{BJ}, page 195, except that $K:=\{z\dvtx
\operatorname{dist} (z, K_1) \le1\}$ and we ignore $K_2$. Since $h(x,y)=
\nabla f(y) g(x)$ is still uniformly continuous, there exists a
$\delta>0$ such that for all $x,y,z$ with $|x-z| < \delta$ and
$|y-z| < \delta$, we have that $|h(x,y) -h(z,z)| < \eps$. Thus the
third display in~\cite{BJ}, page 195, can be modified as
\begin{eqnarray*}
\hspace*{-5pt}&& \biggl| I_t - \int_{\bR^d} b(z) \cdot\nabla f(z) g(z) \,dz \biggr|\\
\hspace*{-5pt}&&\qquad\le
\int_{\bR^d}\int_{\bR^d}\int_{\bR^d}
\int_0^t \hspace*{-0.5pt}\frac{p(t-s, x,z)p(s,z,y)}{t} \,ds |b(z) ||h(x,y) -h(z,z)|
\,dx\,dy\,dz\\
\hspace*{-5pt}&&\qquad\le
2\|h\| \int_{K^c}\int_{K_1}\int_0^t \biggl(\int_{\bR^d}p(t-s, x,z)\,dx\biggr)
\frac{p(s,z,y)}{t} \,ds |b(z) |
\,dy\,dz\\
\hspace*{-5pt}&&\qquad\quad{}+
2\|h\|\int_{K}\int\int_{(B(z, \delta) \times B(z, \delta))^c}
\int_0^t \frac{p(t-s, x,z)p(s,z,y)}{t} \,ds |b(z) |
\,dx\,dy\,dz\\
\hspace*{-5pt}&&\qquad\quad{}+\eps\int_{K}\int\int_{B(z, \delta) \times B(z, \delta)}
\int_0^t \frac{p(t-s, x,z)p(s,z,y)}{t} \,ds |b(z) |
\,dx\,dy\,dz .
\end{eqnarray*}
The remainder of the proof is the same as that of the proof of
\cite{BJ}, Theorem~1.
\end{pf}

The Feller process $X^b$ possesses a L\'evy system (see
\cite{Sh}), which describes how $X^b$ jumps.
Intuitively, since the infinitesimal
generator of $X^b$ is $\sL^b$, $X^b$ should satisfy
\[
dX^b_t=dX_t + b(X^b_t) \,dt.
\]
So $X^b$ should have the same L\'evy system as that of $X$, as the
drift does not contribute to the jumps. This is indeed true, and we
are going to give a rigorous proof.

It is well known that the symmetric stable process $X$ has L\'evy
intensity function
\[
J(x, y)=\mathcal{A}(d, -\alpha)|x-y|^{-(d+\alpha)}.
\]
The L\'evy intensity function gives rise to a L\'evy system $(N, H)$
for $X$, where $N(x, dy)=J(x, y)\,dy$ and $H_t=t$, which describes the
jumps of the process $X$: for any $x\in\bR^d$ and any nonnegative
measurable function $f$ on $\bR_+ \times\bR^d\times\bR^d$
vanishing on $\{(s, x, y)\in\bR_+ \times\bR^d\times\bR^d\dvtx x=y\}$
and stopping time $T$ (with respect to the filtration of $X$),
\[
\E_x \biggl[\sum_{s\le T} f(s,X_{s-}, X_s) \biggr]= \E_x \biggl[
\int_0^T \biggl( \int_{\bR^d} f(s,X_s, y) J(X_s, y)\,dy \biggr) \,ds
\biggr].
\]
(See, e.g.,~\cite{CK}, proof of Lemma 4.7, and~\cite{CK2},
Appendix A.)

We first show that $X^b$ is a solution to the martingale problem of
$\sL^b$.
\begin{theorem}\label{T25}
For every $x\in\bR^d$ and every $f\in C^\infty_c(\bR^d)$,
\[
M^f_t:=f(X^b_t)-f(X^b_0)-\int^t_0 \sL^b f(X^b_s) \,ds
\]
is a martingale under $\PP_x$.
\end{theorem}
\begin{pf}
Define the adjoint operator $P^{b, *}_t$ of $P^b_t$ with
respect to the Lebesgue measure by
\[
P^{b, *}_tf(x):= \int_{\bR^d} p^b(t,y,x) f(y)\,dy.
\]
It follows immediately from~(\ref{e10}) and the continuity of
$p^b(t,x,y)$ that, for any $g\in C_\infty(\bR^d)$ and $s>0$, both
$P^b_sg$ and $P^{b, *}_sg$ are in $C_\infty(\bR^d)$. Thus, for any
$f, g\in C^\infty_c(\bR^d)$ and $s>0$, by applying Proposition
\ref{peffee} with $h=P^{b, *}_sg$ and~(\ref{esemi}), we get that
\begin{eqnarray*}
&&\lim_{t \downarrow0} \int_{\bR^d} t^{-1} \bigl(P^b_{t+s}f(x)-P^b_sf(x)\bigr)
g(x)\,dx\\
&&\qquad=
\lim_{t \downarrow0} \int_{\bR^d} t^{-1} \bigl(P^b_{t}f(x)-f(x)\bigr)
P^{b, *}_sg (x)\,dx \\
&&\qquad=\int_{\bR^d} \sL^b f (x)
P^{b, *}_sg (x)\,dx = \int_{\bR^d}
\E_x[
\sL^b f (X^b_s) ] g(x)\,dx,
\end{eqnarray*}
which implies that
%
%
\begin{equation}\label{enew11}
\int_{\bR^d}\bigl(P^b_{t}f(x)-f(x)\bigr) g(x)\,dx = \int_{\bR^d}\E_x\biggl[
\int^t_0 (\sL^b f) (X^b_s) \,ds \biggr] g(x)\,dx.
\end{equation}
Using the strong Feller property of $P^b_t$, Lemmas~\ref{lkato1}
and~\ref{lkato2}, we can easily see that the function
\[
x\mapsto
P^b_tf(x)-f(x)-\E_x\biggl[ \int_0^t \sL^b f (X^b_s)\,ds \biggr]
=\E_x [ M^f_t ]
\]
is continuous, and thus is identically zero on $\bR^d$ by
(\ref{enew11}). It follows that for any $f\in C^\infty_c(\bR^d)$
and $x\in\bR^d$, $M^f$ is a martingale with respect to $\PP_x$.
\end{pf}

Theorem~\ref{T25} in particular implies that $X^b_t=(X^{b,
1}_t,\ldots, X^{b, d}_t) $ is a semi-martingale. By It\^{o}'s formula,
we have that, for any $f\in C^\infty_c(\bR^d)$,
%
%
\begin{equation}\label{eito}
f(X^b_t)-f(X^b_0)=\sum^d_{i=1}\int^t_0{\partial}_if(X^b_{s-})\,dX^{b,
i}_s +\sum_{s\le t}\eta_s(f) +\frac12 A_t(f),
\end{equation}
where
%
%
\begin{equation}\label{eito2}
\eta_s(f)=f(X^b_s)-f(X^b_{s-})-\sum^d_{i=1}{\partial}_if(X^b_{s-})(X^{b,
i}_s-X^{b, i}_{s-})\vadjust{\goodbreak}
\end{equation}
and
%
%
\begin{equation}\label{eito3}
A_t(f)=\sum^d_{i, j=1}\int^t_0{\partial}_i{\partial}_jf(X^b_{s-})\,d
\langle(X^{b, i})^c, (X^{b, j})^c\rangle_s.
\end{equation}

Now suppose that $A$ and $B$ are two bounded closed sets having a
positive distance from each other. Let $f\in C^\infty_c(\bR^d)$ with
$f=0$ on $A$ and $f=1$ on~$B$. Then we know that
$N^f_t:=\int^t_0{\mathbf 1}_A(X^b_{s-})\,dM^f_s$ is a martingale.
Combining Theorem~\ref{T25} and~(\ref{eito})--(\ref{eito3}) with
(\ref{efraclap}), we get that
\begin{eqnarray*}
N^f_t&=&\sum_{s\le t}{\mathbf1}_A(X^b_{s-})f(X^b_s) -
\int^t_0{\mathbf1}_A(X^b_s)(
\Delta^{\alpha/2}
f(X^b_s))\,ds\\[-2pt]
&=&\sum_{s\le
t}{\mathbf1}_A(X^b_{s-})f(X^b_s)-\int^t_0{\mathbf1}_A(X^b_s)\int_{\bR
^d}f(y)J(X^b_s,
y)\,dy\,ds.
\end{eqnarray*}
By taking a sequence of functions $f_n\in C^\infty_c(\bR^d)$ with
$f_n=0$ on $A$, $f_n=1$ on $B$ and $f_n\downarrow{\mathbf1}_B$, we get
that, for any $x\in\bR^d$,
\[
\sum_{s\le t}{\mathbf1}_A(X^b_{s-}){\mathbf1}_B(X^b_s) -\int^t_0{\mathbf1}_A(X^b_s)\int_BJ(X^b_s, y)\,dy\,ds
\]
is a martingale with respect to $\PP_x$. Thus,
\[
\E_x\biggl[ \sum_{s\le t}{\mathbf1}_A(X^b_{s-}){\mathbf1}_B(X^b_s)
\biggr]=
\E_x\biggl[\int^t_0\int_{\bR^d} {\mathbf1}_A(X^b_s){\mathbf1}_B(y)J(X^b_s, y)\,dy\,ds\biggr].
\]
Using this and a routine measure theoretic argument, we get
\[
\E_x\biggl[ \sum_{s\le
t}f(X^b_{s-}, X^b_s) \biggr]
=\E_x\biggl[\int^t_0\int_{\bR^d}f(X^b_s, y)J(X^b_s, y)\,dy\,ds\biggr]
\]
for any nonnegative measurable function $f$ on $\bR^d\times\bR^d$
vanishing on $\{(x, y)\in\bR^d\times\bR^d: x=y\}$. Finally
following the same arguments as in~\cite{CK}, Lemma~4.7, and
\cite{CK2}, Appendix A, we get:\vspace*{-2pt}
\begin{theorem}\label{tls4xb}
$X^b$ has the same L\'evy system $(N, H)$ as $X$; that is, for any
$x\in\bR^d$ and any nonnegative measurable function $f$ on $\bR_+
\times\bR^d\times\bR^d$ vanishing on $\{(s, x, y)\in\bR_+ \times
\bR^d\times\bR^d: x=y\}$ and stopping time $T$ (with respect to
the filtration of $X^b$)
%
%
\begin{equation}\label{els4xb}
\E_x \biggl[\sum_{s\le T} f(s,X^b_{s-}, X^b_s) \biggr]= \E_x \biggl[
\int_0^T \biggl( \int_{\bR^d} f(s,X^b_s, y) J(X^b_s, y)\,dy \biggr) \,ds
\biggr].\hspace*{-35pt}
\end{equation}
\end{theorem}

For any open subset $E$ of $\bR^d$, let $E_\partial= E\cup\{
\partial\}$, where $\partial$ is the cemetery point. Define for $x,
y\in
E$,
\[
N^E(x, dy):=J(x, y)\,dy,\qquad
N^E(x, \partial):=\int_{E^c}J(x, y)\,dy\vadjust{\goodbreak}
\]
and $H^E_t:=t$. Then it follows from the theorem above that $(N^E,
H^E)$ is a L\'evy system for $X^{b, E}$; that is, for any $x\in E$,
any nonnegative measurable function $f$ on $\bR_+ \times E\times
E_\partial$ vanishing on $\{(s, x, y)\in\bR_+ \times E\times
E\dvtx x=y\}$ and stopping\vspace*{1pt} time $T$ (with respect to the filtration of
$X^{b, E}$)
%
%
\begin{eqnarray}\label{elevysE}
&&\E_x \biggl[\sum_{s\le T} f(s,X^{b,E}_{s-}, X^{b,E}_s) \biggr]\nonumber\\[-8pt]\\[-8pt]
&&\qquad= \E_x
\biggl[ \int_0^T \biggl(
\int_{E_\partial}
f(s,X^{b,E}_s, y) N^E(X^{b,E}_s,
dy) \biggr) \,dH^E_s \biggr].\nonumber
\end{eqnarray}

\section{Subprocess of $X^b$}

In this section we study some basic properties of subprocesses of
$X^b$ in open subsets. These properties will be used in later
sections.
\begin{lemma}\label{lcon}
For any $\delta>0$, we have
\[
\lim_{s\downarrow0} \sup_{x \in\bR^d} \PP_x\bigl( \tau^b_{B(x, \delta)}
\le s\bigr) =0.
\]
\end{lemma}
\begin{pf}
By the strong Markov\vspace*{1pt} property of $X^b$ (see, e.g.,
\cite{BG1}, pages 43 and 44), we have for every $ x \in\bR^d $,
%
%
\begin{eqnarray}\label{enew12}
&&\PP_x\bigl( \tau^b_{B(x, \delta)} \le s\bigr)
\nonumber\\
&&\qquad\le \PP_x\bigl(\tau^b_{B(x,
\delta)} \le s, X^b_s \in B(x, \delta/2) \bigr) +
\PP_x\bigl(
X^b_s \in B(x, \delta/2)^c \bigr)\nonumber\\
&&\qquad\le \E_x\bigl[ \PP_{X_{\tau^b_{B(x, \delta)}}} \bigl(
\bigl|X^b_{s-\tau^b_{B(x, \delta)}} -X^b_0\bigr|\geq\delta/2\bigr)
; \tau^b_{B(x, \delta)} < s \bigr]
\\
&&\qquad\quad{} + \PP_x( |X^b_s -X^b_0|\geq\delta/2
)\nonumber\\
&&\qquad\leq 2 \sup_{t\leq s} \sup_{x\in\bR^d}
\PP_x ( |X^b_t- X^b_0|\geq\delta/2 ).
\nonumber
\end{eqnarray}
Now the conclusion of the lemma follows from~(\ref{econ}).
\end{pf}

A point $z$ on the boundary $\partial G$ of a Borel set $G$ is said
to be a regular boundary point with respect to $X^b$ if
$\PP_z(\tau^b_G=0)=1$. A Borel set $G$ is said to be regular with
respect to $X^b$ if every point in $\partial G$ is a regular
boundary point
with respect to $X^b$.
\begin{prop}\label{preg}
Suppose that $G$ is a Borel set of $\bR^d$ and $z\in
\partial G$. If there is a cone $A$ with vertex $z$ such that
$\operatorname{int}(A)\cap
B(z, r)\subset G^c$ for some $r>0$, then $z$ is a regular boundary
point of $G$ with respect to $X^b$.
\end{prop}
\begin{pf}
This result follows from~(\ref{e10}) and Blumenthal's
zero--one law by a routine argument. For example, the reader can
follow the argument in the\vadjust{\goodbreak} proof of~\cite{KS1}, Proposition 2.2.
Even though~\cite{KS1}, Proposition 2.2, is stated for open sets, the
proof there works for Borel sets. We omit the details.
\end{pf}

This result implies that all bounded Lipschitz open sets, and in
particular, all bounded $C^{1,1}$ open sets, are regular with
respect to $X^b$. Repeating the argument in the second part of the
proof of~\cite{CZ}, Theorem 1.23, we immediately get the following
result.
\begin{prop}\label{pcontb}
Suppose that $D$ is an open set in $\bR^d$, and $f$ is a bounded
Borel function on $\partial D$. If $z$ is a regular boundary point
of $D$ with respect to $X^b$, and $f$ is continuous at $z$, then
\[
\lim_{\bar{D}\ni x\to z}\E_x[f(X^b_{\tau^b_D});
\tau^b_D<\infty]=f(z).
\]
\end{prop}

Let
%
%
\begin{equation}\label{eqdipco}
k_D^b(t,x,y) := \E_x[p^b(t-\tau^b_D, X^b_{\tau^b_D},y);
\tau^b_D < t
]
\end{equation}
and
%
%
\begin{equation}\label{eqdipco1}
p^b_D(t,x,y):= p^b(t,x,y)-k^b_D(t,x,y).
\end{equation}
Then $p^b_D(t,x,y)$ is the transition density of $X^{b, D}$. This is
because by the strong Markov property of $X^b$, for every $t>0$ and
Borel set $A\subset\bR^d$,
%
%
\begin{equation}\label{econti1}
\PP_x(X^{b, D}_t \in A) = \int_{A} p^b_D(t,x,y) \,dy.
\end{equation}
We will use $\{P^{b, D}_t\}$ to denote the semigroup of $X^D$ and
$\sL^b|_D$ to denote the infinitesimal generator of $\{P^{b, D}_t\}$.
Using some standard arguments
(e.g., \cite{B,CZ}),
we can show the following:
\begin{theorem}\label{tconti}
Let $D$ be an open set in $\bR^d$.
The transition density $p^b_D(t,x,y)$
is jointly continuous on $(0, \infty) \times D \times D$. For every
$t>0$ and $s > 0$,
%
%
\begin{equation}\label{eck}
p^b_D(t+s,x,y) = \int_{D} p^b_D(t,x,z)p^b_D(s,z,y) \,dz.
\end{equation}
If $z$ is a regular boundary point of $D$ with respect to $X^b$,
then for any $t>0$ and $y\in D$, $ \lim_{D\ni x\to z}p^b_D(t, x,
y)=0$.
\end{theorem}
\begin{pf}
Note that by~(\ref{e11a}), there exist $c_1, c_2>0$ such that
for all $t_0>0$ and $\delta>0$,
%
%
\begin{eqnarray}\label{ebd1}
\sup_{t\le t_0} \sup_{|x-y|\geq\delta} p^b(t , x,y)
&\le& c_1 e^{c_2 t_0}
\sup_{t \le t_0} \sup_{|x-y| \geq\delta}
\biggl( t^{-d/\alpha} \wedge\frac{t}{|x-y|^{d+\alpha}}\biggr)
\nonumber\\[-8pt]\\[-8pt]
&\le& c_1 e^{c_2 t_0} \frac{t_0}{\delta^{d+\alpha}} < \infty.
\nonumber
\end{eqnarray}
We first show that $k^b_D(t,x, \cdot)$ is jointly continuous on
$(0, \infty)\times D \times D$. For any $\delta>0$, define
$D_\delta=\{x\in D\dvtx \operatorname{dist}(x, D^c) <\delta\}$. For $0\leq
s<r$ and $x, y\in D_\delta$, define
\[
h(s, r, x, y)= \E_x [ p^b(r-\tau^b_D, X^b_{\tau^b_D},y);
s\leq\tau^b_D < r ].
\]
Note that
\begin{eqnarray*}
\E_x [h(s, r, X^b_s, y)] &=&\E_x[ h(s, r, X^b_s, y); s<\tau^b_D]
+\E_x[ h(s, r, X^b_s, y); s\geq\tau^b_D] \\[-2pt]
&=& h(s, r+s, x, y)+\E_x[ h(s, r, X^b_s, y); s\geq\tau^b_D]
\end{eqnarray*}
and
\begin{eqnarray*}
k_D^b(t,x,y) &=& h(0, t, x, y) \\[-2pt]
&=& h(s, t, x, y) + \E_x[p^b(t-\tau^b_D, X^b_{\tau^b_D},y);
\tau^b_D < s ] \\[-2pt]
&=& \E_x [h(s, t-s, X^b_s, y)] - \E_x[ h(s, t-s, X^b_s, y);
\tau^b_D\leq s] \\[-2pt]
&&{} + \E_x[p^b(t-\tau^b_D, X^b_{\tau^b_D},y);
\tau^b_D < s ].
\end{eqnarray*}
For all $t_1, t_2\in(0, \infty)$, by~(\ref{ebd1}), $p^b(r, z, y)$
is bounded on $(0, t_2]\times D^c \times D_\delta$ by a constant
$c_3$. Consequently, $h(s, r, x, y)$ is bounded by $c_3$ for all
$x, y\in D_\delta$ and $s, r\in(0, t_2]$ with $s<r\wedge(t_1/3)$.
Thus we have from the above display as well as~(\ref{ebd1}) that
for all $t\in[t_1, t_2]$, $s<t_1/2$ and $x, y\in D_\delta$,
\begin{eqnarray*}
|k_D^b(t,x,y) - \E_x [h(s, t-s, X^b_s, y)]|
&\leq&2c_3 \PP_x (\tau^b_D \leq s)\\[-2pt]
&\leq& 2c_3 \sup_{z\in\bR^d} \PP_z \bigl(\tau^b_{B(z, \delta)} \leq s\bigr),
\end{eqnarray*}
which by Lemma~\ref{lcon} goes to 0 as $s\to0$ (uniformly in $(t,
x, y)\in[t_1, t_2]\times D_\delta\times D_\delta$). Since $p^b(t,
x, y)$ is jointly continuous, it follows from the bounded
convergence theorem that $\E_x [h(s, t-s, X^b_s, y)]$ is jointly
continuous in $(s, t, y) \in[0, t_1/3]\times[t_1, t_2]\times
D_\delta$. On the other hand, for $(s, t, y)$ in any locally compact
subset of\vspace*{1pt} $(0, t_1/3)\times[t_1, t_2]\times D_\delta$, $\E_x [h(s,
t-s, X^b_s, y)]=\int_{\bR^d} p(s, x, z) h(s, t-s, z, y)\,dy$ is
equi-continuous in $x$. Therefore\vspace*{1pt} $\E_x [h(s, t-s, X^b_s, y)]$ is
jointly continuous in $(s, t, x, y)\in(0, t_1/3)\times[t_1,
t_2]\times D_\delta\times D_\delta$. Consequently, $k^b_D(t, x, y)$
is jointly continuous in $(s, t, y) \in[0, t_1/3]\times[t_1,
t_2]\times D_\delta$ and hence on $(0, \infty) \times D\times D$.
Since $p^b(t,x,y)$ is jointly continuous, we can now conclude that
$p^b_D(t,x,y)$ is jointly continuous on $(0, \infty) \times D\times
D$.

By Proposition~\ref{pcontb}, the last assertion of the theorem can
be proved using the argument in the last paragraph of the proof of
\cite{CZ}, Theorem 2.4. We omit the details.\vspace*{-2pt}
\end{pf}

The next result is a short time lower bound estimate for $p^b_D(t,
x, y)$ near the diagonal. The technique used in its proof is
well known. We give the proof here to demonstrate that symmetry
of the process is not needed.\vspace*{-2pt}
\begin{prop}\label{P35-55}
For any $a_1\in(0, 1)$, $a_2>0$, $a_3>0$ and $R>0$, there is a
constant
$c=c(d, \alpha, a_1, a_2, a_3, R, b)>0$ with the dependence\vadjust{\goodbreak} on
$b$ only via the rate at which $M^\alpha_{|b|}(r)$ goes to zero such
that
such that for all $x_0\in\bR^d$ and $r\in(0, R]$,
%
%
\begin{eqnarray}\label{eaa1}
p^b_{B(x_0,r)}(t, x,y)&\ge& c t^{-d /\alpha} \nonumber\\[-8pt]\\[-8pt]
&&\eqntext{\mbox{for
all } x, y \in B(x_0, a_1 r) \mbox{ and } t\in[ a_2r^\alpha,
a_3 r^\alpha].}
\end{eqnarray}
\end{prop}
\begin{pf}
Let $\kappa:= a_2/(2a_3)$ and $B_r:=B(x_0,r)$. We first show that
there is a constant $c_1\in(0, 1)$ so that~(\ref{eaa1}) holds for
all $r>0$, $x, y\in B(x_0, a_1r)$ and $t\in[\kappa
c_1r^\alpha, c_1 r^\alpha]$.

For $r>0$, $t\in[ \kappa c_1r^\alpha, c_1 r^\alpha]$, and $
x, y\in B(x_0, a_1 r)$, since $|x-y| \le2a_1 r \le2a_1 (\kappa
c_1)^{-1/\alpha} t^{-1/\alpha}$ and $t \le c_1r^\alpha\le R^\alpha
$, we have by
(\ref{e11a}),
(\ref{eqdipco}) and
(\ref{eqdipco1}),
%
%
\begin{eqnarray}\label{eqdipco22}\qquad
p^b_{B_r}(t,x,y)
&\ge&
c_2 c_1^{1+d/\alpha} t^{-d/\alpha}\nonumber\\[-8pt]\\[-8pt]
&&{} - c_3
\E_x\biggl[{\mathbf1}_{\{\tau^b_{B_r}\le
t\}}\biggl((t-\tau^b_{B_r})^{-d/\alpha}\wedge\frac
{t-\tau^b_{B_r}}{|X^b_{\tau^b_{B_r}}-y|^{d+\alpha}}\biggr)\biggr],
\nonumber
\end{eqnarray}
where the positive constants $c_i=c_i(d, \alpha, a_1, a_2, a_3, R,
b), i=2,3,$ are independent of $c_1 \in(0, 1]$. Observe that
\[
|X^b_{\tau^b_{B_r}}-y|\ge(1-a_1)r\qquad
\mbox{for } t-\tau^b_{B_r}\le t\le c_1r^\alpha,
\]
and so
%
%
\begin{equation}\label{eqdipco33}
\frac{t-\tau^b_{B_r}}{|X^b_{\tau^b_{B_r}}-y|^{d+\alpha}} \le\frac
{t-\tau^b_{B_r}}{((1-a_1)r))^{d+\alpha}}\le
\frac{c_1^{1+d/\alpha}}{(1-a_1)^{d+\alpha}} t^{-d/\alpha}.
\end{equation}
Note that if $c_1 < ((1-a_1)/2)^{\alpha}$,
by~(\ref{e11a}),
for $t\le c_1 r^\alpha$,
\begin{eqnarray*}
\PP_x \bigl(X^b_t\notin B\bigl(x, (1-a_1)r/2 \bigr)\bigr) &=& \int_{B(x,
(1-a_1)r/2)^c} p^b(t,x,y) \,dy\\
&\le& c_3 \int_{B(x, (1-a_1)r/2)^c}
\frac{t}{|x-y|^{d+\alpha}}\,dz \\
&\leq& c_{4} \frac{t}{r^\alpha}\le c_4c_1,
\end{eqnarray*}
where $c_4$ is independent of $c_1$. Now by the same argument as
in the proof of Lemma~\ref{lcon}, we have
%
%
\begin{equation}\label{eqdipco23}
\PP_x \bigl( \tau^b_{B(x,(1-a_1)r) }\le t\bigr) \leq2 c_4c_1 .
\end{equation}
Consequently, we have from~(\ref{eqdipco22})--(\ref{eqdipco23}),
\begin{eqnarray*}
p^b_{B_r}(t,x,y)&\ge& \biggl(c_2c_1^{1+d/\alpha}-c_3
\frac{c_1^{1+d/\alpha}}{(1-a_1)^{d+\alpha}} \PP_x ( \tau^b_{B_r
}\le t) \biggr) t^{-d/\alpha}\\
&\ge& \biggl(c_2 c_1^{1+d/\alpha} -c_3 \frac{c_1^{1+d/\alpha}}{(1-
a_1)^{d+\alpha}} \PP_x \bigl( \tau^b_{B(x,(1-a_1)r) }\le
t\bigr)\biggr) t^{-d/\alpha}\\
&\ge& c_1^{1+d/\alpha} \biggl(c_2-2c_4c_3
\frac{c_1}{(1-a_1)^{d+\alpha}} \biggr) t^{-d/\alpha} .
\end{eqnarray*}
Clearly we can choose $c_1< a_3 \wedge((1-a_1)/2)^{\alpha}$ small
so that $ p^b_{B_r}(t,x,y)\ge c_5 t^{-d/\alpha}$. This establishes
(\ref{eaa1}) for any $x_0 \in\bR^d$, $r>0$ and $t\in[ \kappa
c_1r^\alpha, c_1 r^\alpha]$.

Now for $r>0$ and $t\in[ a_2r^\alpha, a_3 r^\alpha]$, define
$k_0 = [a_3/c_1]+1$. Here for $a\geq1$, $[a]$ denotes the largest
integer that does not exceed $a$. Then, since $c_1 <a_3$, $t/k_0 \in
[\kappa c_1 r^\alpha, c_1r^{\alpha}]$. Using the semigroup
property~(\ref{eck}) $k_0$ times, we conclude that for all $x, y\in
B(x_0, a_1 r)$ and
$t\in[ a_2r^\alpha, a_3 r^\alpha]$,
\begin{eqnarray*}
&&p^b_{B(x_0,r)}(t, x,y)\\
&&\qquad= \int_{B(x_0,r)}\cdots\int_{B(x_0,r)} p^b_{B(x_0,r)} (t/k_0,
x,w_1)\cdots\\
&&\qquad\quad\hspace*{79.1pt}{} \times p^b_{B(x_0,r)} (t/k_0, w_{n-1},y) \,dw_1\cdots
dw_{n-1}\\
&&\qquad\geq\int_{B(x_0,a_1r)}\cdots\int_{B(x_0,a_1r)}
p^b_{B(x_0,r)} (t/k_0, x,w_1)\cdots\\
&&\qquad\quad\hspace*{94.6pt}{} \times p^b_{B(x_0,r)} (t/k_0, w_{n-1},y)
\,dw_1\cdots
dw_{n-1}\\
&&\qquad\geq c_5 (t/k_0)^{-d /\alpha} \bigl(c_5 (t/k_0)^{-d /\alpha}
|B(0,1)|(a_1 r)^d\bigr)^{k_0-1} \geq c_{6} t^{-d/\alpha}.
\end{eqnarray*}
The proof of~(\ref{eaa1}) is now complete.
\end{pf}

Using the domain monotonicity of $p^b_D$, the semigroup property
(\ref{eck}) and the L\'evy system of $X^b$, the above proposition
yields the following.
\begin{cor}\label{csp}
For every open subset $D\subset\bR^d$, $p^b_D(t,x,y)$ is strictly
positive.
\end{cor}
\begin{pf}
For $x\in D$, denote by $D(x)$ the connected component of $D$
that contains $x$. If $y\in D(x)$, using a chaining argument and
Proposition~\ref{P35-55}, we have
\[
p^b_D(t, x, y)\geq p^b_{D(x)}(t, x, y) >0.
\]
If $y\notin D(x)$, then by using the strong Markov property
and the L\'evy system~(\ref{els4xb}) of $X^b$,
\begin{eqnarray*}
&&p^b_D(t, x, y)\\
&&\qquad= \E_x \bigl[ p^b_D\bigl(t-\tau^b_{D(x)},
X^b_{\tau^b_{D(x)}}, y\bigr);
\tau^b_{D(x)} < t\bigr] \\
&&\qquad\geq \E_x \bigl[ p^b_D\bigl(t-\tau^b_{D(x)}, X^b_{\tau^b_{D(x)}}, y\bigr);
\tau^b_{D(x)} < t, X^b_{\tau^b_{D(x)}}\in D(y)\bigr]\\
&&\qquad\geq \int_0^t \int_{D(x)} p^b_{D(x)}(s, x, z)
\biggl(\int_{D(y)} J(z, w) p^b_{D(y)} (t-s, w, y) \,dw\biggr) \,dz \,ds>0.
\end{eqnarray*}
The corollary is thus proved.
\end{pf}

In the remainder of this section we assume that $D$ is a bounded
open set in $\bR^d$. The proof of the next lemma is standard; for
example, see~\cite{KS}, Lemma 6.1.
\begin{lemma}\label{dec}
There exist positive constants $C_1$ and $C_2$ depending only on
$d$, $\alpha$, $\operatorname{diam}(D)$ and $b$ with the dependence on $b$
only through the rate at which $ M^\alpha_{|b|}(r)$ goes to zero
such that
\[
p^b_D(t, x, y)\le C_1e^{-C_2t},\qquad (t, x, y)\in(1, \infty)
\times D\times D.
\]
\end{lemma}
\begin{pf}
Put $L:=\operatorname{diam} (D)$. By~(\ref{e10}), for
every $x \in D$ we have
\begin{eqnarray*}
\PP_x ( \tau^b_D \le1)
&\ge&\PP_x ( X^b_1 \in\bR^d \setminus D)
= \int_{ \bR^d \setminus D} p^b(1,x,y) \,dy\\
&\ge& c_1 \int_{ \bR^d \setminus D}
\biggl(1 \wedge\frac{1}{|x-y|^{d+\alpha}} \biggr) \,dy\\
&\ge& c_1 \int_{\{ |z| \ge L \}}
\biggl(1 \wedge\frac{1}{|z|^{d+\alpha}} \biggr) \,dz > 0.
\end{eqnarray*}
Thus
\[
\sup_{x \in D}\int_Dp^b_D(1, x, y)\,dy =
\sup_{x \in D} \PP_x ( \tau^b_D > 1) < 1.
\]
The Markov property of $X^b$ then implies that there
exist positive constants $c_2$ and $c_3$
such that
\[
\int_Dp^b_D(t, x, y)\,dy \le c_2e^{-c_3t}
\qquad\mbox{for } (t, x)\in(0, \infty)\times D.
\]
It follows from
(\ref{e10}) that there exists $c_4>0$ such that
$p^b_D(1, x, y)\le p^b(1, x, y)\le c_4$
for every $(x, y)\in D\times D$.
Thus for any $(t, x, y)\in(1, \infty)\times D\times D$, we have
\begin{eqnarray*}
p^b_D(t, x, y)&=&\int_Dp^b_D(t-1, x, z)p^b_D(1, z, y)\,dz
\\&\le& c_4\int_Dp^b_D(t-1, x, z)\,dz
\le c_2c_4e^{-c_3(t-1)}.
\end{eqnarray*}
\upqed
\end{pf}

Combining the result above with~(\ref{e10}) we know that there
exists a positive constant $c_1 =c_1(d, \alpha, \operatorname{diam}(D), b)$
with the dependence on $b$ only through the rate at which $
M^\alpha_{|b|}(r)$ goes to zero such that for any $(t, x, y)\in(0,
\infty)\times D\times D$,
%
%
\begin{equation}\label{est61}
p^b_D(t, x, y) \le c_1
\biggl(t^{-{d}/{\alpha}}\wedge\frac{t}{|x-y|^{d+\alpha}}\biggr).
\end{equation}
Therefore the Green function
$
G^b_D(x, y)=\int^{\infty}_0 p^b_D(t, x, y)\,dt
$
is finite and continuous off the diagonal of $D\times D$ and
%
%
\begin{equation}\label{Gbd}
G^b_D(x, y) \le c_2 \frac1{|x-y|^{d-\alpha}}
\end{equation}
for some positive constant $c_2=c_2(d, \alpha, \operatorname{diam}(D), b)$
with the dependence on $b$ only through the rate at which $
M^\alpha_{|b|}(r)$ goes to zero.

\section{Uniform estimates on Green functions}

Let
\[
g_D(x,y):=\frac{1} {|x-y|^{d-\alpha}} \biggl(1\wedge\frac{
\delta_D(x) \delta_D(y)}{ |x-y|^{2}}\biggr)^{\alpha/2}.
\]
The following lemma is needed in deriving sharp bounds on the Green
function $G^b_U$ when $U$ is some small $C^{1, 1}$ open set. It can
be regarded as a new type of $3G$ estimates.
\begin{lemma}\label{l41}
There exists a positive constant
$C_3 =C_3(d, \alpha)$ such that for all $x, y, z\in D$,
%
%
\begin{equation} \label{e41}
g_D(x,z) \frac{ g_D(z, y)}{|z-y| \wedge\delta_D(z)} \leq C_3
g_D(x, y) \biggl( \frac{1}{|x-z|^{d+1-\alpha}}+
\frac1{|z-y|^{d+1-\alpha}}\biggr)\hspace*{-28pt}
\end{equation}
and
%
%
\begin{eqnarray}\label{e42}
&&\frac{g_D(x,z)}{|x-z|\wedge\delta_D (x)} \frac{g_D(z, y)}{|z-y|
\wedge\delta_D(z)}\nonumber\\[-8pt]\\[-8pt]
&&\qquad \le C_3 \frac{g_D(x,y)}{|x-y|\wedge
\delta_D(x)} \biggl( \frac{1}{|x-z|^{d+1-\alpha}}+
\frac1{|z-y|^{d+1-\alpha}}\biggr). \nonumber
\end{eqnarray}
\end{lemma}
\begin{pf}
Put $r(x,y)= \delta_D(x) +\delta_D(y)+|x-y|$. Note that for $a,
b>0$,
%
%
\begin{equation}\label{e43}
\frac{ab}{a+b} \leq a\wedge b \leq2\frac{ab}{a+b}.
\end{equation}
Moreover for $x, y\in D$, since
\[
\delta_D(x)^2 \leq\delta_D(x) \bigl(\delta_D(y)+|x-y|\bigr)
\leq\delta_D(x)\delta_D(y)+\delta_D(x)^2/2+|x-y|^2/2,
\]
one has
\[
\delta_D(x)^2 \leq2 \delta_D(x)\delta_D(y)+|x-y|^2.
\]
It follows from these observations that
%
%
\begin{equation}\label{final3}
\frac{ \delta_D(x) \delta_D(y)}{(r(x,y))^2} \le\biggl(1\wedge
\frac{ \delta_D(x) \delta_D(y)}{ |x-y|^{2}}\biggr) \le24
\frac{ \delta_D(x) \delta_D(y)}{(r(x,y))^2}.
\end{equation}
Consequently, we have
%
%
\begin{equation}\label{egg2}
g_D(x, y) \asymp\frac{1} {|x-y|^{d-\alpha}} \frac{
\delta_D(x)^{\alpha/2} \delta_D(y)^{\alpha/2}}{(r(x,y))^{\alpha}}.
\end{equation}
Now
%
%
\begin{eqnarray}\label{e44}\quad
&& g_D(x,z) \frac{ g_D(z, y)}{|z-y| \wedge\delta_D(z)} \nonumber\\
&&\qquad\asymp  g_D(x, y) \frac{|z-y|+\delta_D(z)}{|z-y| \delta_D(z)}
\frac{\delta_D(z)^\alpha r(x, y)^\alpha}{r(x, z)^\alpha
r(z, y)^\alpha} \biggl( \frac{|x-y|}{|x-z| \cdot|z-y|}
\biggr)^{d-\alpha}
\\
&&\qquad\le g_D(x, y) \frac{r(y, z)}{|z-y| }
\frac{\delta_D(z)^{\alpha-1} r(x, y)^\alpha}{r(x, z)^\alpha
r(z, y)^\alpha} \biggl( \frac{|x-y|}{|x-z| \cdot
|z-y|}\biggr)^{d-\alpha}
\nonumber\\
&&\qquad= g_D(x, y) \frac{r(x, y)}{|z-y| r(x, z)}
\biggl(\frac{\delta_D(z) r(x, y)}{r(x, z) r(z,
y)}\biggr)^{\alpha-1} \biggl( \frac{|x-y|}{|x-z| \cdot
|z-y|}\biggr)^{d-\alpha}. \nonumber
\end{eqnarray}
Since $r(x,y) \le r(x,z) + r(z,y)$,
\[
\frac{\delta_D(z) r(x, y)}{r(x, z) r(z, y)} \leq\frac
{\delta_D(z)}{r(x, z)}+\frac{\delta_D(z)}{r(z, y)} \leq2.
\]
On the other hand, since $\delta_D(y)\leq\delta_D(x) +|x-y|$,
\begin{eqnarray*}
\frac{r(x, y)}{|z-y| r(x, z)} &\leq& 2
\frac{|x-y|+\delta_D(x)}{|z-y| r(x, z)}
\leq2 \frac{|x-z|+(|z-y|+\delta_D(x))}{|z-y| r(x, z)} \\
&\leq& \frac2{r(x, z)} + \frac{2}{|z-y|}
\leq\frac{2}{|x-z|}+\frac2{|z-y|}.
\end{eqnarray*}
Hence we deduce from~(\ref{e44}) that
\begin{eqnarray*}
&& g_D(x,z) \frac{ g_D(z, y)}{|z-y| \wedge\delta_D(z)}\\
&&\qquad\leq
2^\alpha g_D(x, y) \biggl(\frac{1}{|x-z|}+\frac1{|z-y|}\biggr)
\biggl( \frac{|x-y|}{|x-z| \cdot|z-y|}\biggr)^{d-\alpha}\\
&&\qquad\leq
c_1 g_D(x, y) \biggl(\frac{1}{|x-z|}+\frac1{|z-y|}\biggr)
\biggl( \frac{1}{|x-z|^{d-\alpha}}+\frac1{|z-y|^{d-\alpha}}
\biggr)\\
&&\qquad\leq c_2 g_D(x, y) \biggl(
\frac{1}{|x-z|^{d+1-\alpha}}+\frac1{|z-y|^{d+1-\alpha}}\biggr),
\end{eqnarray*}
where $c_1$ and $c_2$ are positive constants depending only on $d$
and $\alpha$. This proves~(\ref{e41}).\vadjust{\goodbreak}

Now we show that~(\ref{e42}) holds. Note that by~(\ref{egg2}),
%
%
\begin{eqnarray}\label{e45}
&& \frac{g_D(x,z)}{|x-z|\wedge\delta_D (x)}
\frac{g_D(z, y)}{|z-y| \wedge\delta_D(z)} \nonumber\\
&&\qquad\asymp \frac{\delta_D(x)^{\alpha/2}\delta_D(y)^{\alpha/2}}
{|x-z|^{d+1-\alpha}|z-y|^{d+1-\alpha}}
\frac{|x-z|\cdot|z-y|}{(|x-z|\wedge\delta_D(x))(|z-y|\wedge\delta_D(z))}
\nonumber\\[-8pt]\\[-8pt]
&&\qquad\quad{} \times\frac{\delta_D(z)^\alpha}{r(x, z)^\alpha r(z,
y)^\alpha} \nonumber\\
&&\qquad\asymp \frac{g_D(x, y)}{|x-y|\wedge\delta_D(x)} \cdot
\frac{|x-y|^{d+1-\alpha}}{|x-z|^{d+1-\alpha}|z-y|^{d+1-\alpha}}
\cdot I, \nonumber
\end{eqnarray}
where
\[
I:= \frac{|x-y|\wedge\delta_D(x)}{|x-y|} \cdot
\frac{|x-z|\cdot|z-y|}{(|x-z|\wedge\delta_D(x))(|z-y|\wedge\delta_D(z))}
\frac{\delta_D(z)^\alpha r(x, y)^\alpha}{r(x, z)^\alpha
r(z, y)^\alpha}.
\]
It follows from~(\ref{e43}) and the fact that $|x-z|+\delta_D(z)
\asymp r(x, z)$ that
\begin{eqnarray*}
I &\asymp& \frac{|x-y| \delta_D(x)}{|x-y| (|x-y|+\delta_D(x))}\\
&&{}\times
\frac{|x-z|\cdot|z-y| (|x-z|+\delta_D(x))(|z-y|+\delta_D(z))}
{(|x-z| \delta_D(x))(|z-y| \delta_D(z))}
\frac{\delta_D(z)^\alpha r(x, y)^\alpha}{r(x, z)^\alpha
r(z, y)^\alpha}
\\
&\asymp& \frac{\delta_D(z)^{\alpha-1} r(x, y)^{\alpha-1}}{r(x,
z)^{\alpha-1}
r(z, y)^{\alpha-1}}
\le
\delta_D(z)^{\alpha-1} \biggl( \frac1{r(x, z)^{\alpha-1}}
+\frac1{r(y, z)^{\alpha-1}}\biggr)\\
&\leq&2.
\end{eqnarray*}
Inequality~(\ref{e42}) now follows from~(\ref{e45}).
\end{pf}

Recall that $G_D$ is the Green function of $X^D$. It is known that
%
%
\begin{equation}\label{eg1}
|\nabla_z G_D(z, y)| \le\frac{d}{|z-y| \wedge\delta_D(z)} G_D(z,y);
\end{equation}
see~\cite{BKN}, Corollary 3.3.
Recall also that $b$ is an $\bR^d$-valued function on
$\bR^d$ such that $|b|\in\bK_{d, \alpha-1}$.
\begin{prop}\label{tdh1}
If $D$ is a bounded open set, and ${\mathbf1}_D b$ has compact support
in $D$, then $G^b_D$ satisfies
%
%
\begin{equation}\label{Du2}
G^b_{D} (x, y)
=G_D(x, y)+
\int_{D}G^b_{D}(x, z) b(z)\cdot\nabla_z G_D(z, y)\,dz.
\end{equation}
\end{prop}
\begin{pf}
Recall that by\vspace*{1pt}
Theorem~\ref{T25}, for every $f \in C^\infty_c(\bR^d
)$, $M^f_t:=f(X^b_t)-f(X^b_0)-\int^t_0 \sL^b f(X^b_s) \,ds$ is a
martingale with respect to $\PP_x$. Since\vspace*{-1pt} ${\mathbf1}_D b$ has compact support
in $D$, in view of
(\ref{Gbd}),~(\ref{eg1}) and the fact that $|b|\in\bK_{d,
\alpha-1}$,
$M^f_{t\wedge
\tau_D}$ is a uniformly integrable martingale.\vadjust{\goodbreak}

Define $D_j:=\{x\in D\dvtx \operatorname{dist}(x, D^c)> 1/j\}$. Let $\phi\in
C^\infty_c(\bR^d)$ with $\phi\geq1$, $\operatorname{supp} [\phi]\subset
B(0, 1)$ and $\int_{\bR^d} \phi(x) \,dx=1$. For any $\psi\in
C_c(D)$, define $f=G_D\psi$ and $f_n:=\phi_n * f$, where $\phi_n
(x):=n^d \phi(nx)$. Clearly $f_n\in C^\infty_c(\bR^d)$ and $f_n$
converges uniformly to $f=G_D\psi$. Fix $j\geq1$. Since
$\E_x[M^{f_n}_0] =\E_x[M^{f_n}_{\tau_{D_j}}]$, and for every $y\in
D_j$ and sufficiently large $n$,
\[
\phi_n *( \Delta^{\alpha/2}f) (y)
= \int_{B(0, 1/n)}
\phi_n(z)\Delta^{\alpha/2}(G_D\psi)(y-z)\,dz,
\]
we have, by Dynkin's formula, that for sufficiently large $n$,
\begin{eqnarray*}
&&\E_x[f_n(X^b_{\tau_{D_j}})] -f_n(x)\\
&&\qquad= \int_{D_j} G^b_{D_j}(x,y)
\bigl(\Delta^{\alpha/2} f_n (y)+ b(y) \cdot\nabla f_n (y)\bigr)\,dy
\\
&&\qquad= \int_{D_j} G^b_{D_j}(x,y)
\bigl(\phi_n *( \Delta^{\alpha/2}f) (y)+ b(y) \cdot\phi_n*(\nabla
f) (y) \bigr)\,dy\\
&&\qquad= \int_{D_j} G^b_{D_j}(x,y)
\bigl(-\phi_n *\psi(y)+ b(y) \cdot\phi_n*(\nabla
(G_D\psi) (y) )\bigr)\,dy.
\end{eqnarray*}
Taking $n\to\infty$, we get, by~(\ref{Gbd}),~(\ref{eg1}) and the
fact that $|b|\in\bK_{d, \alpha-1}$,
%
%
\begin{equation}\label{e49}
\E_x[f (X^b_{\tau_{D_j}})] -f (x) = \int_{D}
G^b_{D_j}(x,y) \bigl( -\psi(y)+ b(y) \cdot\nabla(G_D\psi) (y)
\bigr)\,dy.\hspace*{-32pt}
\end{equation}
Now using the fact that ${\mathbf1}_D b$ has compact support in $D$, taking
$j\to\infty$, we have by
(\ref{Gbd}),~(\ref{eg1}) and the fact that $|b|\in\bK_{d,
\alpha-1}$,
\[
-f (x)= \int_{D } G^b_{D }(x,y)
\bigl(-\psi(y)+ b(y) \cdot\nabla(G_D\psi) (y) \bigr)\,dy.
\]
Hence we have
\[
-G_D \psi(x)= - G^b_D \psi+ G^b_D (b\cdot\nabla G_D \psi).
\]
This shows that for each $x\in D$,~(\ref{Du2}) holds for a.e. $y\in
D$.
Since $G^b_D$ is continuous off the diagonal of $D\times D$, we get
that~(\ref{Du2}) holds for all $x, y\in D$.
\end{pf}

We will derive two-sided estimates on the Green function of $X^b$ on
certain nice open sets when the diameter of such open sets are less
than or equal to some constant depending on $b$ only through the
rate at which $M^\alpha_{|b|}(r)$ goes to zero.
\begin{prop}
\label{tufgfnest1} There exists a positive constant $r_* =r_*(d,
\alpha, b)$ with the dependence on $b$ only via the rate at which
$M^\alpha_{|b|}(r)$ goes to zero such that for any ball $B=B(x_0,
r)$ of radius $r\le r_*$ and any $n\ge1$,
\[
{2}^{-1}G_B(x, y)\le G^{b_n}_B(x, y)\le2 G_B(x, y),\qquad x, y\in
B,\vadjust{\goodbreak}
\]
where
%
%
\begin{equation}\label{ebB}
b_n(x)=b(x){\mathbf1}_{B^c}(x)+b(x){\mathbf1}_{K_n}(x),\qquad
x\in\bR^d,
\end{equation}
with $K_n$ being an increasing sequence of compact subsets of $B$
such that $\bigcup_{n}K_n=B$.
\end{prop}
\begin{pf}
It is well known that there exists a constant
$c_1 =c_1(d, \alpha)>1$ such that
%
%
\begin{eqnarray}\label{egege}
&&c^{-1}_1\frac{1} {|x-y|^{d-\alpha}} \biggl(1\wedge\frac{ \delta_B(x)
\delta_B(y)}{ |x-y|^{2}}\biggr)^{\alpha/2}\nonumber\\[-8pt]\\[-8pt]
&&\qquad\le G_B(x, y) \le c_1\frac{1}
{|x-y|^{d-\alpha}} \biggl(1\wedge\frac{ \delta_B(x) \delta_B(y)}{
|x-y|^{2}}\biggr)^{\alpha/2}.\nonumber
\end{eqnarray}
Define $\wt{I}^n_k(x,y)$ recursively for $n\ge1$, $k \ge0$ and
$(x , y) \in B \times B$ by
\begin{eqnarray*}
\wt{I}^n_0(x,y)&:=&G_B(x, y), \\
\wt{I}^n_{k+1}(x,y)&:=& \int_{B} \wt{I}^n_k(x,z)
b_n(z)\cdot\nabla_z G_B( z, y)\,dz.
\end{eqnarray*}
Iterating~(\ref{Du2}) gives that for each $m\geq2$ and for every $(x
, y)\in B \times B$,
%
%
\begin{equation}\label{iter}
G^{b_n}_B( x, y)= \sum_{k=0}^{m} \wt{I}^n_k(x,y)
+ \int_B G^{b_n}_B(x, z) b_n(z) \cdot
\nabla_z \wt I^n_m (z, y) \,dz.
\end{equation}
Using induction, Lemma~\ref{l41},~(\ref{eg1}) with $D=B$ and
(\ref{egege}), we see that there exists a positive constant $c_2$
(in fact, one can take $c_2=2dC_3c_1^3$ where $C_3$ is the constant
in Lemma~\ref{l41}) depending only on $d$ and $\alpha$ such that
for $n, k\ge1$ and $(x , y)\in B \times B$,
%
%
\begin{equation}\label{inequalJ}
|\wt{I}^n_k(x,y)| \le
c_2 G_B(x,y) \bigl(c_2 M^\alpha_{|b|} (2r) \bigr)^k
\end{equation}
and
%
%
\begin{equation}\label{e412}
| \nabla_x \wt I^n_k (x, y)|\leq c_2 \frac{G_B(x, y)}{|x-y|\wedge
\delta_B(x)} \bigl(c_2 M^\alpha_{|b|} (2r) \bigr)^k.
\end{equation}
There exists an $\wh r_1>0 $ depending on $b$ only via the rate at
which $M^\alpha_{|b|}(r)$ goes to zero such that
%
%
\begin{equation}\label{N1}
c_2 M^\alpha_{|b|} (r) <
\frac{1}{1+2c_2 }
\qquad\mbox{for every } 0<r \le
\wh r_1.
\end{equation}
Equations~(\ref{Gbd}) and~(\ref{e412}),~(\ref{N1}) imply that if $r\le
\wh r_1/2$, then for $n\ge1$ and $(x , y)\in B \times B$,
\begin{eqnarray*}
&& \biggl|\int_B G^{b_n}_B(x, z) b_n(z) \cdot\nabla_z \wt
I^n_m (z, y) \,dz \biggr|\\
&&\qquad\le c_2\biggl(\int_B
G^{b_n}_{B }(x, z)
|b_n(z)| \frac{G_B(z, y)}{|z-y|\wedge
\delta_B(z)} \,dz \biggr) \bigl(c_2 M^\alpha_{|b|} (2r) \bigr)^m\\
&&\qquad\le c_3\biggl(\int_B \frac{1}{|x-z|^{d-\alpha}} \frac{G_B(z,
y)}{|z-y| }|b(z)|
\,dz \biggr)\biggl(\frac{1}{1+2c_2 } \biggr)^m\\
&&\qquad\le c_4\biggl(\int_B \frac{1}{|x-z|^{d+1-\alpha}}
\frac{|b(z)|}{|z-y|^{d+1-\alpha}}
\,dz \biggr)\biggl(\frac{1}{1+2c_2 } \biggr)^m\\
&&\qquad\le c_5 (1+2c_2)^{-m}
|x-y|^{-(d+1-\alpha)} \int_B \biggl( \frac
{|b(z)|}{|x-z|^{d+1-\alpha}}
+\frac{|b|(z)}{|y-z|^{d+1-\alpha}} \biggr)\,dz \\
&&\qquad\le c_6 (1+2c_2 )^{-(m+1)}
|x-y|^{-(d+1-\alpha)},
\end{eqnarray*}
which goes to zero as $m \to\infty$. In the second inequality, we
have used the fact that $b_n$ is compactly supported in $B$. Thus,
by~(\ref{iter}), $G^{b_n}_B( x, y)= \sum_{k=0}^{\infty} \wt
{I}^n_k(x,y)$. Moreover, by~(\ref{inequalJ}),
\[
\sum_{k=1}^{\infty} |\wt{I}^n_k(x,y)| \le c_2 G_B( x, y)
\sum_{k=1}^{\infty}(1+2c_2)^{-k} \le G_B( x, y)/2.
\]
It follows that for any $x_0 \in\bR^d $ and $B=B(x_0, r)$ of
radius $r\le\wh r_1/2$,
\[
G_B(x, y)/2\leq G^{b_n}_B(x, y) \leq3G_B(x, y)/2 \qquad\mbox{for
all } n\ge1 \mbox{ and } x, y\in B.
\]
This proves the theorem.
\end{pf}

For any bounded $C^{1, 1}$ open set $D$ with characteristic $(R_0,
\Lambda_0)$, it is well known (see, e.g.,~\cite{So}, Lemma 2.2) that there exists $L=L(R_0, \Lambda_0, d)>0$ such that for
every $z \in\partial D$ and $r \le R_0$, one can find a $C^{1,1}$
open set $U_{(z, r)}$ with characteristic $(rR_0/L, \Lambda_0L/r)$
such that $D \cap B(z, r/2) \subset U_{(z, r)} \subset D \cap B(z,
r) $. For the remainder of this paper, given a bounded $C^{1, 1}$
open set $D$, $U_{(z, r)}$ always refers to the $C^{1, 1}$ open set
above.

For $U_{(z, r)}$, we also have a result similar to Proposition
\ref{tufgfnest1}.
\begin{prop}
\label{tufgfnest2} For every $C^{1,1}$ open set $D$ with the
characteristic $(R_0, \Lambda_0)$, there exists $r_0 =r_0(d, \alpha,
R_0, \Lambda_0, b) \in(0, (R_0 \wedge1)/8]$ with the dependence
on $b$ only via the rate at which $M^\alpha_{|b|}(r)$ goes to zero
such that for all $n\ge1$, $z\in\partial D$ and $r\le
r_0$, we have
%
%
\begin{equation}\label{eempaok}\qquad
2^{-1} G_{U_{(z, r)}}(x,y) \le G^{b_n}_{U_{(z, r)}}(x,y) \le2
G_{U_{(z, r)}}(x,y),\qquad x, y\in U_{(z, r)},
\end{equation}
where
%
%
\begin{equation}\label{ebD}
b_n(x)=b(x){\mathbf1}_{U_{(z, r)}^c}(x)+b(x){\mathbf1}_{K_n}(x),\qquad
x\in\bR^d,
\end{equation}
with $K_n$ being an increasing sequence of compact subsets of
$U_{(z, r)}$ such that $\bigcup_{n}K_n=U_{(z, r)}$.
\end{prop}
\begin{pf}
It is well known (see~\cite{J}, e.g.) that, for any bounded
$C^{1, 1}$ open set~$U$,
there exists $c_1
=c_1(R_0, \Lambda_0, \operatorname{diam}(U))>1$ such that
%
%
\begin{eqnarray}\label{e411}
&&c^{-1}_1\frac{1} {|x-y|^{d-\alpha}} \biggl(1\wedge\frac{ \delta_U(x)
\delta_U(y)}{ |x-y|^{2}}\biggr)\nonumber\\[-8pt]\\[-8pt]
&&\qquad\le G_U(x, y) \le
c_1\frac{1}
{|x-y|^{d-\alpha}} \biggl(1\wedge\frac{ \delta_U(x) \delta_U(y)}{
|x-y|^{2}}\biggr).\nonumber
\end{eqnarray}
It follows from this, the fact that $r^{-1}U_{(z,r)}$ is a $C^{1,1}$
open set with characteristic $(R_0/L, \Lambda_0L)$ and scaling that,
for any bounded $C^{1, 1}$ open set $D$ with characteristics $(R_0,
\Lambda_0)$, there exists a constant $c_2=c_2(R_0, \Lambda_0, d)>1$
such that for all $z \in\partial D$, $r \le R_0$ and $x, y\in
U_{(z, r)}$,
\begin{eqnarray*}
&&c^{-1}_2\frac{1} {|x-y|^{d-\alpha}} \biggl(1\wedge\frac{
\delta_{U_{(z, r)}}(x) \delta_{U_{(z, r)}}(y)}{ |x-y|^{2}}\biggr)\\
&&\qquad\le
G_{U_{(z, r)}}(x, y) \le c_2\frac{1} {|x-y|^{d-\alpha}}
\biggl(1\wedge\frac{ \delta_{U_{(z, r)}}(x) \delta_{U_{(z, r)}}(y)}{
|x-y|^{2}}\biggr).
\end{eqnarray*}
Now we can repeat the argument of Theorem~\ref{tufgfnest1} to
complete the proof.
\end{pf}

Now we are going to extend Propositions
\ref{tufgfnest1} and~\ref{tufgfnest2} to $G^b_B$ and $G^b_{U(z, r)}$.
For the remainder of this section, we let $U$ be either a ball
$B=B(x_0, r)$ with $r \le r_*$ where $r_*$ is the constant in
Proposition~\ref{tufgfnest1} or $U(z, r)$ [for a $C^{1,1}$ open set
$D$ with the characteristic $(R_0, \Lambda_0)$] with $r\le r_0$
where $r_0$ is the constant in Proposition~\ref{tufgfnest2}. We
also let $b_n$ be defined by either~(\ref{ebB}) or~(\ref{ebD}), and
we will take care of the two cases simultaneously.

By~\cite{BJ}, Lemma 13, and its proof, there exists a constant $C_4>0$
such that
\[
\int_{\bR^d} \int_0^t p(t-s, x,z)|b(z)| |\nabla_z p(s, z, y)|\,ds\,dz
\le C_4 p(t,x,y) \bN_b(t)
\]
and so
%
%
\begin{equation}\label{eN0}\qquad
\int_{\bR^d} \int_0^t p(t-s, x,z)|b_n(z)| |\nabla_z p(s, z, y)|\,ds\,dz
\le C_4 p(t,x,y) \bN_b(t),
\end{equation}
where
\[
\bN_b(t):= \sup_{w \in\bR^d} \int_{\bR^d}\int_0^t
|b(z)|\bigl(
|w-z|^{-d-1}\wedge s^{-(d+1)/\alpha}\bigr) \,ds\,dz,
\]
which is finite and goes to zero as $t \to0$ by~\cite{BJ}, Corollary 12.
We remark that the constant $C_4$ here
is independent of $t$ and is not the same constant $C_4$
from~\cite{BJ}, Lemma 13.
Moreover,
%
%
\begin{eqnarray}\label{eN1}\qquad
&&\int_{\bR^d} \int_0^t p(t-s, x,z)|b(z)-b_n(z)| |\nabla_z p(s, z,
y)|\,ds\,dz \nonumber\\
&&\qquad\le C_4
p(t,x,y) \bN_{b-b_n}(t)\nonumber\\[-8pt]\\[-8pt]
&&\qquad= C_4 p(t,x,y) \nonumber\\
&&\qquad\quad{}\times\sup_{w \in\bR^d}
\int_{U \setminus K_n}\int_0^t |b(z)| \bigl(
|w-z|^{-d-1} \wedge s^{-(d+1)/\alpha}\bigr) \,ds\,dz.\nonumber
\end{eqnarray}
Now, by~\cite{BJ}, (27),
%
%
\begin{equation}\label{inequalIk}
|p^b_k(t,x,y)| \vee|p^{b_n}_k(t,x,y)| \le (C_4\bN_{b}(t))^k
p(t,x,y).
\end{equation}
Choose $T_1>0$ small so that
%
%
\begin{equation}\label{eN2}
C_4\bN_{b}(t) <\tfrac12,\qquad t \le T_1.
\end{equation}
We will fix this constant $T_1$ until the end of this section.
\begin{lemma}\label{l33k}
For all $k \ge1$ and $(t,x,y) \in(0, T_1] \times\bR^d \times
\bR^d$,
\begin{eqnarray*}
&&| p^{b_n}_k(t,x,y)-p^b_k(t,x,y)| \\
&&\qquad\le k C_4 2^{-(k-1)} p(t,x,y) \\
&&\qquad\quad{}\times\sup_{w
\in\bR^d} \int_{U \setminus K_n}\int_0^t |b(z)|\bigl(
|w-z|^{-d-1} \wedge s^{-(d+1)/\alpha}\bigr) \,ds\,dz.
\end{eqnarray*}
\end{lemma}
\begin{pf}
We prove the lemma by induction. For $k=1$, we have
\begin{eqnarray*}
&&| p^{b_n}_{1}(t,x,y)-p^b_{1}(t,x,y)| \\
&&\qquad\le\int^{t}_0 \int_{\bR^d}
p(s,x,z) | \nabla_z p(t-s, z, y)| |b-b_n|(z)\,dz\,ds.
\end{eqnarray*}
Thus by~(\ref{eN1}), the lemma is true for $k=1$.

Next we assume that the lemma holds for $k\geq1$.
We will show that the lemma hods for $k+1$.
Let
\[
I(n,t,x,y):= \int^{t}_0 \int_{\bR^d}| p^b_k(s,x,z)|| \nabla_z
p(t-s, z, y)| |b-b_n|(z)\,dz\,ds
\]
and
\begin{eqnarray*}
&&\mathit{II}(n,t,x,y)\\
&&\qquad:=\int^{t}_0 \int_{\bR^d}|p^{b_n}_k(s,x,z)-
p^b_k(s,x,z)| |\nabla_z p(t-s, z, y)| |b_n(z)|\,dz\,ds.
\end{eqnarray*}
Then we have
\[
| p^{b_n}_{k+1}(t,x,y)-p^b_{k+1}(t,x,y)| \le
I(n,t,x,y) + \mathit{II}(n,t,x,y).
\]
By~(\ref{eN1})--(\ref{eN2}),
%
%
\begin{eqnarray}\label{enn1}
&&I(n, t, x, y)\nonumber\\
&&\qquad\le (C_4\bN_{b}(t))^k
\int_{\bR^d} \int_0^t p(t-s, x,z)|b(z)-b_n(z)| |\nabla_z p(s, z,
y)|\,ds\,dz \\
&&\qquad= C_4 2^{-k} p(t,x,y) \sup_{w \in\bR^d} \int_{U \setminus
K_n}\int_0^t |b(z)|\bigl(
|w-z|^{-d-1} \wedge s^{-(d+1)/\alpha}\bigr)
\,ds\,dz. \hspace*{-14pt}\nonumber
\end{eqnarray}
On the other hand, by the induction assumption,~(\ref{eN0}) and
(\ref{eN2}),
%
%
\begin{eqnarray}\label{enn2}\quad
&&\mathit{II}(n, t, x, y) \nonumber\\
&&\qquad\le kC_4 2^{-(k-1)} \biggl(\sup_{w \in\bR^d} \int_{U
\setminus K_n}
\int_0^t |b(z)|\bigl(
|w-z|^{-d-1} \wedge s^{-(d+1)/\alpha}\bigr) \,ds\,dz \biggr) \nonumber\\
&&\qquad\quad{}\times\int_{\bR^d} \int_0^t p(s,x,z) |\nabla_z p(t-s, z,
y)| |b_n(z)|\,dz\,ds
\nonumber\\
&&\qquad\le kC_4 2^{-(k-1)} (C_4\bN_{b}(t))p(t,x,y) \\
&&\qquad\quad{}\times\sup_{w \in\bR^d}
\int_{U \setminus K_n}
\int_0^t |b(z)|\bigl(
|w-z|^{-d-1} \wedge s^{-(d+1)/\alpha}\bigr) \,ds\,dz \nonumber\\
&&\qquad\le kC_4
2^{-k} p(t,x,y) \nonumber\\
&&\qquad\quad{}\times\sup_{w \in\bR^d} \int_{U \setminus K_n}\int_0^t |b(z)|
\bigl(
|w-z|^{-d-1} \wedge s^{-(d+1)/\alpha}\bigr) \,ds\,dz. \nonumber
\end{eqnarray}
Combining~(\ref{enn1}) and~(\ref{enn2}),
we see that the lemma holds for $k+1$, and thus by induction,
the lemma holds for every $k\geq1$.
\end{pf}
\begin{theorem}\label{t34}
$p^{b_n}(t, x, y)$ converges uniformly to $p^{b}(t, x, y)$ on any
$[t_0$, $T] \times\bR^d \times\bR^d$, where $0<t_0<T<\infty$.
Moreover,
%
%
\begin{equation}\label{egcon}
\lim_{n\to\infty} G_U^{b_n}f = G_U^bf \qquad\mbox{for every } f
\in C_b(\overline U) .
\end{equation}
\end{theorem}
\begin{pf}
Without of loss of generality, we may assume that
$0<t_0\leq T_1/2$, where $T_1$ is the constant in~(\ref{eN2}).
We first consider the case $(t,x,y) \in[t_0, T_1] \times\bR^d
\times\bR^d$. By Theorem~\ref{T11}(i) and Lemma~\ref{l33k},
\begin{eqnarray*}
&& \sup_{(t,x,y) \in[t_0, T_1] \times\bR^d \times\bR^d}
|p^{b}(t, x, y)-p^{b_n}(t, x, y)| \\
&&\qquad\le \sup_{(t,x,y) \in[t_0, T_1] \times\bR^d \times\bR^d}
\sum_{k=1}^\infty
| p^{b_n}_k(t,x,y)-p^b_k(t,x,y)| \\
&&\qquad\le C_4 \sup_{(t,x,y) \in[t_0, T_1] \times\bR^d \times\bR^d}
\sum_{k=1}^\infty k2^{-(k-1)} p(t,x,y) \\
&&\qquad\quad{} \times\sup_{w \in\bR^d} \int_{U
\setminus K_n}\int_0^t |b(z)|\bigl( |w-z|^{-d-1} \wedge
s^{-(d+1)/\alpha}
\bigr) \,ds\,dz \\
&&\qquad\le cC_4t_0^{-d/\alpha} \sup_{w \in\bR^d} \int_{U \setminus
K_n}\int_0^{T_1} |b(z)|\bigl(|w-z|^{-d-1} \wedge
s^{-(d+1)/\alpha}\bigr) \,ds\,dz,
\end{eqnarray*}
which goes to zero as $n \to\infty$.

If $(t,x,y) \in(T_1, 3T_1/2] \times\bR^d \times\bR^d$, using the
semigroup property
(\ref{esemi})
with
$t_1=T_1/2$,
\begin{eqnarray*}
\hspace*{-4pt}&& \mathop{\sup_{(t,x,y) \in}}_{(T_1, 3T_1/2] \times\bR^d \times\bR^d}
|p^{b}(t, x, y)-p^{b_n}(t, x, y)|\\
\hspace*{-4pt}&&\qquad\le
\mathop{\sup_{(t,x,y) \in}}_{(T_1, 3T_1/2] \times\bR^d \times\bR^d}
\biggl|\int_{\bR^d} p^{b}(t_1, x, z)p^{b}(t-t_1, z, y)\,dz \\
\hspace*{-4pt}&&\hspace*{72pt}\qquad\quad{} -
\int_{\bR^d} p^{b_n}(t_1, x, z)p^{b_n}(t-t_1, z, y)\,dz \biggr|\\
\hspace*{-4pt}&&\qquad\le
\mathop{\sup_{(t,x,y) \in}}_{(T_1, 3T_1/2] \times\bR^d \times\bR^d}
\int_{\bR^d} p^{b}(t_1, x, z)|p^{b}(t-t_1, z, y)-p^{b_n}(t-t_1, z,
y)|\,dz \\
\hspace*{-4pt}&&\qquad\quad{} +
\mathop{\sup_{(t,x,y) \in}}_{(T_1, 3T_1/2] \times\bR^d \times\bR^d}
\int_{\bR^d} |p^{b_n}(t_1, x, z)-p^{b}(t_1, x, z)| p^{b_n}(t-t_1, z,
y)\,dz,
\end{eqnarray*}
which is, by~(\ref{e10}), less than or equal to
$c_1 t_1^{-d/\alpha}$
times
\begin{eqnarray*}
&&\sup_{(t,y) \in(T_1, 3T_1/2] \times\bR^d}
\int_{\bR^d}|p^{b}(t-t_1, z, y)-p^{b_n}(t-t_1, z, y)|\,dz \\
&&\qquad{}+
\sup_{x \in\bR^d }
\int_{\bR^d} |p^{b_n}(t_1, x, z)-p^{b}(t_1, x, z)| \,dz \,ds.
\end{eqnarray*}
Now, by the first case, we see that the above goes to zero as $n \to
\infty$. Iterating the above argument one can deduce that the
theorem holds for $L=[t_0, kT_0/2]$ for any integer $k\geq2$. This
completes the proof of the first claim of the theorem.

First observe that by~(\ref{e10}), for each fixed $x\in\bR^d$ and
for every $0\leq t_1<t_2<\cdots< t_k$,
the distributions of $\{(X^{b_n}_{t_1},\ldots, X^{b_n}_{t_k}),
\PP_x\}$ form a tight sequence.
Next, by the same argument as that for~(\ref{enew12}),
\[
\PP_x\bigl(X^{b_n}_s \notin B(x,r)\bigr) \le p \qquad\mbox{for all } n \ge1,
0 \le s \le t\mbox{ and } x \in\bR^d\vadjust{\goodbreak}
\]
implies
\[
\PP_x\Bigl(
\sup_{s\leq t} |X^{b_n}_t-X^{b_n}_0| \ge2r
\Bigr)= \PP_x\bigl(
\tau^{b_n}_{B(x, 2r)} \le t \bigr) \le2p \qquad\mbox{for all }
n \ge1, x \in\bR^d.
\]
Hence by~(\ref{e10}) and the same argument leading to
(\ref{econ}), we have for every $r>0$,
\[
\lim_{t \downarrow0}\sup_{n \ge1, x \in\bR}\PP_x\Bigl( \sup
_{s\leq
t} |X^{b_n}_t-X_0^{b_n}| \ge2r \Bigr) =0.
\]
Thus it follows from the Markov property and~\cite{Gr}, Theorem 2
(see also~\cite{EK}, Corollary 3.7.4, and~\cite{Bi}, Theorem 3) that,
for
each $x\in\bR^d$, the laws of $\{X^{b_n}, \PP_x\}$ form a tight sequence
in the Skorohod space $D([0, \infty), \bR^d)$. Combining\vspace*{2pt} this and
Theorem~\ref{t34} with~\cite{EK}, Corollary 4.8.7, we get that
$X^{b_n}$ converges\vspace*{1pt} to $X^b$ weakly. It follows directly from the
definition of Skorohod topology on
$D([0, \infty), \bR^d)$
(see, e.g.,~\cite{EK}, Section 3.5) that $\{t<\tau^b_U\}$ and
$\{t>\tau^b_{\bar U}\}$ are disjoint open subsets in
$D([0, \infty), \bR^d)$. Thus the boundary of
$\{t<\tau^b_U\}$ in $D([0, \infty), \bR^d)$ is contained in
$\{\tau^b_U\leq t\leq\tau^b_{\bar U} \}$. Note that, by the strong
Markov property,
\begin{eqnarray*}
\PP_x(\tau^b_U< \tau^b_{\bar U})
&=&\PP_x(\tau^b_U< \tau^b_U+\tau^b_{\bar U} \circ
\theta_{ \tau^b_U}, X^b_{\tau^b_U}\in\partial U)\\
&=& \PP_x(0< \tau^b_{\bar U} \circ\theta_{ \tau^b_U},
X^b_{\tau^b_U}\in\partial U)\\
&=&\PP_x\bigl(\PP_{ X^b_{\tau^b_U}}
(0< \tau^b_{\bar U}); X^b_{\tau^b_U}\in\partial U
\bigr) = 0.
\end{eqnarray*}
The last equality follows from the regularity of
$\bar U$; that is, $\PP_{z}( \tau^b_{\bar U}=0)=1$
for every $z \in\partial U$; see Proposition~\ref{preg}.
Therefore, using the L\'evy system for~$X^b$,
\begin{eqnarray*}
\PP_x(\tau^b_U\leq t\leq\tau^b_{\bar U})&=&
\PP_x(\tau^b_U= t =\tau^b_{\bar U})\\
&\leq& \PP_x (X^b_t \in\partial U)+ \PP_x(t=\tau^b_U \mbox{
and }
X^b_{\tau_U-} \not= X^b_{\tau_U}) \\
&=&\int_{\partial U}p^b(t, x, y) \,dy+0=0,
\end{eqnarray*}
which implies that the boundary of $\{t<\tau^b_U\}$ in $D([0,
\infty), \bR^d)$ is $\PP_x$-null for every $x\in U$. For every $f\in
C_b(\overline U)$, $f(X^{b}_t){\mathbf1}_{\{t<\tau^{b}_U\}}$ is a
bounded function on $D([0, \infty), \bR^d)$ with discontinuity
contained in the boundary of $\{t<\tau^b_U\}$. Thus we have (cf.
Theorem 2.9.1(vi) in~\cite{Du})
%
%
\begin{equation}\label{enew13}
\lim_{n\to\infty}\E_x\bigl[f(X^{b_n}_t){\mathbf1}_{\{t<\tau^{b_n}_U\}
}\bigr]
=\E_x\bigl[f(X^{b}_t){\mathbf1}_{\{t<\tau^{b}_U\}}\bigr].
\end{equation}

Given $f\in C_b(\overline U)$ and $\eps>0$, choose $T>1$ large such
that
\[
2C_1C_2^{-1}\|f\|_\infty e^{-C_2T} < \eps,
\]
where $C_1$ and
$C_2$ are the constants in Lemma~\ref{dec} with $D=U$. By the bounded
convergence theorem and Fubini's theorem, from~(\ref{enew13}) we
have
\begin{eqnarray*}
\lim_{n\to\infty}\E_x\biggl[\int_0^T f(X^{b_n}_t){\mathbf1}_{\{t<\tau
^{b_n}_U\}} \,dt\biggr]
&=&
\lim_{n\to\infty}\int_0^T \E_x\bigl[f(X^{b_n}_t){\mathbf1}_{\{t<\tau
^{b_n}_U\}}\bigr]\,dt\\
&=&\E_x\biggl[\int_0^T f(X^{b}_t){\mathbf1}_{\{t<\tau^{b}_U\}}\,dt\biggr].
\end{eqnarray*}
On the other hand, by the choice of $T$ and the fact that $C_1$ and
$C_2$ depending only on $d$, $\alpha$, $\operatorname{diam}(U)$ and $b$, with
the dependence on $b$ only through the rate at which $
M^\alpha_{|b|}(r)$ goes to zero, we have by Lemma~\ref{dec}
\begin{eqnarray*}
&&\E_x\biggl[\int_T^\infty f(X^{b_n}_t){\mathbf1}_{\{t<\tau^{b_n}_U\}}
\,dt\biggr]
+ \E_x\biggl[ \int_T^\infty f(X^{b}_t){\mathbf1}_{\{t<\tau^{b}_U\}}
\,dt\biggr] \\
&&\qquad \le \|f\|_\infty\int_T^\infty\biggl( \int_D \bigl(p^{b_n}_D (t,x,y)+
p^{b}_D (t,x,y)\bigr) \,dy \biggr) \,dt \\
&&\qquad \le 2 C_1\|f\|_\infty\int_T^\infty
e^{-C_2t } \,dt < \eps.
\end{eqnarray*}
This completes the proof of~(\ref{egcon}).
\end{pf}

As immediate consequences of~(\ref{egcon}) and Propositions
\ref{tufgfnest1} and~\ref{tufgfnest2}, we get the following:
\begin{theorem}\label{tufgfnest11}
There exists a constant $r_*
=r_*(d, \alpha, b)>0$
with the dependence on $b$ only via the rate at which
$M^\alpha_{|b|}(r)$ goes to zero such that for any ball $B=B(x_0,
r)$ of radius $r\le r_*$,
\[
{2}^{-1}G_B(x, y)\le G^{b}_B(x, y)\le2 G_B(x, y),\qquad x, y\in B.
\]
\end{theorem}
\begin{theorem}\label{tufgfnest21}
For every $C^{1,1}$ open set $D$ with the characteristic $(R_0,
\Lambda_0)$, there exists a constant $r_0 =r_0(d, \alpha, R_0,
\Lambda_0, b) \in(0, (R_0 \wedge1)/8]$ with the dependence on $b$
only via the rate at which $M^\alpha_{|b|}(r)$ goes to zero such
that for any for any $z\in\partial D$ and $r\le r_0$, we have
%
%
\begin{equation}\label{eempaok1}\qquad
2^{-1} G_{U_{(z, r)}}(x,y) \le G^{b}_{U_{(z, r)}}(x,y) \le2
G_{U_{(z, r)}}(x,y),\qquad x, y\in U_{(z, r)}.
\end{equation}
\end{theorem}

We will need the above two results later on.

\section{Duality}

In this section we assume that $E$ is an arbitrary bounded open set
in $\bR^d$. We will discuss some basic properties of $X^{b,E}$ and
its dual process under some reference measure. The results of this
section will be used later in this paper.

By Theorem~\ref{tconti} and Corollary~\ref{csp}, $X^{b,E}$ has a
jointly continuous and strictly positive transition density
$p^b_E(t,x,y)$. Using\vadjust{\goodbreak} the continuity of $p^b_E(t,x,y)$ and the
estimate
\[
p^b_E(t,x,y) \le p^b(t,x,y) \le c_1e^{c_2t} \biggl( t^{-d/\alpha}
\wedge
\frac{t}{|x-y|^{d+\alpha}}\biggr),
\]
the proof of the next proposition is easy. We omit the details.
\begin{prop}\label{pHC}
$X^{b,E}$ is a Hunt process, and it satisfies the strong Feller
property, that is, for every $f \in L^{\infty}(E)$,
$P^E_tf(x):=\E_x[f(X^{b,E}_t)]$ is bounded and continuous in $E$.
\end{prop}

Define
\[
h_E(x) := \int_E G^b_E(y,x) \,dy \quad\mbox{and}\quad \xi
_E(dx):=h_E(x) \,dx.
\]
The following result says that $\xi_E$ is a reference measure for $X^{b,E}$.
\begin{prop}\label{em}
$\xi_E$ is an excessive measure for $X^{b,E}$, that is, for every
Borel function $f \ge0$,
\[
\int_E f(x) \xi_E(dx) \ge\int_E \E_x[f(X^{b,E}_{t})]
\xi_E(dx).
\]
Moreover, $h_E$ is a strictly positive, bounded continuous function on $E$.
\end{prop}
\begin{pf}
By the Markov property, we have for any Borel function $f \ge0$
and $x\in E$,
\begin{eqnarray*}
\int_E\E_y[f(X^{b,E}_{t})] G^b_E(x,y) \,dy
&=& \E_x \int_0^\infty\E_{X^{b,E}_s}[f(X^{b,E}_{t})]
\,ds\\
&=& \int_0^\infty\E_{x}[f(X^{b,E}_{t+s})] \,ds \\
&\le&
\int_E f(y) G^b_E(x,y) \,dy.
\end{eqnarray*}
Integrating with respect to $x$, we get by Fubini's theorem
\[
\int_E\E_y[f(X^{b,E}_{t})] h_E(y) \,dy
\le\int_E f(y) h_E(y) \,dy.
\]
The second claim follows from~(\ref{Gbd}), the continuity of
$G^b_E$ and the strict positivity of $p^b_E$ (Corollary
\ref{csp}).
\end{pf}

We define a transition density with respect to the reference measure
$\xi_E$ by
\[
\overline{p}^b_E(t,x,y) := \frac{ p^b_E(t,x,y)}{h_E(y)}.
\]
Let
\[
\overline{G}{}^b_E(x,y) :=
\int^{\infty}_0 \overline{p}^b_E(t, x, y)\,dt = \frac{ G^b_E(x,y)}{h_E(y)}.
\]
Then $\overline{G}{}^b_E(x,y)$ is
the Green function of $X^{b,E}$ with respect to the reference measure~$\xi_E$.

Before we discuss properties of $ \overline{G}{}^b_E(x,y)$,
let us first recall some definitions.

\begin{defn}\label{dhar}
Suppose that $U$ is an open subset of $E$. A Borel function $u$ on
$E$ is said to be:
\begin{longlist}
\item harmonic in $U$ with respect to $X^{b, E}$ if
%
%
\begin{equation}\label{AV}
u(x) = \E_x[u(X^{b, E}_{\tau^b_{B}})],\qquad
x\in B,
\end{equation}
for every bounded open set $B$ with $\overline{B}\subset U$;

\item excessive with respective to $X^{b, E}$ if $u$ is
nonnegative and for every $t > 0$ and $x\in E$
\[
u(x) \geq \E_x[u(X^{b, E}_{t})]
\quad\mbox{and}\quad u(x) = \lim_{t \downarrow0}
\E_x[u(X^{b, E}_t)];
\]

\item a potential with respect to $X^{b, E}$ if it is excessive
with respect to $X^{b, E}$ and for every sequence $\{U_n\}_{n \ge1}$
of open sets with
$\overline{U}_n \subset U_{n+1}$ and $\bigcup_n U_n = E$,
\[
\lim_{n \to\infty}
\E_x[u(X^{b, E}_{\tau^b_{U_n}})] = 0,\qquad
\xi_E\mbox{-a.e. } x\in E;
\]

\item a pure potential with respect to $X^{b, E}$ if
it is a potential with respect to $X^{b, E}$ and
\[
\lim_{t \to\infty} \E_x[u(X^{b, E}_t)] = 0,\qquad
\xi_E \mbox{-a.e. } x\in E;
\]

\item
regular harmonic with
respect to $X^{b, E}$ in $U$ if $u$ is harmonic with respect to
$X^{b, E}$ in $U$ and~(\ref{AV}) is true for $B=U$.
\end{longlist}
\end{defn}

We list some properties of the Green function $\overline{G}{}^b_E(x,y)$
of $X^{b, E}$ that we will need later:
\begin{longlist}[(A1)]
\item[(A1)]
$\overline{G}{}^b_E(x,y) > 0$ for all $(x,y) \in E \times E$;
$\overline{G}{}^b_E(x,y) =\infty$ if and only if $x=y \in E$.
\item[(A2)]
For every $x \in E$, $\overline{G}{}^b_E(x, \cdot)$ and
$\overline{G}{}^b_E( \cdot,x)$ are extended continuous in~$E$.
\item[(A3)]
For every compact subset $K$ of $E$, $\int_K \overline{G}{}^b_E(x,y)
\xi_E(dy) < \infty$.
\end{longlist}

(A3) follows from~(\ref{Gbd}) and Proposition~\ref{em}. Both (A1)
and (A2) follow from~(\ref{Gbd}), Proposition~\ref{em},
domain monotonicity of Green functions and the lower bound in
(\ref{egege}).

From (A1)--(A3), we know that the process $X^{b, E}$ satisfies the
condition~($R$) on~\cite{CW}, page 211, and the conditions\vadjust{\goodbreak} (a) and (b) of
\cite{CW}, Theorem~5.4. It follows from~\cite{CW}, Theorem~5.4, that
$X^{b, E}$ satisfies Hunt's Hypothesis ($B$). Thus by~\cite{CW},
Theorem 13.24, $X^{b, E}$ has a dual process $\wh X^{b, E}$, which is
a standard process.

In addition, we have the following.

\begin{longlist}[(A4)]
\item[(A4)]
For each $y$, $ x \mapsto\overline{G}{}^b_E(x,y)$ is excessive with
respect to $X^{b, E}$ and harmonic with respect to $X^{b, E}$ in $E
\setminus\{y\}$. Moreover, for every open subset $U$ of $E$, we
have
%
%
\begin{equation}\label{eGH}
\E_x[\overline{G}{}^b_E(X^{b, E}_{T^b_{U}},y)]
=\overline{G}{}^b_E(x,y),\qquad (x,y) \in E \times U,
\end{equation}
where $T^b_U:=\inf\{ t > 0\dvtx X^{b, E}_t \in U\}$. In particular,
for every $y \in E$ and $\eps> 0$, $\overline{G}{}^b_E( \cdot, y)$
is regular harmonic in $E \setminus B(y, \eps)$ with respect to
$X^{b, E}$.
\end{longlist}
\begin{pf*}{Proof of (A4)}
It follows from
\cite{CR}, Proposition 3, and~\cite{L}, Theorem~2 on page 373,
that, to prove
(A4),
it suffices to show
that, for any $x\in E\setminus U$, the
function
\[
y\mapsto\E_x[\overline{G}{}^b_E(X^{b, E}_{T_{U}},y
)]
\]
is continuous on $U$. (See the proof of~\cite{P}, Theorem 1.)
Fix $x\in E\setminus U$ and $y\in U$. Put
$r:=\delta_U(y)$.
Let $\wh y \in B(y, r/4)$.
It follows from~(\ref{elevysE}) and~(\ref{Gbd}) that, for any
$\delta\in(0,
\frac{r}2)$,
\begin{eqnarray*}
&&\E_x[\overline{G}{}^b_E(X^{b, E}_{T^b_{U}},\wh y);
X^{b, E}_{T^b_{U}}\in B(y, \delta)]\\
&&\qquad= \int_{B(y,
\delta)}\biggl(\int_{E\setminus U}G^b_{E\setminus U}(x, w)\frac
{\mathcal{A}(d, - \alpha)}{|w-z|^{d+\alpha}}\,dw\biggr)\overline
{G}^b_E(z, \wh y)\,dz\\
&&\qquad\le \frac{c_1}{\inf_{\wtt y \in\overline{B(y, r/4)}} h_E(\wt y)}
\\
&&\qquad\quad{}\times\int_{B(y, \delta)}\biggl(\int_{E\setminus
U}\frac1{|x-w|^{d-\alpha}}\frac{1}{|w-z|^{d+\alpha}}\,dw\biggr)\,
\frac{dz}{|z-\wh y|^{d-\alpha}}.
\end{eqnarray*}
Thus, for any $\epsilon>0$, there is a $\delta\in(0, \frac{r}2)$
such that
%
%
\begin{equation}\label{enew15}\quad
\E_x[\overline{G}{}^b_E(X^{b, E}_{T_{U}},y); X^{b,
E}_{T_{U}}\in B(y,
\delta)] \le \frac{\eps}4 \qquad\mbox{for every } \wh y
\in B(y, r/4).
\end{equation}
Now we fix this $\delta$ and let $\{y_n\}$ be a sequence of points
in $B(y, r/4)$ converging to $y$. Since the function $(z, u)\mapsto
\overline{G}{}^b_E(z, u)$ is bounded and continuous in $(E\setminus
B(y, \delta))\times B(y, \frac{\delta}2)$, we have\vspace*{1pt} by the
bounded
convergence theorem that there exists $n_0>0$ such that for all
$n\ge n_0$,
%
%
\begin{eqnarray}\label{enew16}
&&|\E_x[\overline{G}{}^b_E(X^{b, E}_{T_{U}},y);
X^{b, E}_{T_{U}}\in B(y,
\delta)^c ]\nonumber\\[-8pt]\\[-8pt]
&&\qquad{}-\E_x[\overline{G}{}^b_E (X^{b,
E}_{T_{U}},y_n); X^{b,
E}_{T_{U}}\in B(y, \delta)^c]|
\le\frac{\eps}2. \nonumber
\end{eqnarray}
Since $\eps>0$ is arbitrary, combining~(\ref{enew15}) and(\ref{enew16}),
the proof of (A4) is now complete.
\end{pf*}
\begin{theorem}\label{tp}
For each $y \in E$, $ x \mapsto\overline{G}{}^b_E(x,y)$ is a pure
potential with respect to $X^{b, E}$. In fact, for every sequence
$\{U_n\}_{n \ge1}$ of open sets with $\overline{U_n} \subset
U_{n+1}$ and $\bigcup_n U_n = E$,
$
\lim_{n \to\infty}
\E_x[\overline{G}{}^b_E(X^{b, E}_{\tau^b_{U_n}},y)]= 0$ for every $ x
\not= y$ in $E$.
Moreover, for every $x,y \in E$, we have
$
\lim_{t \to\infty} \E_x[ \overline{G}{}^b_E(X^{b, E}_t, y)] = 0.
$
\end{theorem}
\begin{pf}
For $y\in E$, let $X^{b, E,y}$ denote the $h$-conditioned process
obtained from $X^{b, E}$ with $h(\cdot)= \overline{G}{}^b_E(\cdot,
y)$, and let $\E_x^{y}$ denote the expectation for $X^{b, E,y}$
starting from $x\in E$.

Let $x \not= y \in E$. Using (A1), (A2), (A4) and the strict
positivity of $\overline{G}{}^b_E$, and applying~\cite{LZ}, Theorem 2,
we get that the lifetime $\zeta^{b, E, y}$ of $X^{b, E,y}$ is finite
$\PP_x^{y}$-a.s. and
%
%
\begin{equation}\label{elim}
\lim_{t\uparrow\zeta^{b, E, y}}X^{b, E,y}_t =y,\qquad
\PP_x^{y}\mbox{-a.s.}
\end{equation}
Let $\{E_k, k \ge1\}$ be an increasing sequence of
relatively compact open subsets of $E$ such that
$E_k \subset\overline{E}_k \subset E$ and
$\bigcup_{k=1}^\infty E_k =E$. Then
\[
\E_x [ \overline{G}{}^b_E(X^{b, E}_{\tau^b_{E_k}}, y)
] =
\overline{G}{}^b_E(x, y) \PP_x^y ( \tau^b_{E_k}<\zeta^{b, E,
y}).
\]
By~(\ref{elim}), we have $\lim_{k \to\infty} \PP_x^y (
\tau^b_{E_k}<\zeta^{b, E, y}) =0$. Thus
\[
\lim_{k \to\infty} \E_x [
\overline{G}{}^b_E(X^{b, E}_{\tau^b_{E_k}}, y)] =0.
\]

The last claim of the theorem is easy. By~(\ref{est61}) and
(\ref{Gbd}), for every $x,y \in E$, we have
\[
\E_x[ \overline{G}{}^b_E(X^{b, E}_t, y)] \le
\frac{c}{t^{{d}/{\alpha}}h_E(y)}
\int_E \frac{dz}{|z-y|^{d-\alpha}},
\]
which converges to zero as $t$ goes to $\infty$.
\end{pf}

We note that
\[
\int_E \overline{G}{}^b_E (x,y) \xi_E (dx) \le
\frac{\|h_E\|_{\infty}}{h_E(y)}
\int_E G^b_E (x,y) \,dx = \|h_E\|_{\infty} < \infty.
\]
So we have
\begin{longlist}[(A5)]
\item[(A5)] for every compact subset $K$ of $E$, $\int_K
\overline{G}{}^b_E(x,y) \xi_E(dx) < \infty$.
\end{longlist}
Using (A1)--(A5),~(\ref{Gbd}) and Theorem~\ref{tp} we get from
\cite{L1,L} that $X^{b, E}$ has a Hunt process as a dual.
\begin{theorem} \label{dualpr}
There exists a transient Hunt process $\wh X^{b, E} $ in $E$ such
that $\wh X^{b, E} $ is a strong dual\vadjust{\goodbreak} of $X^{b, E}$ with respect to
the measure $\xi_E$; that is, the density of the semigroup $\{ \wh
P^E_t\}_{t\geq0}$ of $\wh X^{b, E}$ is given by
$\overline{p}^b_E(t,y,x)$ and thus
\[
\int_E f(x) P^E_t g(x) \xi_E(dx)=\int_E g(x) \wh P^E_t f(x) \xi_E(dx)
\qquad\mbox{for all } f, g\in L^2(E, \xi_E).
\]
\end{theorem}

\begin{pf}
The existence of a dual Hunt process $\wh X^{b, E}$ is proved
in
\cite{L1,L}. To show $\wh X^{b, E}$ is transient, we need to show
that for every compact subset~$K$ of $E$, $\int_K \overline{G}{}^b_E
(x,y) \xi_E (dx)$ is bounded. This is just (A5) above.
\end{pf}

In Theorem~\ref{tls4xb}, we have determined a L\'evy system $(N,
H)$ for $X^b$ with respect to the Lebesgue measure $dx$. To derive
a L\'evy system for $\wh X^{b, E}$, we need to consider a L\'evy
system for $X^{b, E}$ with respect to the reference measure
$\xi_E(dx)$. One can easily check that, if
\begin{eqnarray*}
\overline{N}{}^E(x,dy)&:=&\frac{J(x,y)}{h_E(y)} \xi_E(dy) \qquad\mbox{for
} (x,y) \in E \times E, \\
\overline{N}{}^E(x,\partial)&:=&\int_{E^c}J(x,y)\,dy \qquad\mbox{for }
x\in E,
\end{eqnarray*}
and $\overline{H}{}^E_t:=t$,\vspace*{1pt} then $(\overline{N}{}^E, \overline{H}{}^E)$ is
a L\'evy system for $X^{b, E}$ with respect to the reference
measure $\xi_E(dx)$.
It follows from~\cite{G} that a L\'evy system $(\wh N^E, \wh H^E)$
for $\wh X^{b,E}$ satisfies $\wh H^E_t=t$ and
\[
\wh N^E(y,dx) \xi_E(dy) =\overline{N}{}^E(x,dy)\xi_E (dx).
\]
Therefore, using $J(x,y)=J(y,x)$, we have
for every stopping time $T$ with respect to the filtration of
$\wh X^{b, E}$,
%
%
\begin{eqnarray}\label{elevy2}
&&\E_x \biggl[\sum_{s\le T} f(s,\wh X^{b, E}_{s-}, \wh X^{b, E}_s)
\biggr] \nonumber\\
&&\qquad= \E_x \biggl[ \int_0^T \biggl( \int_{E} f (s,\wh X^{b,E}_s,
y) \frac{J( \wh X^{b, E}_s, y )} {h_E (\wh X^{b, E}_s)}\xi_E(dy)
\biggr) \,d\wh H^E_s \biggr]\\
&&\qquad= \E_x \biggl[ \int_0^T \biggl( \int_{E} f (s,\wh X^{b,E}_s, y)
\frac{J( \wh X^{b, E}_s, y )h_E (y)} {h_E (\wh X^{b, E}_s)}\,dy
\biggr) \,ds \biggr].\nonumber
\end{eqnarray}
That is,
\[
\wh N^E (x, dy)=\frac{J(x, y)h_E(y)}{h_E(x)} \,dy.
\]

Let
\[
P^{b, E}_t f(x):= \int_E \overline{p}^b_E(t,x,y) f(y)\xi_E(dy)
\]
and
\[
\wh P^{b, E}_t f(x):= \int_E \overline{p}^b_E(t,y,x)
f(y)\xi_E(dy).
\]

For any open subset $U$ of $E$, we use $\wh X^{b, E, U}$ to denote
the subprocess of $\wh X^{b, E}$ in $U$,
that is, $\wh X^{b, E, U}_t(\omega)=\wh X^{b, E}_t(\omega)$ if $t<
\wh\tau^{b, E}_U(\omega)$ and $\wh X^{b, E, U}_t(\omega)=\partial$
if $t\geq\wh\tau^{b, E}_U(\omega)$, where $\wh\tau^{b,
E}_U:=\inf\{t >0 \dvtx \wh X^{b, E}_t \notin U\}$, and $\partial$ is the
cemetery state. Then by~\cite{S}, Theorem 2, and Remark 2 following
it,
$X^{b, U}$
and $\wh X^{b, E, U}$ are dual processes with respect to
$\xi_E$. Now we let
%
%
\begin{equation}\label{whp}
\wh p^{b, E}_{U}(t,x,y) := \frac{p^b_U(t,y,x)h_E(y)}{h_E(x)}.
\end{equation}
By the joint continuity of $p^b_U(t,x,y)$ (Theorem~\ref{tconti})
and the continuity and positivity of $h_E$ (Proposition~\ref{em}),
we know that $
\wh p^{b, E}_{U}(t,\cdot,\cdot)$ is jointly continuous
on $U \times U$. Thus we have the following.
\begin{theorem}\label{tkilldual}
For every open subset $U$, $ \wh p^{b, E}_{U}(t,x,y)$ is strictly
positive and jointly continuous on $U \times U$ and is the
transition density of $\wh X^{b, E, U}$ with respect to the Lebesgue
measure. Moreover,
%
%
\begin{equation}\label{edGk}
\wh G^{b, E}_{U}(x,y) := \frac{G^b_U(y,x)h_E(y)}{h_E(x)}
\end{equation}
is the Green function of $\wh X^{b, E, U}$ with respect to the Lebesgue measure
so that for every nonnegative Borel function $f$,
\[
\E_x\biggl[\int_0^{\wh\tau^{b, E}_U} f(\wh
X^{b,E}_t)\,dt\biggr] = \int_U \wh G^{b, E}_{U}(x,y) f(y) \,dy.
\]
\end{theorem}

\section{Scaling property and uniform boundary Harnack principle}

In this section, we first study the scaling property of $X^b$, which
will be used later in this paper.

For $\lambda>0$, let $Y^{b, \lambda}_t:=\lambda
X^b_{\lambda^{-\alpha}t}$. For any function $f$ on $\bR^d$, we
define $f^\lambda( \cdot)=f(\lambda\cdot)$. Then we have
\[
\E_x[f(Y^{b, \lambda}_t)]= \int_{\bR^d}
p^b(\lambda^{-\alpha}t,\lambda^{-1}x,y) f^{\lambda}(y)\,dy.
\]
It follows from Theorem~\ref{T11}(iii) that for any $f,g \in
C^\infty_c(\bR^d)$,
\begin{eqnarray*}
&&\lim_{t \downarrow0} \int_{\bR^d} t^{-1} \bigl(\E_x[f(Y^{b,
\lambda}_t)]-f(x)\bigr) g(x)\,dx \\[-0.5pt]
&&\qquad=\lim_{t \downarrow0} \int_{\bR^d}\lambda^{-\alpha}
(\lambda^{\alpha} t)^{-1} \bigl(P^b_{\lambda^{-\alpha}t}
f^\lambda(\lambda^{-1}
x)- f^\lambda(\lambda^{-1} x)\bigr) g^\lambda(\lambda^{-1} x)\,dx
\\[-0.5pt]
&&\qquad=\lim_{t \downarrow0}
\int_{\bR^d}\lambda^{d-\alpha}(\lambda^{\alpha} t)^{-1}
\bigl(P^b_{\lambda^{-\alpha}t}f^\lambda(z)-f^\lambda(z)\bigr) g^\lambda(z)\,dz
\\[-0.5pt]
&&\qquad= \lambda^{d-\alpha}\int_{\bR^d}
\bigl(-(-\Delta)^{\alpha/2}f^\lambda(z)+b(z) \cdot
\nabla f^\lambda(z)\bigr) g^\lambda(z)\,dz\\[-0.5pt]
&&\qquad= \lambda^{d-\alpha}\int_{\bR^d}
\bigl(-(-\Delta)^{\alpha/2}f^\lambda(z)+\lambda b(z)
\cdot\nabla f(\lambda z)\bigr) g(\lambda z)\,dz\\[-0.5pt]
&&\qquad= \int_{\bR^d}
\bigl(-(-\Delta)^{\alpha/2}f(x)+\lambda^{1-\alpha}
b(\lambda^{-1} x) \cdot\nabla f(x)\bigr) g(x)\,dx.
\end{eqnarray*}
Thus $\{\lambda X^{b, D}_{\lambda^{-\alpha} t}, t\geq0\}$
is the subprocess of $X^{\lambda^{1-\alpha}b(\lambda^{-1} \cdot)}$
in $\lambda D$.
So for any $\lambda>0$, we have
%
%
\begin{eqnarray}\label{escaling}
p^{\lambda^{1-\alpha} b(\lambda^{-1} \cdot)}_{\lambda D} ( t, x, y)&=&
\lambda^{-d} p^{b}_D (\lambda^{-\alpha}t, \lambda^{-1} x, \lambda
^{-1} y)\hspace*{25.8pt}\hspace*{25.8pt}\nonumber\\[-8pt]\\[-8pt]
&&\eqntext{\mbox{for } t>0 \mbox{ and } x, y \in\lambda D,}
\end{eqnarray}
%
%
\begin{equation}\label{escaling1}
\qquad G^{\lambda^{1-\alpha} b(\lambda^{-1} \cdot)}_{\lambda D} (x, y)=
\lambda^{\alpha-d} G^{b}_D (\lambda^{-1} x, \lambda^{-1} y)
\qquad\mbox{for } x, y \in\lambda D.
\end{equation}

Define
%
%
\begin{equation}\label{e72}
b_\lambda(x):=\lambda^{1-\alpha} b(x /\lambda) \qquad\mbox{for }
x\in\bR^d.
\end{equation}
Then we have
\begin{eqnarray*}
M^\alpha_{|b_{\lambda}|}(r)
&=&\lambda^{1-\alpha}\sum^d_{i=1}\sup_{x\in\bR^d}\int_{|x-y|\le
r}\frac{|b^i|(\lambda^{-1}y)\,dy}
{|x-y|^{d+1-\alpha}}\\
&=&\sum^d_{i=1}\sup_{\wh x\in\bR^d}\int_{|\wh x-z|\le\lambda^{-1}
r}\frac{|b^i|(z)\,dz} {|\wh x-z|^{d+1-\alpha}}=M^\alpha_{
|b|}(\lambda^{-1}r).
\end{eqnarray*}
Therefore
for every $\lambda\ge1$ and $r>0$,
%
%
\begin{equation}\label{eMgoOK}
M^\alpha_{|b_\lambda|}(r) = M^\alpha_{ |b|}(\lambda^{-1}r) \le
M^\alpha_{ |b|}(r).
\end{equation}

In the remainder of this paper, we fix a bounded $C^{1,1}$ open set
$D$ in $\bR^d$ with $C^{1,1}$ characteristics $(R_0, \Lambda_0)$ and
a ball $E\subset\bR^d$ centered at the origin so that $D\subset
\frac14 E$. Define
%
%
\begin{equation}\label{MMM}
M:=M(b, E):= \sup_{x,y \in3E/4}\frac{h_E(x)}{h_E(y)},
\end{equation}
which is a finite positive constant no less than 1.
Note that, in
view of scaling property~(\ref{escaling1}), we have
%
%
\begin{equation}\label{e73}
M(b, E)= M(b_\lambda, \lambda E).
\end{equation}

Although $E$ and $D$ are fixed, the constants in all the results of
this section will depend only on $d, \alpha, R_0, \Lambda_0$, $b$
and $M$ (not the diameter of $D$ directly) with the dependence on
$b$ only via the rate at which $M^\alpha_{|b|}(r)$ goes to zero. In
view of~(\ref{eMgoOK}) and~(\ref{e73}), the results of this
section in particular hold for $\sL^{b_\lambda}$ (equivalently,
for~$X^{b_\lambda}$) and the pair $(\lambda D, \lambda E)$ for every
$\lambda\geq1$.

In the remainder of this section, we will establish a uniform
boundary Harnack principle on $D$ for certain harmonic functions for
$X^{b, E}$ and $\wh X^{b, E}$. Since\vspace*{1pt} the arguments are mostly
similar for $X^{b, E}$ and $\wh X^{b, E}$, we will only give the
proof for~$\wh X^{b, E}$.

A real-valued function $u$ on $E$ is said to be harmonic in an open
set $U\subset E$ with respect to $\wh X^{b, E}$ if for every
relatively compact open subset $B$ with $\overline B\subset U$,
%
%
\begin{equation}\label{ehar}\qquad
\E_x [ | u(\wh X^{b, E}_{\wh\tau^{b,E}_{B}}
)|
]<\infty\quad\mbox{and}\quad u(x)= \E_x [ u(\wh
X^{b, E}_{\wh\tau^{b,E}_{B}})]
\qquad\mbox{for every } x\in
B.
\end{equation}
A real-valued function $u$ on $E$ is said to be regular
harmonic in an open set $U\subset E$ with respect to $\wh X^{b, E}$
if
(\ref{ehar}) is true with $B=U$.
Clearly, a regular harmonic function in $U$ is harmonic in $U$.

For any bounded open set $U$, define the Poisson kernel for $X^b$ of
$U$ as
\[
K_{U}^b (x,z):= \int_{ U} G^b_{U}(x,y)J(y,z)\,dy,\qquad (x,z) \in U
\times(\bR^d \setminus\overline{U}).
\]
When $U\subset E$, we define the Poisson kernel for $\wh X^{b, E} $
of $U\subset E$ as
%
%
\begin{equation}\label{eKK1}
\wh K_{U}^{b, E} (x,z):= \frac{h_E (z)}{h_E(x)}\int_{ U} G^b_{U}(y,
x)J(z, y) \,dy,\qquad (x,z) \in U
\times(E\setminus\overline{U}).\hspace*{-32pt}
\end{equation}
By~(\ref{elevysE}) and~(\ref{elevy2}), we have
\[
\E_x[f(X^{b, E}_{\tau^{b}_U});
X^b_{\tau^{b}_U-} \not=
X^b_{\tau^{b}_U} ] =\int_{\bar{U}^c} K^b_U(x,z)f(z)\,dz
\]
and
%
%
\begin{eqnarray}\label{elevyp}
&&\E_x[f(\wh X^{b,E}_{\wh
\tau^{b, E}_U}); \wh X^{b, E}_{\wh\tau^{b, E}_U-}
\not= \wh X^{b, E}_{\wh\tau^{b, E}_U} ] \nonumber\\
&&\qquad= \E_x \int_0^{\wh\tau^{b, E}_U} \biggl(
\int_{\bar{U}^c} f (z) \frac{J( \wh X^{b, E}_s, z )h_E (z)}
{h_E (\wh X^{b, E}_s)}\,dz
\biggr) \,ds \nonumber\\[-8pt]\\[-8pt]
&&\qquad= \int_U \frac{G^b_{U}(y, x)h_E (y)}{h_E (x)}
\int_{\bar{U}^c} f (z) \frac{J( y, z )h_E (z)} {h_E (y)}\,dz \,dy
\nonumber\\
&&\qquad= \int_{\bar{U}^c}\wh K^{b, E}_U(x,z)f(z)\,dz. \nonumber
\end{eqnarray}
\begin{lemma}\label{L20}
Suppose that $U$ is a bounded $C^{1,1}$ open set in $\bR^d$ with $U
\subset\frac12 E$ and $\operatorname{diam} (U) \le3r_*$ where $r_*$ is the
constant in Theorem~\ref{tufgfnest11}. Then
%
%
\begin{equation}\label{elpz1}
\PP_x(X^b_{\tau^b_U} \in\partial U)=0 \qquad\mbox{for every
} x\in U
\end{equation}
and
%
%
\begin{equation}\label{elpz2}
\PP_x(\wh X^{b, E}_{\wh\tau^{b, E}_U} \in\partial U)=0
\qquad\mbox{for every } x\in U.
\end{equation}
\end{lemma}
\begin{pf}
The proof is similar to that of~\cite{B}, Lemma 6. For our
readers' convenience, we are going to spell out the details of the
proof of~(\ref{elpz2}). Let $B_x:=B(x, \delta_U(x)/3)$. By
(\ref{elevy2}) we have for $x\in U$,
\begin{eqnarray*}
&&\PP_x\biggl(\wh X^{b, E}_{\wh\tau^{b, E}_{B_x}} \in
\biggl(\frac34E\biggr) \Bigm\backslash U \biggr)\\
&&\qquad= \int_{B_x} \frac{G^b_{B_x}(y,x)
h_E(y)}{h_E(x)}
\biggl( \int_{ (\frac34 E) \setminus U}\frac{J(y, z)h_E(z)}{h_E(y)}
\,dz \biggr) \,dy.
\end{eqnarray*}
Since $\operatorname{diam} (U) \le3r_*$, $\delta_U(x)/3 \le r_*$, thus by
Theorem~\ref{tufgfnest11},
for $x\in U$,
%
%
\begin{eqnarray}\label{esds1}
&&\PP_x\biggl(\wh X^{b, E}_{\wh\tau^{b, E}_{B_x}} \in\biggl(\frac34 E\biggr)
\Bigm\backslash U \biggr)\nonumber\\
&&\qquad\geq c_1\biggl(\inf_{u,v \in3E/4}\frac{h_E (u)}{h_E(v)}
\biggr)
\int_{B_x} G_{B_x}(x,y) \biggl( \int_{ (3E/4) \setminus U}
J(y,z) \,dz \biggr) \,dy\\
&&\qquad\ge c_1 M^{-1} \PP_x\biggl(X_{ \tau_{B_x}} \in\biggl(\frac34 E\biggr)
\Bigm\backslash U \biggr) , \nonumber
\end{eqnarray}
where $M$ is the constant defined in~(\ref{MMM}).
Let $V_x:=B(\delta_U(x)^{-1} x, 1/3)$.
By the scaling
property of
$X$,
%
%
\begin{eqnarray}\label{esds}
&& \PP_x\biggl(X_{ \tau_{B_x}} \in\biggl(\frac34 E\biggr) \Bigm\backslash U \biggr)
\nonumber\\
&&\qquad = \PP_{\delta_U(x)^{-1}x} \biggl(X_{ \tau_{\delta_U(x)^{-1}B_x}}
\in\delta_U(x)^{-1} \biggl(\frac34 E\biggr) \Bigm\backslash U \biggr) \\
&&\qquad = \int_{V_x} G_{V_x}(\delta_U(x)^{-1}x,a) \biggl( \int_{ \delta
_U(x)^{-1}(3E/4)
\setminus U} J(a,b) \,db \biggr) \,da.\nonumber
\end{eqnarray}
Let $z_x \in\partial U$ be such that $\delta_U(x)=|x-z_x|$. Since
$U$ is $C^{1,1}$, $\delta_U(x)^{-1}((\frac34 E) \setminus U) \supset
\delta_U(x)^{-1}(\frac34 E \setminus\frac12 E)$ and $\delta_U(x)
\le3 r_*$, there exists $\eta>0$ such that, under an appropriate
coordinate system, we have $z_x+\wh C \subset
\delta_U(x)^{-1}((\frac34 E) \setminus U)$ where
\[
\wh C:= \bigl\{y=(y_1,\ldots,y_d) \in\bR^d \dvtx
0<y_d<(12r_*)^{-1},
\sqrt{y_1^2+ \cdots+y_{d-1}^2 } < \eta y_d \bigr\}.
\]
Thus there is a constant $c_2>0$ such that
\[
\inf_{a \in V_x}\int_{
\delta_U(x)^{-1}((3E/4) \setminus U )} J(a,b) \,db \geq
c_2 > 0 \qquad\mbox{for every } x \in U.
\]
Combining this with~(\ref{esds1}) and~(\ref{esds}),
%
%
\begin{equation}\label{eLb}\quad
\inf_{x \in U} \PP_x\biggl(\wh X^{b, E}_{\wh\tau^{b, E}_{B_x}} \in
\biggl(\frac34 E\biggr) \Bigm\backslash U \biggr)
\ge{c_1c_2} M^{-1}
\E_w \bigl[\tau_{B(0, 1/3)} \bigr] \ge c_3 >0.
\end{equation}
On the other hand, since by~(\ref{elevy2})
$\PP_x(\wh X_{\wh\tau^{b, E}_{B_x}} \in
\partial U)=0$ for every $x \in U$, we have
\[
\PP_x(\wh X^{b, E}_{\wh\tau^{b, E}_{U}} \in\partial U)
=\E_x[ \PP_{\wh X^{b, E}_{\wh\tau^{b, E}_{B_x}}}( \wh X^{b,
E}_{\wh\tau^{b, E}_U} \in\partial U); \wh X^{b, E}_{\wh\tau
^{b, E}_{B_x}} \in U ].
\]
Thus inductively,
$
\PP_x(\wh X^{b, E}_{\wh\tau^{b, E}_{U}} \in\partial U)
=\lim_{k \to\infty} p_k(x)$,
where
\[
p_0(x):=\PP_x(\wh X^{b, E}_{\wh\tau^{b, E}_{U}} \in\partial
U)
\]
and
\[
p_{k }(x):= \E_x [ p_{k-1}
(\wh X^{b, E}_{\wh\tau^{b, E}_{B_x}}); \wh X^{b, E}_{\wh
\tau^{b, E}_{B_x}} \in U ] \qquad\mbox{for } k\geq1.
\]
By~(\ref{eLb}),
\[
\sup_{x \in U}p_{k+1}(x) \le (1-c_3)\sup_{x \in U}p_{k}(x)
\le (1-c_3)^{k+1} \to 0.
\]
Therefore,
$ \PP_x(\wh X^{b, E}_{\wh\tau^{b, E}_U} \in\partial U)=0
$ for every $x \in U$.
\end{pf}

Let $z\in\partial D$. We will say that a function $u\dvtx\bR^d\to\bR$
vanishes continuously on $ D^c \cap B(z, r)$ if $u=0$ on $ D^c \cap
B(z, r)$, and $u$ is continuous at every point of $\partial D\cap
B(z,r)$.
\begin{theorem}[(Boundary Harnack principle)]\label{ubhp} There exist positive
constants $c_1=c_1(d, \alpha, R_0, \Lambda_0, b)$ and $r_1=r_1(d,
\alpha, R_0, \Lambda_0, b)$ with the dependence on $b$ only via the
rate at which $M^\alpha_{|b|}(r)$ goes to zero such that for all
$z\in\partial D$, $r\in(0, r_1]$ and all function $u\ge0$ on
$\bR^d$ that is positive harmonic with respect to $X^{b}$ (or $\wh
X^{b,E}$, resp.) in $D\cap B(z, r)$, and vanishes
continuously on $ D^c \cap B(z, r)$ (or $ D^c$, resp.) we
have
\[
\frac{u(x)}{u(y)} \le c_1 M^2
\frac{\delta_D(x)^{\alpha/2}}{\delta_D(y)^{\alpha/2}},\qquad x,
y\in D\cap B(z, r/4).
\]
\end{theorem}
\begin{pf}
We only give the proof for $\wh X^{b,E}$. Recall that $r_*$ and
$r_0$ are the constants from Theorems~\ref{tufgfnest11} and
\ref{tufgfnest21}, respectively. Let $r_1= r_* \wedge r_0$, and fix
$r\in(0,r_1]$ throughout this proof. Recall that there exists
$L=L(R_0, \Lambda_0, d)$ such that for\vadjust{\goodbreak} every $z \in\partial D$ and
$r \le R_0/2$, one can find a $C^{1,1}$ open set $U=U_{(z, r)}$ with
$C^{1, 1}$ characteristic $(rR_0/L, \Lambda_0L/r)$ such that $D \cap
B(z, r/2) \subset U \subset D \cap B(z, r) $. Without loss of
generality, we assume $z=0$.

Note that, by the same proof as that of~\cite{CKSV}, Lemma 4.2,
every nonnegative function $u$ in $\bR^d$ that is harmonic with
respect to $\wh X^{b,E}$ in $D\cap B(0, r)$ and
vanishes
continuously on $ D^c$
is regular harmonic in
$D\cap B(0, r)$ with respect to $\wh
X^{b,E}$.

For all
functions $u \ge0$ on $E$ that is positive regular harmonic for
$\wh X^{b,E}$ in $D\cap B(0, r)$ and vanishing on $D^c$, by
(\ref{elevy2}) and Lemma~\ref{L20}, we have
%
%
\begin{eqnarray}\label{u12}
u (x) & = & \E_x[u (
\wh X^{b, E}_{\wh\tau^{b, E}_{U}}) ; \wh X^{b, E}_{\wh\tau^{b, E}_{U}}
\in D\setminus U]
\nonumber\\
& = &
\int_{D \setminus U} \wh K^{b, E}_{U} (x,w) u (w) \,dw \\
& = & \int_U G^b_{U}(y, x) \biggl( \int_{D \setminus U}
\frac{h_E(w)}{h_E(x)} J(w, y) u (w) \,dw \biggr) \,dy. \nonumber
\end{eqnarray}
Define
\begin{eqnarray*}
h_u (x) &:=&\E_x[u( X_{\tau_{U}})
; X_{\tau_{U}} \in D \setminus U]\\
&=& \int_U G_{U}(y, x) \biggl( \int_{D \setminus U}
J(w, y) u (w) \,dw \biggr) \,dy,
\end{eqnarray*}
which is positive regular harmonic for $X$ in $D\cap B(0, r/2)$ and
vanishing on $D^c$. Applying Theorem~\ref{tufgfnest21} to
(\ref{u12}), we get
%
%
\begin{equation}\label{e714}
c_1^{-1} M^{-1} h_u(x) \le
u (x) \le c_1 M h_u(x) \qquad\mbox{for } x\in D.
\end{equation}
By the boundary Harnack principle for $X$ in $C^{1,1}$ open sets
(see~\cite{CS1,SW}), there is a constant $c_2>1$ that
depends only on $R_0, \Lambda_0$, $d$ and $\alpha$ so that
\[
\frac{h_u(x)}{h_u(y)} \leq c_2
\qquad\mbox{for } x, y\in D\cap B(0, r/4) .
\]
Combining this with~(\ref{e714})
and the two-sided estimates on $G_U(x,y)$ we arrive at the conclusion of
the theorem.
\end{pf}

\section{Small time heat kernel estimates}

Our strategy is to first establish sharp two-sided estimates
on $p^b_D(t, x, y)$ at time $t=1$. We then use a scaling argument
to establish estimates for
$t \le T$.

We continue to fix a ball $E$ centered at the origin and a
$C^{1,1}$ open set $D \subset\frac14E$ with characteristics $(R_0,
\Lambda_0)$. Recall that $M>1$ is the constant defined in
(\ref{MMM}).

The next result follows from Proposition~\ref{P35-55},~(\ref{whp})
and~(\ref{MMM}).
\begin{prop}\label{P35-55-1}
For all $a_1\in(0, 1)$, $a_2, a_3, R>0$, there is a constant $c_1
=c_1(d, \alpha, a_1, a_2, a_3, R, M, b)>0$ with the dependence on
$b$ only via the rate at which $M^\alpha_{|b|}(r)$ goes to zero\vadjust{\goodbreak} such
that for all open ball $B(x_0, r) \subset\frac34 E$ with $r \le R$,
\[
\wh p^{b, E}_{B(x_0,r)}(t, x,y)\ge c_1 t^{-d /\alpha}
\qquad\mbox{for all } x, y\in B(x_0, a_1 r) \mbox{ and } t\in[
a_2r^\alpha, a_3 r^\alpha].
\]
\end{prop}

Again, we emphasize that the constants in all the results of the
remainder of this section (except Theorem~\ref{thmoppz} where the
constant also depends on $T$ for an obvious reason) will depend only
on $d, \alpha, R_0, \Lambda_0$, $M$ (not the diameter of $D$
directly) and $b$ with the dependence on $b$ only through the rate
at which $M^\alpha_{|b|}(r)$ goes to zero. In view of
(\ref{e72}),~(\ref{eMgoOK}) and~(\ref{e73}), in particular, all
the results of this section are applicable to $\sL^{b_\lambda}$ and
the pair $(\lambda D, \lambda E)$ for every $\lambda\geq1$.

Recall that $r_*$ and $r_0$ are the constants from Theorems
\ref{tufgfnest11} and~\ref{tufgfnest21},
respectively, which
depend only on
$d$, $\alpha$,
$R_0$, $\Lambda_0$ and
$b$ with the dependence on $b$ only
via the rate at which $M^\alpha_{|b|}(r)$ goes to zero.
\begin{lemma}\label{lemetd1r}
There is $c_1 =c_1(d, \alpha, R_0, r, M, \Lambda_0, b)>0$ with the
dependence on $b$ only via the rate at which $M^\alpha_{|b|}(r)$
goes to zero such that
for all $x \in D$
%
%
\begin{equation}\label{eqcw0}
\PP_x(\tau^{b}_{D}>1/4) \le c_1\bigl( 1\wedge
\delta_{D}(x)^{\alpha/2}\bigr)
\end{equation}
and
%
%
\begin{equation}\label{eqcw}
\PP_x(\wh\tau^{b, E}_{D}>1/4) \le c_1 \bigl( 1\wedge
\delta_{D}(x)^{\alpha/2}\bigr).
\end{equation}
\end{lemma}
\begin{pf}
We only give the proof of~(\ref{eqcw}). The proof of
(\ref{eqcw0}) is similar. Recall that there exists $L=L(R_0,
\Lambda_0, d)$ such that for every $z \in\partial D$ and $r \le
R_0$, one can find a $C^{1,1}$ open set $U_{(z, r)}$ with $C^{1, 1}$
characteristic $(rR_0/L, \Lambda_0L/r)$ such that $D \cap B(z, r/2)
\subset U_{(z, r)} \subset D \cap B(z, r) $. Clearly it suffices to
prove~(\ref{eqcw}) for $x \in D$ with $\delta_D(x) < (r_0\wedge r_*) /8$.

Choose $Q_x \in\partial D$ such that $\delta_D(x)=|x-Q_x|$, and
choose a $C^{1,1}$ open set $U:=U_{(Q_x,
(r_0\wedge r_*)/2)}$ with $C^{1, 1}$
characteristic $(
(r_0\wedge r_*)R_0/(2L), 2\Lambda_0L/ (r_0\wedge r_*))$ such that $D
\cap
B(Q_x,
(r_0\wedge r_*)/4) \subset U \subset D \cap B(Q_x,
(r_0\wedge r_*)/2) $.

Note that by~(\ref{edGk}),~(\ref{eKK1}) and Lemma~\ref{L20},
\begin{eqnarray*}
&&\PP_x(\wh\tau^{b, E}_{D}>1/4)\\
&&\qquad\leq
\PP_x(\wh\tau^{b, E}_{U}>1/4)+
\PP_x(\wh X^{b, E}_{\wh\tau^{b, E}_{U}}\in D)\\
&&\qquad\leq 4\E_x[\wh\tau^{b, E}_{U}]+
\PP_x(\wh X^{b, E}_{\wh\tau^{b, E}_{U}}\in D)\\
&&\qquad= 4\int_U G^b_U(y,x) \frac{h_E(y)}{h_E(x)} \,dy \\
&&\qquad\quad{}+ \int_{D\setminus
U} \int_{U} G^b_U(y,x) \frac{h_E(z)}{h_E(x)} J(y,z) \,dy \,dz.
\end{eqnarray*}
Now using Theorem
\ref{tufgfnest21}, we get
\begin{eqnarray*}
&&\PP_x(\wh\tau^{b, E}_{D}>1/4)\\
&&\qquad
\leq 4c_1M \int_U G_U(y,x) \,dy\\
&&\qquad\quad{}+c_1M \int_{D\setminus U} \int_{U}
G_U(y,x) J(y,z) \,dy \,dz \\
&&\qquad= 4c_1M \int_U G_U(x,y) \,dy
+c_1M \PP_x( X_{\tau_U} \in D\setminus\overline U) \\
&&\qquad\le c_2 \delta_U(x)^{\alpha/2} =c_2 \delta_D(x)^{\alpha/2}.
\end{eqnarray*}
The last inequality is due to~(\ref{e411}) and the boundary
Harnack principle for
$X$ in $C^{1,1}$ open sets.
\end{pf}
\begin{lemma}\label{lgen}
Suppose that $U_1,U_3, U$ are open subsets of $\bR^d$ with $U_1,
U_3\subset U \subset\frac34E$ and $\operatorname{dist}(U_1,U_3)>0$. Let
$U_{2} :=U\setminus(U_1\cup U_3)$. If $x\in U_1$ and $y \in U_3$,
then for all $t >0$,
%
%
\begin{eqnarray}\label{equb0}
p^{b}_{U}(t, x, y)
&\le& \PP_x(X^{b}_{\tau^{b}_{U_1}}\in U_{2})
\cdot\sup_{s<t, z\in U_{2}} p_U^b(s, z, y) \nonumber\\[-8pt]\\[-8pt]
&&{} +
( t \wedge\E_x [\tau^{b}_{U_1}] ) \cdot
\sup_{u\in U_1, z\in U_3}J(u,z) ,\nonumber\\
\label{equb}
p^{b}_{U}(t, y, x)
&\le& M\PP_x(\wh X^{b, E}_{\wh\tau^{b, E}_{U_1}}\in U_{2})
\cdot\sup_{s<t, z\in U_{2}}
p^{b}_U(s, y, z) \nonumber\\[-8pt]\\[-8pt]
&&{} +
M
( t \wedge\E_x [\wh\tau^{b, E}_{U_1}] ) \cdot
\sup_{u\in U_1, z\in U_3}J(u,z)\nonumber
\end{eqnarray}
and
%
%
\begin{equation}\label{eqlb}\qquad
p^{b}_{U}(1/3, x, y)\ge\frac1{3M} \PP_x(\tau^{b}_{U_1}>1/3)
\PP_y(\wh\tau^{b, E}_{U_3}>1/3) \cdot\inf_{u\in U_1,
z\in U_3}J(u,z) .
\end{equation}
\end{lemma}
\begin{pf}
The proof of~(\ref{equb0}) is similar to the proof of~\cite{BGR}, Lemma 2, which is a variation of the proof of~\cite{CKS},
Lemma 2.2.
Hence we omit its proof. We will present a proof for (\ref
{equb}) and~(\ref{eqlb}).
Using the strong Markov property and~(\ref{whp}), we have
\begin{eqnarray*}
p^{b}_{U}(t, y, x) &=& \frac{h_E(x)}{h_E(y)}\wh p^{b, E}_{U}(t, x,
y) \\
&=&\frac{h_E(x)}{h_E(y)}\E_x[\wh p^{b, E}_{U}(t-\wh
\tau^{b, E}_{U_1}, \wh X^{b, E}_{\wh\tau^{b, E}_{U_1}}, y)
; \wh\tau^{b, E}_{U_1}<t ]\\
&=&\frac{h_E(x)}{h_E(y)} \E_x[\wh p^{b, E}_{U}(t-\wh
\tau^{b, E}_{U_1}, \wh X^{b, E}_{\wh\tau^{b, E}_{U_1}},
y); \wh\tau^{b, E}_{U_1}<t, \wh X^{b, E}_{\wh\tau^{b,
E}_{U_1}}\in U_{2} ] \\
&&{} + \frac{h_E(x)}{h_E(y)}\E_x[\wh p^{b, E}_{U}(t-\wh\tau
^{b, E}_{U_1}, \wh X^{b, E}_{\wh\tau^{b, E}_{U_1}}, y); \wh
\tau^{b, E}_{U_1}<t, \wh X^{b, E}_{\wh\tau^{b, E}_{U_1}}\in
U_3] \\
&=&\!: I+\mathit{II}.
\end{eqnarray*}
Using~(\ref{whp}) again,
\begin{eqnarray*}
I &\le& \frac{h_E(x)}{h_E(y)}\PP_x(\wh\tau^{b, E}_{U_1}<t,
\wh X^{b, E}_{\wh\tau^{b, E}_{U_1}}\in U_{2}) \Bigl(
\sup_{s<t, z\in U_{2}} \wh p^{b, E}_{U}(s, z, y)\Bigr) \nonumber\\
&=&\frac{h_E(x)}{h_E(y)}\PP_x(\wh\tau^{b, E}_{U_1}<t, \wh X^{b,
E}_{\wh\tau^{b, E}_{U_1}}\in U_{2}) \biggl( \sup_{s<t, z\in
U_{2}}
p^{b}_{U}(s, y, z)\frac{h_E(y)}{h_E(z)}\biggr) \nonumber\\
&\le& \biggl(\sup_{a,b \in3E/4} \frac{h_E(a)}{h_E(b)}\biggr)
\PP_x( \wh X^{b, E}_{\wh\tau^{b, E}_{U_1}}\in
U_{2}) \Bigl( \sup_{s<t, z\in U_{2}}
p^{b}_{U}(s, y, z)\Bigr).
\end{eqnarray*}
On the other hand, by~(\ref{elevy2}) and~(\ref{whp}),
\begin{eqnarray*}
\mathit{II}
&=& \frac{h_E(x)}{h_E(y)}\int_0^{t} \int_{U_1} \wh p^{b,
E}_{{U_1}}(s, x, u) \int_{U_3} J(u,z) \frac{h_E(z)}{h_E(u)}
p^{b}_{U}(t-s, y, z) \\
&&\hspace*{136.8pt}{}\times\frac{h_E(y)}{h_E(z)}
\,dz \,du \,ds\\
&\le&\biggl(\sup_{a,b \in3E/4} \frac{h_E(a)}{h_E(b)}\biggr)
\int_0^{t} \int_{U_1} \wh p^{b, E}_{{U_1}}(s, x, u)
\int_{U_3} J(u,z) p^{b}_{U}(t-s, y, z) \,dz \,du \,ds\\
&\le& M\Bigl(\sup_{u\in U_1, z\in U_3}J(u,z)\Bigr)
\int_0^{t} \PP_x(\wh\tau^{b, E}_{U_1}>s) \biggl(\int_{U_3}
p^{b}_{U}(t-s, y,
z)\,dz \biggr) \,ds\\
&\le& M\int_0^{t}\PP_x(\wh\tau^{b, E}_{U_1}>s) \,ds \cdot\sup_{u\in
U_1,
z\in U_3}J(u,z) \\
&\le& M (t \wedge\E_x [\wh\tau^{b, E}_{U_1}])
\cdot\sup_{u\in U_1, z\in U_3}J(u,z) .
\end{eqnarray*}

Now we consider the lower bound.
By~(\ref{elevysE}) and~(\ref{whp}),
\begin{eqnarray*}
\hspace*{-2pt}&& p^{b}_{U}(1/3, x, y)\\
\hspace*{-2pt}&&\qquad\ge \E_x[ p^{b}_{U}(1/3- \tau^{b}_{U_1}, X^{b}_{
\tau^{b}_{U_1}}, y);\tau^{b}_{U_1}<1/3, X^{b}_{ \tau
^{b}_{U_1}}\in
U_3 ]\\
\hspace*{-2pt}&&\qquad=\int_0^{1/3} \biggl( \int_{U_1}
p^{b}_{{U_1}}(s, x, u) \biggl( \int_{U_3} J(u,z)
p^{b}_{U}(1/3-s, z, y)
\,dz\biggr) \,du\biggr) \,ds\\
\hspace*{-2pt}&&\qquad\ge \inf_{u\in U_1, z\in U_3}J(u,z)
\int_0^{1/3}\int_{U_3} p^{b}_{U}(1/3-s, z, y)\PP_x( \tau
^{b}_{U_1}>s)\,dz \,ds\\
\hspace*{-2pt}&&\qquad\ge \PP_x( \tau^{b}_{U_1}>1/3) \inf_{u\in U_1, z\in U_3}J(u,z)
\int_0^{1/3}\int_{U_3}p^{b}_{U_3}(1/3-s, z, y)\,dz \,ds\\
\hspace*{-2pt}&&\qquad= \PP_x( \tau^{b}_{U_1}>1/3) \inf_{u\in U_1, z\in U_3}J(u,z)
\int_0^{1/3}\int_{U_3}\wh p^{b, E}_{U_3}(1/3-s, y, z)
\frac{h_E(y)}{h_E(z)} \,dz \,ds\\
\hspace*{-2pt}&&\qquad\ge M^{-1} \PP_x( \tau^{b}_{U_1}>1/3) \inf_{u\in U_1, z\in
U_3}J(u,z)
\int_0^{1/3}\PP_y( \wh\tau^{b, E}_{U_3}
>1/3-s) \,ds\\
\hspace*{-2pt}&&\qquad\ge \frac1{3M} \PP_x( \tau^{b}_{U_1}>1/3) \inf_{u\in U_1, z\in
U_3}J(u,z) \PP_y( \wh\tau^{b, E}_{U_3} > 1/3) .
\end{eqnarray*}
\upqed
\end{pf}
\begin{lemma}\label{lemppu}
There is a positive constant $c_1 =c_1(d, \alpha, R_0, \Lambda_0, M,
b)$ with the dependence on $b$ only via the rate at which
$M^\alpha_{|b|}(r)$ goes to zero such that
for all $x, y\in D$,
%
%
\begin{equation}\label{eqppu1}
p^{b}_{D}(1/2,x,y)\leq c_1 \bigl( 1\wedge
\delta_{D}(x)^{\alpha/2}\bigr) \biggl( 1 \wedge
\frac{1}{|x-y|^{d+\alpha}}\biggr)
\end{equation}
and
%
%
\begin{equation}\label{eqppu2}
p^{b}_{D}(1/2,x,y)\leq c_1 \bigl( 1\wedge
\delta_{D}(y)^{\alpha/2}\bigr) \biggl( 1 \wedge
\frac{1}{|x-y|^{d+\alpha}}\biggr) .
\end{equation}
\end{lemma}
\begin{pf}
We only give the proof of~(\ref{eqppu2}). Recall that there
exists $L=L(R_0, \Lambda_0, d)$ such that for every $z \in\partial
D$ and $r \le R_0/2$, one can find a $C^{1,1}$ open set
$U_{(z,r)}$
with $C^{1, 1}$ characteristic $(rR_0/L, \Lambda_0L/r)$ such
that $D \cap B(z,\break r/2) \subset U_{(z,r)} \subset D \cap B(z, r) $.

It follows from~(\ref{e11a}) that
\[
p^{b}_{D}(1/2,x,y) \leq p^{b}(1/2,x,y) \leq
c_1\biggl( 1 \wedge
\frac{1}{|x-y|^{d+\alpha}}\biggr),
\]
so it suffices to prove of~(\ref{eqppu2}) for $y\in D$ with
$\delta_{D} (y) <r_0/(32)$.

When $|x-y|\le r_0 $,
by the semigroup property~(\ref{eck}),
(\ref{e10}) and~(\ref{whp}),
\begin{eqnarray}
p^{b}_{D}(1/2,x,y)&=&
\int_D p^b_D(1/4,x,z)p^b_D(1/4,z,y)\,dz \nonumber\\
&\le& \int_D p^b(1/4,x,z)\wh p^{b, E}_{D}(1/4,y, z)
\frac{h_E(y)}{h_E(z)}\,dz \nonumber\\
&\leq&c_2 M \int_D\biggl( 1 \wedge\frac{1}{|x-z|^{d+\alpha
}}\biggr)
\wh p^{b, E}_{D}(1/4,y, z) \,dz \nonumber\\
&\leq& c_2 M \PP_y(\wh\tau^{b, E}_{D}>1/4).\nonumber
\end{eqnarray}
Applying~(\ref{eqcw}), we get
\begin{eqnarray*}
p^{b}_{D}(1/2,x,y) & \leq&c_3 \bigl( 1\wedge
\delta_{D}(y)^{\alpha/2}\bigr)\\
& \le& c_3 ( 1 \vee r_0^{d+\alpha})
\bigl( 1\wedge\delta_{D}(y)^{\alpha/2}\bigr) \biggl( 1
\wedge\frac{1}{|x-y|^{d+\alpha}}\biggr).
\end{eqnarray*}

Finally we consider the case that $|x-y|> r_0$ [and $\delta_{D} (y)
<r_0/(32)$]. Fix $y \in D$ with $\delta_{D} (y)<r_0/(32)$, and let
$Q\in
\partial D$ be such that $|y-Q|=\delta_{D} (y)$. Choose a $C^{1,1}$
open set $U_y:=U_{(Q, r_0/8)}$ with $C^{1, 1}$ characteristic $(r_0
R_0/(8L),\break 8\Lambda_0L/r_0)$ such that $D \cap B(Q, r_0/(16)) \subset U_y
\subset D \cap B(Q, r_0/8 ) $.

Let $D_3:= \{z\in D\dvtx |z-y|>|x-y|/2\}$ and $D_2:=D\setminus(U_y\cup
D_3)$. Note that $|z-y| >r_0/2 $ for $z\in D_3$. So, if
$u\in U_y$ and $z\in D_3$, then
\[
|u-z| \ge|z-y|-|y-u| \ge|z-y|-r_0/4 \ge
\tfrac{1}{2}|z-y| \ge\tfrac{1}{4}|x-y|.
\]
Thus
%
%
\begin{equation} \label{en0}\quad
\sup_{u\in U_y, z\in D_3}J(u,z) \le\sup_{(u,z)\dvtx|u-z| \ge
|x-y|/4}J(u,z) \le c_4 \biggl( 1 \wedge
\frac{1}{|x-y|^{d+\alpha}}\biggr).
\end{equation}
If $z \in D_2$, then $|z-x| \ge|x-y| -|y-z| \ge|x-y|/2$. Thus by
(\ref{e10}),
%
%
\begin{eqnarray} \label{en1}
\sup_{s<1/2, z\in D_2} p_D^b(s, x, z)
&\le&
\sup_{s<1/2, z\in D_2}p^b(s, x, z)\nonumber \\
&\le&
c_5\sup_{s<1/2, z \in
D_2} \biggl( 1
\wedge\frac{1}{|x-z|^{d+\alpha}}\biggr)\\
&\le& c_6 \biggl( 1
\wedge\frac{1}{|x-y|^{d+\alpha}}\biggr)
\nonumber
\end{eqnarray}
for some $c_5, c_6>0$.
Applying Lemma~\ref{lgen} with~(\ref{en0}) and~(\ref{en1}), we
obtain
\[
p^{b}_{D}(1/2, x, y) \le c_7 \biggl( 1 \wedge
\frac{1}{|x-y|^{d+\alpha}}\biggr) \bigl(\PP_y (\wh X^{b,
E}_{\wh
\tau^{b, E}_{U_y}}\in D ) + \E_y [\wh\tau^{b,
E}_{U_y}]
\bigr).
\]
On the other hand, by~(\ref{edGk}),~(\ref{eKK1}), Lemma
\ref{L20} and Theorem
\ref{tufgfnest21},
\begin{eqnarray*}
&&\E_y[\wh\tau^{b, E}_{U_y}]+ \PP_y(\wh X^{b,
E}_{\wh\tau^{b, E}_{U_y}}\in D)\\
&&\qquad= \int_{U_y} G^b_{U_y}(z,y)
\frac{h_E(z)}{h_E(y)} \,dz + \int_{D\setminus{U_y}} \int_{{U_y}}
G^b_{U_y}(w,y) \frac{h_E(z)}{h_E(y)} J(w,z) \,dw \,dz\\
&&\qquad \leq c_8M \int_{U_y} G_{U_y}(z,y) \,dz +c_8M \int_{D\setminus
{U_y}} \int_{{U_y}}
G_{U_y}(w,y) J(w,z) \,dw \,dz \\
&&\qquad\le c_9 \delta_{U_y}(y)^{\alpha/2} =c_9 \delta_D(y)^{\alpha/2}.
\end{eqnarray*}
Therefore
\[
p^{b}_{D}(1/2, x, y) \le c_{10}
\delta_D (y)^{\alpha/2}
\biggl( 1 \wedge\frac{1}{|x-y|^{d+\alpha}}\biggr) .
\]
Equation~(\ref{eqppu1}) can be proved in a similar way.
\end{pf}
\begin{lemma}\label{lemppu1}
There is a positive constant $c_1 =c_1(d, \alpha, R_0, \Lambda_0, M,
b)$ with the dependence on $b$ only via the rate at which
$M^\alpha_{|b|}(r)$ goes to zero such that
for all $x, y\in D$,
%
%
\begin{equation}\label{eqppu}\quad
p^{b}_{D}(1,x,y)\leq c_1 \bigl( 1\wedge
\delta_{D}(x)^{\alpha/2}\bigr) \bigl( 1\wedge
\delta_{D}(y)^{\alpha/2}\bigr) \biggl( 1 \wedge
\frac{1}{|x-y|^{d+\alpha}}\biggr) .
\end{equation}
\end{lemma}
\begin{pf}
Using~(\ref{eqppu1}) and~(\ref{eqppu2}), the semigroup
property~(\ref{eck}) and the two-sided estimates of $p(t, x, y)$,
\begin{eqnarray*}
p^{b}_D(1,x,y) &=& \int_{\bR^d} p^{b}_D(1/2,x,z)p^{b}_D(1/2,z,y)\,dz\\
&\leq& c\bigl( 1\wedge\delta_{D}(x)^{\alpha/2}\bigr)\bigl(
1\wedge
\delta_{D}(y)^{\alpha/2}\bigr)\\
&&{} \times\int_{\bR^d} \biggl( 1 \wedge
\frac{1}{|x-z|^{d+\alpha}}\biggr) \biggl( 1
\wedge\frac{1}{|z-y|^{d+\alpha}}\biggr)\,dz\\
&\leq& c\bigl( 1\wedge\delta_{D}(x)^{\alpha/2}\bigr)\bigl(
1\wedge
\delta_{D}(y)^{\alpha/2}\bigr)\int_{\bR^d}p(1/2, x, z)p(1/2, z,
y)\,dz\\
&=&c\bigl( 1\wedge\delta_{D}(x)^{\alpha/2}\bigr) p(1, x, y)\\
&\leq&c\bigl( 1\wedge\delta_{D}(x)^{\alpha/2}\bigr)\bigl(
1\wedge
\delta_{D}(y)^{\alpha/2}\bigr)\biggl( 1 \wedge
\frac{1}{|x-y|^{d+\alpha}}\biggr) .
\end{eqnarray*}
\upqed
\end{pf}
\begin{lemma}\label{lowerbound12}
If $r>0$, then there is a constant $c_1
=c_1(d, \alpha, r, M, b)>0$
with the dependence on $b$ only via the rate at
which $M^\alpha_{|b|}(r)$ goes to zero such that for every $B(u,r),
B(v,r) \subset\frac{3}{4}E$,
\[
p^{b}_{B(u,r)\cup B(v,r)}(1/3, u, v) \ge
c_1 \biggl( 1
\wedge\frac{1}{|u-v|^{d+\alpha}}\biggr).
\]
\end{lemma}
\begin{pf}
If $|u-v|\le r/2$, by Proposition~\ref{P35-55},
\begin{eqnarray*}
p^{b}_{B(u,r)\cup B(v,r)}(1/3, u, v) &\ge& \inf_{ |u-v|<r/2}
p^{b}_{B(u,r)}(1/3, u, v) \\&\ge& c_1 \ge c_2 \biggl( 1
\wedge\frac{1}{|u-v|^{d+\alpha}}\biggr).
\end{eqnarray*}

If $|u-v|\ge r/2$, with $U_1= B(u,r/8)$ and $U_3=B(v,r/8)$, we have,
by~(\ref{eqlb}),
\begin{eqnarray*}
&& p^{b}_{B(u,r)\cup B(v,r)}(1/3, u, v) \\
&&\qquad\ge \frac13
\PP_u(\tau^{b}_{U_1}>1/3)
\PP_v(\wh\tau^{b, E}_{U_3}>1/3)\inf_{w\in U_1, z\in U_3}J(w,z)\\
&&\qquad \ge c \int_{ B(u,r/16)} p^{b}_{B(u,r/8)}(1/3, u, z)\,dz \int_{
B(v,r/16)} \wh p^{b, E}_{B(u,r/8)}(1/3, v, z)\,dz \\
&&\qquad\quad{} \times\biggl( 1
\wedge\frac{1}{|u-v|^{d+\alpha}}\biggr) \\[-2pt]
&&\qquad\ge c \Bigl(\inf_{z \in B(u,r/16)} p^{b}_{B(u,r/8)}(1/3, u,
z)\Bigr)
\Bigl(\inf_{z \in B(v,r/16)} \wh p^{b, E}_{B(u,r/8)}(1/3, v,
z)\Bigr)\\[-2pt]
&&\qquad\quad{} \times\biggl( 1 \wedge\frac{1}{|u-v|^{d+\alpha}}\biggr).
\end{eqnarray*}
Now applying Propositions~\ref{P35-55} and~\ref{P35-55-1}, we
conclude that
\[
p^{b}_{B(u,r)\cup B(v,r)}(1/3, u, v) \ge c \biggl( 1
\wedge\frac{1}{|u-v|^{d+\alpha}}\biggr).
\]
\upqed
\end{pf}
\begin{lemma}\label{lemppl}
There is a positive constant $c_1 =c_1(d, \alpha, R_0, \Lambda_0, M,
b)$ with the dependence on $b$ only via the rate at which
$M^\alpha_{|b|}(r)$ goes to zero such that
\[
p^{b}_{D}(1,x,y)\geq c_1\bigl( 1\wedge
\delta_{D}(x)^{\alpha/2}\bigr) \bigl( 1\wedge
\delta_{D}(y)^{\alpha/2}\bigr) \biggl( 1 \wedge
\frac{1}{|x-y|^{d+\alpha}}\biggr) .
\]
\end{lemma}
\begin{pf}
Recall that $r_0 \le R_0/8$ is the constant from Theorem
\ref{tufgfnest21}
which depends only on $d$, $\alpha$,
$R_0$, $\Lambda_0$ and
$b$ with the dependence on $b$
only via the rate at which $M^\alpha_{|b|}(r)$
goes to zero. Since $D$ is $C^{1,1}$ with $C^{1,1}$ characteristics
$(R_0, \Lambda_0)$, there exist $\delta=\delta(R_0, \Lambda_0) \in
(0, r_0/8)$ and $L=L(R_0, \Lambda_0)>1$
so that for all \mbox{$x,y \in D$},
there are
$\xi_x \in D\cap B(x, L\delta)$ and $\xi_y \in D
\cap B(y, L\delta)$ with $ B(\xi_x, 2\delta) \!\cap\! B(x, 2\delta)=
\varnothing$, $B(\xi_y, 2\delta) \!\cap\! B(y, 2\delta)= \varnothing$ and
\mbox{$B(\xi_x, 8\delta) \!\cup\! B(\xi_y, 8\delta) \subset D$}.

Note that by the semigroup property~(\ref{eck}) and Lemma~\ref{lowerbound12},
%
%
\begin{eqnarray}\label{eloww21}
&&p^{b}_{D}(1,x,y)\nonumber\\[-2pt]
&&\qquad\geq \int_{B(\xi_y, \delta)}\int_{B(\xi_x, \delta)}
p^{b}_{D}(1/3,x,u)
p^{b}_{D}(1/3,u,v)p^{b}_{D}
(1/3,v,y)\,du\,dv \nonumber\\[-2pt]
&&\qquad\geq \int_{B(\xi_y, \delta)}\int_{B(\xi_x, \delta)}
p^{b}_{{D}}(1/3,x,u)p^{b}_{B(u, \delta/2) \cup
B(v,\delta/2)}(1/3,u,v)\nonumber\\[-2pt]
&&\qquad\quad\hspace*{64.5pt}{}\times p^{b}_
{D}(1/3,v,y)\,du\,dv\\[-2pt]
&&\qquad\geq c_1\int_{B(\xi_y, \delta)}\int_{B(\xi_x, \delta)}
p^{b}_{{D}}(1/3,x,u)\bigl(J(u,v)\wedge
1\bigr)p^{b}_{D}(1/3,v,y)\,du\,dv\nonumber\\[-2pt]
&&\qquad\geq c_1 \Bigl(\inf_{(u,v) \in B(\xi_x, \delta) \times
B(\xi_y,\delta)} \bigl(J(u,v)\wedge1\bigr) \Bigr) \nonumber\\[-2pt]
&&\qquad\quad{} \times
\biggl( \int_{B(\xi_x, \delta)}
p^{b}_{{D}}(1/3,x,u)\,du\biggr)
\biggl(\int_{B(\xi_y, \delta)}
p^{b}_{D}(1/3,v,y)\,dv\biggr)
.\nonumber
\end{eqnarray}
If $|x-y| \ge\delta/8$,
$|u-v| \le2(1+L)\delta+|x-y| \le(17 +16L)|x-y|$, and
we have
%
%
\begin{equation}\label{eloww2}
\inf_{(u,v) \in B(\xi_x, \delta) \times B(\xi_y,\delta)}
\bigl(J(u,v)\wedge1\bigr)\ge c_2\biggl( 1
\wedge\frac{1}{|x-y|^{d+\alpha}}\biggr).\vadjust{\goodbreak}
\end{equation}
If $|x-y| \le\delta/8$,
$|u-v| \le2(2+L)\delta$ and
%
%
\begin{equation}\label{eloww21p}
\inf_{(u,v) \in B(\xi_x, \delta) \times B(\xi_y,\delta)}
\bigl(J(u,v)\wedge1\bigr)\ge c_3 \ge c_4\biggl( 1
\wedge\frac{1}{|x-y|^{d+\alpha}}\biggr).
\end{equation}

We claim that
%
%
\begin{equation}\label{eloww31}
\int_{B(\xi_x, \delta)}
p^{b}_{{D}}(1/3,x,u)\,du \ge c_5 \bigl( 1\wedge
\delta_{D}(x)^{\alpha/2}\bigr)
\end{equation}
and
%
%
\begin{equation}\label{eloww3}
\int_{B(\xi_y, \delta)}
p^{b}_{D}(1/3,v,y)\,dv \ge c_5 \bigl( 1\wedge
\delta_{D}(y)^{\alpha/2}\bigr),
\end{equation}
which, combined with~(\ref{eloww21})--(\ref{eloww21p}), proves the
theorem.

We only give the proof of
(\ref{eloww3}).
If $\delta_{D}(y)> \delta$, since $\operatorname{dist}(B(\xi_y, \delta), B(y,
\delta))>0$,
by~(\ref{eqlb}),
%
%
\begin{eqnarray} \label{eddd0}
&&\int_{B(\xi_y, \delta)} p^{b}_{D}(1/3,v,y)\,dv\nonumber\\
&&\qquad\ge \frac{1}{3M}
\biggl(\int_{B(\xi_y, \delta)} \PP_v\bigl(\tau^{b}_{B(\xi_y,
\delta)}>1/3\bigr)\,dv \biggr) \PP_y\bigl(\wh\tau^{b, E}_{B(y,
\delta)}>1/3\bigr)\\
&&\qquad\quad{} \times\inf_{w\in B(\xi_y, \delta), z\in{B(y,
\delta)}}J(w,y) \nonumber,
\end{eqnarray}
which is greater than or equal to some positive constant depending
only on $d, \alpha, R_0,
\Lambda_0, $
$M$ and $b$, with the dependence on
$b$ only via the rate at which $M^\alpha_{|b|}(r)$ goes to zero
by Propositions~\ref{P35-55} and~\ref{P35-55-1}.

If $\delta_{D}(y) \le\delta$, choose a $Q\in
\partial D$ be such that $|y-Q|=\delta_{D} (y)$, and
choose a $C^{1,1}$ open set
$U_y:=U_{(Q,4\delta)}$
with $C^{1, 1}$ characteristic
$(4\delta R_0/L, \Lambda_0L/(4\delta))$
such that
\[
D \cap B(Q, 2\delta)
\subset U_y \subset D \cap B(Q,4\delta) \subset D \cap B(Q, 6\delta)
=:V_y.
\]
Then,
since $\operatorname{dist}(B(\xi_y, \delta), V_y)>0$,
by~(\ref{eqlb}),
%
%
\begin{eqnarray}\label{eddd3}
&&\int_{B(\xi_y, \delta)} p^{b}_{D}(1/3,v,y)\,dv\nonumber\\
&&\qquad\ge \frac{1}{3M} \biggl(\int_{B(\xi_y, \delta)}
\PP_v\bigl(\tau^{b}_{B(\xi_y, \delta)}>1/3\bigr)\,dv \biggr) \PP
_y(\wh\tau^{b,
E}_{V_y}>1/3)\\
&&\qquad\quad{} \times\inf_{w\in B(\xi_y, \delta), z\in
V_y}J(w,y)\nonumber,
\end{eqnarray}
which is greater than or equal to $c_6 \PP_y(\wh\tau^{b,
E}_{V_y}>1/3)$ for some positive constant $c_6$ depending
only on $d, \alpha, R_0, \Lambda_0$,
$M$ and
$b$ with the dependence on
$b$ only via the
rate at which $M^\alpha_{|b|}(r)$ goes to zero by Propositions
\ref{P35-55} and~\ref{P35-55-1}.\vspace*{1pt}\vadjust{\goodbreak}

Let $B(y_0, 2c_7
\delta)$ be a ball in $ D \cap(B(Q,
6\delta) \setminus
B(Q,
4\delta))$ where $c_7=c_7(\Lambda_0$, $
d)>0$. By the strong Markov property,
\begin{eqnarray*}
&&\Bigl(\inf_{w \in B(y_0, c_7
\delta/2)} \PP_{w}\bigl(\wh\tau^{b, E}_{B(w, c_7
\delta)}
>1/3\bigr)\Bigr) \PP_y\bigl(
\wh X^{b, E}_{\wh\tau^{b, E}_{U_y}} \in B(y_0, c_7
\delta/2)\bigr)\\[-0.4pt]
&&\qquad\le \E_y\bigl[ \PP_{\wh X^{b, E}_{\wh\tau^{b, E}_{U_y}}}
\bigl(\wh\tau^{b, E}_{B(\wh X^{b, E}_{\wh\tau^{b, E}_{U_y}},c_7
\delta)}
>1/3\bigr); \wh X^{b, E}_{\wh\tau^{b, E}_{U_y}} \in B(y_0, c_7
\delta/2)\bigr]\\[-0.4pt]
&&\qquad\le \E_y[ \PP_{\wh X^{b, E}_{\wh\tau^{b, E}_{U_y}}}
(\wh\tau^{b, E}_{V_y}
>1/3); \wh X^{b, E}_{\wh\tau^{b, E}_{U_y}}
\in B(y_0, c_7
\delta/2)]\\[-0.4pt]
&&\qquad\le \PP_y\bigl( \wh\tau^{b, E}_{V_y}
>1/3, \wh X^{b, E}_{\wh\tau^{b, E}_{U_y}} \in B(y_0, c_7
\delta/2)\bigr) \le\PP_y( \wh\tau^{b, E}_{V_y}
>1/3 ).
\end{eqnarray*}
Using Proposition
\ref{P35-55-1}, we get
%
%
\begin{equation}
\PP_y( \wh\tau^{b, E}_{V_y}
>1/3) \ge c_8 \PP_y\bigl( \wh
X^{b, E}_{\wh\tau^{b, E}_{U_y}} \in B(y_0, c_7
\delta/2)\bigr).
\end{equation}
Now applying~(\ref{edGk}),~(\ref{eKK1}) and Theorem
\ref{tufgfnest21},
%
%
\begin{eqnarray}\label{eddd1}
&&\PP_y\bigl( \wh X^{b, E}_{\wh\tau^{b, E}_{U_y}} \in B(y_0, c_7
\delta/2)\bigr)\nonumber\\[-0.4pt]
&&\qquad= \int_{B(y_0, c_7
\delta/2)} \int_{{U_y}}
G^b_{U_y}(w,y) \frac{h_E(z)}{h_E(y)} J(w,z) \,dw \,dz
\nonumber\\[-8pt]\\[-8pt]
&&\qquad \ge c_{9}M^{-1} \int_{B(y_0, c_7
\delta/2)} \int_{{U_y}}
G_{U_y}(w,y) J(w,z) \,dw \,dz \nonumber\\[-0.4pt]
&&\qquad\ge c_{10} \delta_{U_y}(y)^{\alpha/2} =c_{10} \delta_D(y)^{\alpha
/2}.\nonumber
\end{eqnarray}
Combining~(\ref{eddd0})--(\ref{eddd1}), we have proved
(\ref{eloww3}).
\end{pf}
\begin{theorem}\label{thmoppz}
There exists $c=c(d, \alpha, R_0, \Lambda_0, T, M,b)>0$
with the dependence on $b$ only via the rate at which
$M^\alpha_{|b|}(r)$ goes to zero such that
for $0<t\leq T$, $x,y\in D$,
%
%
\begin{eqnarray}\label{eqget}
&& c^{-1} \biggl( 1\wedge\frac{\delta_D(x)^{\alpha/2}}{\sqrt
{t}}\biggr)
\biggl( 1\wedge\frac{\delta_D(y)^{\alpha/2}}{\sqrt{t}}\biggr)
\biggl(
t^{-d/\alpha} \wedge\frac{t}{|x-y|^{d+\alpha}}\biggr)
\nonumber\\
&&\qquad\leq
{p^b_{D}(t, x, y)} \\
&&\qquad\le
c \biggl( 1\wedge
\frac{\delta_D(x)^{\alpha/2}}{\sqrt{t}}\biggr) \biggl( 1\wedge
\frac{\delta_D(y)^{\alpha/2}}{\sqrt{t}}\biggr) \biggl(
t^{-d/\alpha}
\wedge\frac{t}{|x-y|^{d+\alpha}}\biggr). \nonumber
\end{eqnarray}
\end{theorem}
\begin{pf}
Let $D_{t}:={t}^{-1/\alpha}D$ and $E_{t}:={t}^{-1/\alpha}E$.
By the scaling property in~(\ref{escaling}),
(\ref{eqget}) is equivalent to
\begin{eqnarray*}
&& c^{-1} \bigl( 1\wedge\delta_{D_t}(x)^{\alpha/2}\bigr) \bigl(
1\wedge\delta_{D_t}(y)^{\alpha/2}\bigr) \biggl( 1 \wedge
\frac{1}{|x-y|^{d+\alpha}}\biggr)
\\
&&\qquad\leq
p^{t^{(\alpha-1)/\alpha} b({t}^{1/\alpha} \cdot)}_{D_{t}}(1,x,y)
\\
&&\qquad\leq
c \bigl( 1\wedge\delta_{D_t}(x)^{\alpha/2}\bigr) \bigl( 1\wedge
\delta_{D_t}(y)^{\alpha/2}\bigr) \biggl( 1 \wedge
\frac{1}{|x-y|^{d+\alpha}}\biggr).
\end{eqnarray*}
The above holds in view of
(\ref{e72}),~(\ref{eMgoOK}),~(\ref{e73})
and the fact that for $t \le T$, the $D_{t}$'s are $C^{1,1}$ open sets
in $\bR^d$ with the same $C^{1,1}$ characteristics $(R_0
(T)^{-1/\alpha},
\Lambda_0(T)^{-1/\alpha})$. The theorem is thus proved.
\end{pf}

\section{Large time heat kernel estimates}

Recall that we have fixed a ball $E$ centered at the origin, and
$M>1$ is the constant in~(\ref{MMM}). Let $U$ be an arbitrary open
set $U \subset\frac14E$, and we let
\[
\overline{p}^{b, E}_U(t,x,y):= \frac{ p^b_U(t,x,y)}{h_E(y)},
\]
which is strictly positive, bounded and continuous on \mbox{$(t,x,y) \in
(0, \infty) \times U \times U$} because $p^b_U(t,x,y)$ is strictly
positive, bounded and continuous on $(t,x,y) \in(0, \infty) \times
U \times U$, and $h_E(y)$ is strictly positive and continuous on $E$.
For each $x \in U$, $(t,y) \mapsto\overline{p}^{b, E}_U(t,x,y)$ is
the transition density of $(X^{b,U}, \PP_x)$ with respect to the
reference measure $\xi_E$, and, for each $y \in U$, $(t,x) \mapsto
\overline{p}^{b, E}_U(t,x,y)$ is the transition density of $(\wh
X^{b,E, U}, \PP_y)$, the dual process of $X^{b,U}$ with respect to
the reference measure $\xi_E$.

Let
\[
P^{b, E, U}_tf(x):=\int_U\overline{p}^{b, E}_U(t,x,y) f(y) \xi_E(dy)
\]
and
\[
\wh{P}^{b, E, U}_tf(x):=\int_U\overline{p}^{b,
E}_U(t,y,x) f(y) \xi_E(dy).
\]
Let $\sL^{b, E}_U$ and $\wh{\sL}^{b, E}_U$ be the infinitesimal
generators of the semigroups $\{P^{b, E, U}_t\}$ and $\{\wh{P}^{b,
E, U}_t\}$ on $L^2(U, \xi_E)$, respectively.

Note that, since for each $t>0$, $\overline{p}^{b, E}_U(t,x,y)$ is
bounded in $U\times U$, it follows from Jentzsch's theorem
(\cite{Sc}, Theorem V.6.6, page 337) that the common value
$-\lambda^{b, E, U}_0:= \sup\operatorname{Re} (\sigma(\sL^{b, E}_U))= \sup
\operatorname{Re} (\sigma(\wh{\sL}^{b, E}_U))$ is an eigenvalue of
multiplicity~1 for both $\sL^{b, E}_U$ and $\wh{\sL}^{b, E}_U$, and that an
eigenfunction $\phi^{b, E}_U$ of $\sL^{b, E}_U$ associated\vspace*{1pt} with
$\lambda_0^{b, E, U}$ can be chosen to be strictly positive with
$\|\phi^{b, E}_U\|_{L^2(U, \xi_E(dx))}=1$, and an eigenfunction
$\psi^{b, E}_U$ of $\wh\sL^{b, E}_U$ associated\vspace*{1pt} with $\lambda_0^{b,
E, U}$ can be chosen to be strictly positive with $\|\psi^{b,
E}_U\|_{L^2(U, \xi_E(dx))}=1$.\vspace*{1pt}

It is clear from the definition that, for any Borel function $f$,
\[
P^{b, E, U}_tf(x)=P^{b, U}_tf(x) \qquad\mbox{for every } x\in U
\mbox{ and } t>0.
\]
Thus the operators $\sL^b |_U$ and
$\sL^{b, E}_U$ have the same eigenvalues. In particular, the
eigenvalue $\lambda_0^{b, E, U}$ does not depend on $E$, and so from
from now on we will denote it by $\lambda_0^{b, U}$.\vspace*{-2pt}
\begin{defn}
The semigroups $\{P^{b, E, U}_t\}$ and $\{\wh{P}^{b, E, U}_t\}$ are
said to be intrinsically ultracontractive if, for any $t>0$, there
exists a constant $c_t>0$ such that
\[
\overline{p}^{b, E}_U(t,x,y)\le c_t\phi^{b, E}_U(x)\psi^{b, E}_U(y)
\qquad\mbox{for } x, y\in U .\vspace*{-2pt}
\]
\end{defn}

It follows from~\cite{KS2}, Theorem 2.5, that if $\{P^{b, E, U}_t\}$
and $\{\wh{P}^{b, E, U}_t\}$ are intrinsically ultracontractive, then
for any $t>0$ there exists a positive constant $c_t>1$ such that
%
%
\begin{equation}\label{e12lb}
\overline{p}^{b, E}_U(t,x,y) \geq c_t^{-1}
\phi^{b, E}_U (x) \psi^{b, E}_U (y) \qquad\mbox{for } x, y\in U.\vspace*{-2pt}
\end{equation}

\begin{theorem}\label{tiu0}
For every $B(x_0, 2r) \subset U$ there exists a constant
$c=c(d, \alpha, r$, $\operatorname{diam}(U), M)>0$ such that for every $x\in
D$,
%
%
\begin{equation}\label{eIU1} \E_x \biggl[
\int^{\tau^b_U}_0{\mathbf1}_{B(x_0,r)}(X^{b, U}_t)\,dt \biggr]
\ge c \E_x [ \tau^b_U ]
\end{equation}
and
%
%
\begin{equation}\label{eIU2} \E_x \biggl[ \int^{\wh\tau^{b,
E}_U}_0{\mathbf1}_{B(x_0,r)}(\wh X^{b, E, U}_t)\,dt \biggr]
\ge c \E_x [ \wh\tau^{b, E}_U].\vspace*{-2pt}
\end{equation}
\end{theorem}
\begin{pf}
The method of the proof to be given below is now well known; see
\cite{CKS1,Ku2}. For the reader's convenience, we present the
details here. We give the proof of~(\ref{eIU2}) only. The proof for
(\ref{eIU1}) is similar. Fix a ball $B(x_0, 2r)\subset U$ and put
\[
B_0:=B(x_0, r/4),\qquad
K_1:=\overline{B(x_0, r/2)} \quad\mbox{and}\quad
B_2:=B(x_0, r).
\]
Let $\{\theta_t, t>0\}$ be the shift operators of $\wh X^{b, E}$, and
we define stopping times $S_n$ and $T_n$ recursively by
\begin{eqnarray*}
S_1(\omega) &:=& 0,\\
T_n(\omega)&:=& S_n(\omega) + \wh\tau^{b, E}_{U\setminus
K_1} \circ
\theta_{S_n}(\omega) \qquad\mbox{for }S_n(\omega) < \wh\tau^{b, E}_U
\end{eqnarray*}
and
\[
S_{n+1}(\omega):= T_{n}(\omega) +
\wh\tau^{b, E}_{B_2} \circ\theta_{T_{n}}(\omega) \qquad\mbox{for
}T_n(\omega) < \wh\tau^{b, E}_U.
\]
Clearly\vspace*{1pt} $S_n \le\wh\tau^{b, E}_U$. Let $S:=\lim_{n \to\infty} S_n
\le\wh\tau^{b, E}_U$. On $\{S< \wh\tau^{b, E}_U\}$, we must have
$S_n<T_n<S_{n+1}$ for every $n\ge0$. Using the fact that $\PP_x (
\wh\tau^{b, E}_U <\infty)=1$ for every $x\in U$ and the quasi-left
continuity of $\wh X^{b, E, U}$, we have $\PP_x( S< \wh\tau^{b,
E}_U)=0$. Therefore, for every $x \in U$,
%
%
\begin{equation}\label{eiuSTZ}
\PP_x\Bigl(\lim_{n
\to\infty}S_n=\lim_{n \to\infty} T_n = \wh\tau^{b, E}_U\Bigr)
= 1.\vadjust{\goodbreak}
\end{equation}

For any $x\in
K_1$, by Proposition~\ref{P35-55-1} we have
\[
\E_x [\wh\tau^{b, E}_{B_2}] \ge c_0 \int_{B(x_0, r/2)}
\int_{r^\alpha}^{2 r^\alpha} \wh p^{b, E}_{B_2} (t,x,y)\,dt\,dy \ge c_1
\qquad\mbox{for every } x\in K_1.
\]
Now it follows from the strong Markov property that
\begin{eqnarray*}
\E_x [ S_{n+1}-T_n ] &=& \E_x \bigl[ \E_{\wh X^{b,
E, U}_{T_n} }[
\wh\tau^{b, E}_{B_2}] ; T_n<\wh\tau^{b, E}_U \bigr]\\[-2pt]
&\geq& c_1\PP_x(\wh X^{b, E, U}_{T_n}\in B_0) \\[-2pt]
&=& c_1\E
_x [ \PP_{\wh X^{b, E, U}_{S_n}}
( \wh X^{b, E, U}_{\wh\tau^{b, E}_{U\setminus K_1}} \in B_0
)].
\end{eqnarray*}
Note that for any $x\in U\setminus B_2$,
by~(\ref{elevyp}), we have
\begin{eqnarray*}
&&\PP_x ( \wh X^{b, E, U}_{\wh\tau^{b, E}_{U\setminus K_1}}\in
B_0 )\\[-2pt]
&&\qquad=
\int_{U\setminus K_1}\frac{G^b_{U\setminus
K_1}(y,x)}{h_E(x)}\int_{B_0}
\biggl( \frac{J(y,z)h_E(z)}{h_E(y)}\,dz \biggr) \xi_E(dy)\\[-2pt]
&&\qquad\ge M^{-1}\mathcal{A}(d, -\alpha)\int_{U\setminus
K_1}\frac{G^b_{U\setminus
K_1}(y,x)}{h_E(x)}\int_{B_0}
\biggl(
\frac{dz}{(\operatorname{diam}(U))^{d+\alpha}} \biggr) \xi_E(dy)\\[-2pt]
&&\qquad= c_2\E_x [\wh\tau^{b, E}_{U\setminus
K_1}]
\end{eqnarray*}
for some constant $c_2=c_2(\alpha, r, \operatorname{diam}(U), M)>0$.
It follows
then
%
%
\begin{equation}\label{eiu2}
\E_x [ S_{n+1}-T_n ] \geq c_1 c_2 \E_x [
\E_{\wh X^{b, E, U}_{S_n}} [\wh\tau^{b, E}_{U\setminus
K_1}] ] = c_1c_2\E_x
[T_n-S_n ].
\end{equation}
Since $\wh X^{b, E, U}_t\in B_2$ for $T_n<t<S_{n+1}$, we have, by
(\ref{eiuSTZ}),
\begin{eqnarray*}
&&\E_x \biggl[ \int^{\wh\tau^{b, E}_U}_0{\mathbf1}_{B_2}(\wh
X^{b, E, U}_t)\,dt \biggr]\\[-2pt]
&&\qquad=\E_x\Biggl[
\sum^{\infty}_{n=1}
\biggl(\int^{T_n}_{S_n}{\mathbf1}_{B_2}(\wh X^{b, E, U}_t
)\,dt+\int^{S_{n+1}}_{T_n}{\mathbf1}_{B_2}
(\wh X^{b, E, U}_t)\,dt\biggr)\Biggr]\\[-2pt]
&&\qquad\ge\E_x\Biggl[ \sum^{\infty}_{n=1}\biggl(\int
^{S_{n+1}}_{T_n}{\mathbf1}_{B_2}
(\wh X^{b, E, U}_t)\,dt\biggr)\Biggr]
= \E_x\Biggl[\sum^{\infty}_{n=1}(S_{n+1}-T_n)\Biggr].
\end{eqnarray*}
\mbox{Using~(\ref{eiuSTZ}) and~(\ref{eiu2}) and noting that
$\wh X^{b, E, U}_t\notin U \setminus B_2$
for $t\in[T_n,
S_{n+1})$, we~get}
\begin{eqnarray*}
&&\E_x \biggl[ \int^{\wh\tau^{b, E}_U}_0{\mathbf1}_{B_2}(\wh
X^{b, E, U}_t)\,dt \biggr] \\[-2pt]
&&\qquad\ge c_1c_2\E_x\Biggl[
\sum^{\infty}_{n=1}(T_{n}-S_n)\Biggr]\\[-2pt]
&&\qquad\ge c_1c_2\E_x\Biggl[\sum^{\infty}_{n=1}\biggl(\int
^{T_n}_{S_n}{\mathbf1}_{U\setminus
B_2}(\wh X^{b, E, U}_t)\,dt+\int^{S_{n+1}}_{T_n}{\mathbf1}_{U\setminus B_2}
(\wh X^{b, E, U}_t)\,dt\biggr)\Biggr]\\[-2pt]
&&\qquad=c_1c_2\E_x\biggl[\int^{\wh\tau^{b, E}_U}_0{\mathbf1}_{U\setminus B_2}
(\wh X^{b, E, U}_t)\,dt\biggr].
\end{eqnarray*}
Thus
\[
\E_x \biggl[ \int^{\wh\tau^{b, E}_U}_0{\mathbf1}_{B_2}(\wh X^{b,
E, U}_t)\,dt \biggr]
\ge\frac{c_1 c_2}{1+c_1 c_2} \E_x [ \wh\tau^{b, E}_U
].
\]
\upqed
\end{pf}
\begin{theorem}\label{tiu}
$\{P^{b, E, U}_t\}$ and $\{\wh{P}^{b, E, U}_t\}$ are intrinsically
ultracontractive.
\end{theorem}
\begin{pf}
Since $\psi^{b, E}_U = e^{\lambda_0^{b, U}} \wh P^{b, E,U}_1
\psi^{b, E}_U $, it follows that $\psi^{b, E}_U$ is strictly
positive, bounded and continuous in $U$. Theorem~\ref{tiu0} implies that
%
%
\begin{eqnarray}\label{eiiu}
\E_x [ \wh\tau^{b, E}_U] &\le& c_1 \int_{B_2} \frac{G^{b,
E}_U(z,y)}{h_E(y)}\psi^{b, E}_U(z) \xi_E(dz)\nonumber\\[-8pt]\\[-8pt]
&\le& c_1 \int_{U}
\frac{G^{b, E}_U(z,y)}{h_E(y)}\psi^{b, E}_U(z) \xi_E(dz)
=\frac{c_1}{\lambda_0^{b, U}} \psi^{b, E}_U (y).
\nonumber
\end{eqnarray}
Similarly,
%
%
\begin{equation}\label{eiiu-1}
\E_x [ \tau^{b}_U]\le\frac{c_2}{\lambda^{b, U}_0}\phi
^{b, E}_U (x).
\end{equation}
By the semigroup property and
(\ref{e10}),
\begin{eqnarray*}
&&\overline{p}^{b, E}_U(t, x, y)\\
&&\qquad= \int_U \overline{p}^{b, E}_U(t/3, x, z)
\int_U \overline{p}^{b, E}_U(t/3, z, w)
\overline{p}^{b, E}_U(t/3, w, y) \xi_E(dw)\xi_E(dz) \\
&&\qquad\le c_3 t^{-d/\alpha}
\int_U \overline{p}^{b, E}_U(t/3, x, z) \xi_E(dz)\int_U
\overline{p}^{b, E}_U(t/3, w, y) \xi_E(dw)\\
&&\qquad= c_3 t^{-d/\alpha} \PP_x ( \tau^{b, E}_U > t/3)
\PP_y ( \wh\tau^{b, E}_U > t/3) \\
&&\qquad\le (9c_3/t^2) t^{-d/\alpha} \E_x [ \tau^{b}_U]
\E_y [ \wh\tau^{b, E}_U].
\end{eqnarray*}
This\vspace*{2pt} together with~(\ref{eiiu}) and~(\ref{eiiu-1}) establishes the
intrinsic ultracontractivity of $\{P^{b, E, U}_t\}$ and
$\{\wh{P}^{b, E, U}_t\}$.
\end{pf}

Applying~\cite{KS2}, Theorem 2.7, we obtain:
\begin{theorem}\label{c2e}
There exist positive constants $c$ and $\nu$ such that
%
%
\begin{equation}\label{h}\qquad
\biggl|\frac{M^{b, E}_Ue^{t \lambda_0^{b, U}}\overline{p}^{b, E}_U(t,
x, y)}{\phi^{b, E}_U(x)\psi^{b, E}_U(y)}-1\biggr| \le ce^{-\nu
t},\qquad
(t, x, y)\in(1, \infty)\times U\times U,
\end{equation}
where $M^{b, E}_U:=\int_U\phi^{b, E}_U(y)\psi^{b, E}_U(y)\xi_E(dy)
\le1$.\vadjust{\goodbreak}
\end{theorem}

Now we can present the following:
\begin{pf*}{Proof of Theorem~\ref{Tmain}\textup{(ii)}}
Assume that the ball $E$ is large enough so that $D\subset
\frac14E$.
Since\vspace*{1pt} $\phi^{b, E}_D =
e^{\lambda_0^{b, D}}
P^{b, D}_1 \phi^{b, E}_D$ and $\psi^{b, E}_D
= e^{\lambda_0^{b, D}} \wh P^{b, E, D}_1 \psi^{b, E}_D$, we have
from Theorem~\ref{Tmain}(i) that on $D$,
%
%
\begin{eqnarray} \label{eed}\qquad
\phi^{b, E}_D(x) &\asymp&\bigl( 1\wedge
\delta_D(x)^{\alpha/2}\bigr) \int_D \bigl( 1\wedge
\delta_D(y)^{\alpha/2}\bigr) \biggl(1\wedge\frac1{|x-y|^{d+\alpha}}
\biggr) \phi^{b, E}_D(y) \,dy \nonumber\\[-8pt]\\[-8pt]
&\asymp&\delta_D(x)^{\alpha/2} \nonumber
\end{eqnarray}
and
%
%
\begin{eqnarray} \label{eed-1}
\psi^{b, E}_D(x)
&\asymp&\bigl( 1\wedge
\delta_D(x)^{\alpha/2}\bigr) \int_D \bigl( 1\wedge
\delta_D(y)^{\alpha/2}\bigr) \biggl(1\wedge\frac1{|x-y|^{d+\alpha}}
\biggr) \nonumber\\
&&\hspace*{78pt}{}\times\frac{h_E(y)}{h_E(x)}\psi^{b, E}_D(y) \,dy
\\
&\asymp&
\delta_D(x)^{\alpha/2}. \nonumber
\end{eqnarray}
Theorem~\ref{tiu},~(\ref{eed}) and~(\ref{eed-1}) imply that
\[
c_t^{-1} \delta_D (x)^{\alpha/2} \delta_D (y)^{\alpha/2} \leq
\overline{p}^{b, E}_D(t, x, y) \leq c_t \delta_D (x)^{\alpha/2}
\delta_D (y)^{\alpha/2}
\]
for every $(t, x, y) \in(0,
\infty)\times D\times D$,
and so
\[
c_1^{-1}c_t^{-1} \delta_D (x)^{\alpha/2} \delta_D
(y)^{\alpha/2} \leq p^{b}_D(t, x, y) \leq c_1 c_t \delta_D
(x)^{\alpha/2}
\delta_D (y)^{\alpha/2}
\]
for every $(t, x, y) \in(0,
\infty)\times D\times D$.

Furthermore, by Theorem~\ref{c2e} and~(\ref{eed}),
there exist $c_2>1$ and $T_1 > 0$ such that for all $(t, x, y)\in
[T_1, \infty)\times D\times D$,
\[
c_2^{-1} e^{-t \lambda_0^{b, D}} \delta_D (x)^{\alpha/2}
\delta_D (y)^{\alpha/2} \leq\overline{p}^{b, E}_D(t, x, y)
\leq c_2 e^{-t \lambda_0^{b, D}} \delta_D
(x)^{\alpha/2} \delta_D (y)^{\alpha/2},
\]
which implies that
\[
c_3^{-1} e^{-t \lambda_0^{b, D}} \delta_D (x)^{\alpha/2}
\delta_D (y)^{\alpha/2} \leq p^{b}_D(t, x, y) \leq c_3
e^{-t \lambda_0^{b, D}} \delta_D (x)^{\alpha/2} \delta_D
(y)^{\alpha/2}.
\]
If $T< T_1$, by Theorem~\ref{Tmain}(i), there is a constant $c_2\geq
1$ such that
\[
c_2^{-1} \delta_D(x)^{\alpha/2} \delta_D (y)^{\alpha/2} \leq
p^b_D(t, x, y) \leq c_2 \delta_D(x)^{\alpha/2} \delta_D
(y)^{\alpha/2}
\]
for every $t\in[T, T_1)$ and $x, y
\in D$.
This establishes Theorem~\ref{Tmain}(ii).
\end{pf*}
\begin{remark}
(i) Using Corollary~\ref{C12} and the argument of the proof of
Lemma~\ref{L20},~(\ref{elpz1}) is, in fact, true for all bounded
open sets $U$ with exterior cone condition.

(ii) In view of Corollary~\ref{C12}, estimate~(\ref{eg1})
and Lemma~\ref{l41},
we can deduce from~(\ref{e49}) by the dominated convergence
theorem that Proposition~\ref{tdh1} holds for general $b$
with $|b|\in\bK_{d, \alpha-1}$.
\end{remark}

\section*{Acknowledgments}
The main results of this paper were
presented by the authors at the \textit{Sixth International Conference
on L\'evy Processes}: \textit{Theory and Applications} held in Dresden,
Germany from July 26 to 30, 2010 and at the \textit{34th Stochastic
Processes and their Applications} conference held in Osaka, Japan
from September 6 to 10, 2010.

K. Bogdan announced at the \textit{Sixth International Conference on
L\'evy Processes: Theory and Applications} in Dresden that he and T.
Jakubowski have also obtained the same sharp estimates on $G^b_D$
in bounded $C^{1,1}$ domains as given in Corollary~\ref{C12} of
this paper. They have also obtained one part (for harmonic functions
of $X^b$ only) of the boundary Harnack principles established in
Theorem~\ref{ubhp} of our paper. Their preprint~\cite{BJ2}
containing these two results (Theorem 1 and Lemma 18 there) appeared
in the \href{http://arxiv.org/abs/arXiv:1009.2471}{arXiv:1009.2471}
on September 14, 2010.

We thank the referees for their helpful comments on the first version
of this paper.


%

%
\printaddresses


\begin{thebibliography}{38}

\bibitem{Bi}
%
\begin{barticle}[mr]
\bauthor{\bsnm{Billingsley},~\bfnm{Patrick}\binits{P.}}
(\byear{1974}).
\btitle{Conditional distributions and tightness}.
\bjournal{Ann. Probab.}
\bvolume{2}
\bpages{480--485}.
\bid{mr={0368095}}
\end{barticle}
%
\endbibitem

\bibitem{BG}
%
\begin{barticle}[mr]
\bauthor{\bsnm{Blumenthal},~\bfnm{R.~M.}\binits{R.~M.}} \AND
\bauthor{\bsnm{Getoor},~\bfnm{R.~K.}\binits{R.~K.}}
(\byear{1960}).
\btitle{Some theorems on stable processes}.
\bjournal{Trans. Amer. Math. Soc.}
\bvolume{95}
\bpages{263--273}.
\bid{issn={0002-9947}, mr={0119247}}
\end{barticle}
%
\endbibitem

\bibitem{BG1}
%
\begin{bbook}[mr]
\bauthor{\bsnm{Blumenthal},~\bfnm{R.~M.}\binits{R.~M.}} \AND
\bauthor{\bsnm{Getoor},~\bfnm{R.~K.}\binits{R.~K.}}
(\byear{1968}).
\btitle{Markov Processes and Potential Theory}.
\bseries{Pure and Applied Mathematics}
\bvolume{29}.
\bpublisher{Academic Press}, \baddress{New York}.
\bid{mr={0264757}}
\end{bbook}
%
\endbibitem

\bibitem{B}
%
\begin{barticle}[mr]
\bauthor{\bsnm{Bogdan},~\bfnm{Krzysztof}\binits{K.}}
(\byear{1997}).
\btitle{The boundary {H}arnack principle for the fractional {L}aplacian}.
\bjournal{Studia Math.}
\bvolume{123}
\bpages{43--80}.
\bid{issn={0039-3223}, mr={1438304}}
\end{barticle}
%
\endbibitem

\bibitem{BGR}
%
\begin{barticle}[mr]
\bauthor{\bsnm{Bogdan},~\bfnm{Krzysztof}\binits{K.}},
\bauthor{\bsnm{Grzywny},~\bfnm{Tomasz}\binits{T.}} \AND
\bauthor{\bsnm{Ryznar},~\bfnm{Micha{\l}}\binits{M.}}
(\byear{2010}).
\btitle{Heat kernel estimates for the fractional {L}aplacian with {D}irichlet
conditions}.
\bjournal{Ann. Probab.}
\bvolume{38}
\bpages{1901--1923}.
\bid{doi={10.1214/10-AOP532}, issn={0091-1798}, mr={2722789}}
\end{barticle}
%
\endbibitem

\bibitem{BJ}
%
\begin{barticle}[mr]
\bauthor{\bsnm{Bogdan},~\bfnm{Krzysztof}\binits{K.}} \AND
\bauthor{\bsnm{Jakubowski},~\bfnm{Tomasz}\binits{T.}}
(\byear{2007}).
\btitle{Estimates of heat kernel of fractional {L}aplacian perturbed by
gradient operators}.
\bjournal{Comm. Math. Phys.}
\bvolume{271}
\bpages{179--198}.
\bid{doi={10.1007/s00220-006-0178-y}, issn={0010-3616}, mr={2283957}}
\end{barticle}
%
\endbibitem

\bibitem{BJ2}
%
\begin{barticle}[auto:STB|2011-03-03|12:04:44]
\bauthor{\bsnm{Bogdan},~\bfnm{K.}\binits{K.}} \AND
\bauthor{\bsnm{Jakubowski},~\bfnm{T.}\binits{T.}}
(\byear{2012}).
\btitle{Estimates of the Green function for the fractional Laplacian
perturbed by gradient}.
\bjournal{Potential Anal.}
\bvolume{36}
\bpages{455--481}.
\end{barticle}
%
\endbibitem

\bibitem{BKN}
%
\begin{barticle}[mr]
\bauthor{\bsnm{Bogdan},~\bfnm{K.}\binits{K.}},
\bauthor{\bsnm{Kulczycki},~\bfnm{T.}\binits{T.}} \AND
\bauthor{\bsnm{Nowak},~\bfnm{Adam}\binits{A.}}
(\byear{2002}).
\btitle{Gradient estimates for harmonic and {$q$}-harmonic functions of
symmetric stable processes}.
\bjournal{Illinois J. Math.}
\bvolume{46}
\bpages{541--556}.
\bid{issn={0019-2082}, mr={1936936}}
\end{barticle}
%
\endbibitem

\bibitem{CKS}
%
\begin{barticle}[mr]
\bauthor{\bsnm{Chen},~\bfnm{Zhen-Qing}\binits{Z.-Q.}},
\bauthor{\bsnm{Kim},~\bfnm{Panki}\binits{P.}} \AND
\bauthor{\bsnm{Song},~\bfnm{Renming}\binits{R.}}
(\byear{2010}).
\btitle{Heat kernel estimates for the {D}irichlet fractional {L}aplacian}.
\bjournal{J. Eur. Math. Soc. (JEMS)}
\bvolume{12}
\bpages{1307--1329}.
\bid{doi={10.4171/JEMS/231}, issn={1435-9855}, mr={2677618}}
\end{barticle}
%
\endbibitem

\bibitem{CKS1}
%
\begin{barticle}[mr]
\bauthor{\bsnm{Chen},~\bfnm{Zhen-Qing}\binits{Z.-Q.}},
\bauthor{\bsnm{Kim},~\bfnm{Panki}\binits{P.}} \AND
\bauthor{\bsnm{Song},~\bfnm{Renming}\binits{R.}}
(\byear{2010}).
\btitle{Two-sided heat kernel estimates for censored stable-like processes}.
\bjournal{Probab. Theory Related Fields}
\bvolume{146}
\bpages{361--399}.
\bid{doi={10.1007/s00440-008-0193-3}, issn={0178-8051}, mr={2574732}}
\end{barticle}
%
\endbibitem

\bibitem{CKSV}
%
\begin{bmisc}[auto:STB|2011-03-03|12:04:44]
\bauthor{\bsnm{Chen},~\bfnm{Z.~Q.}\binits{Z.~Q.}},
\bauthor{\bsnm{Kim},~\bfnm{P.}\binits{P.}},
\bauthor{\bsnm{Song},~\bfnm{R.}\binits{R.}} \AND
\bauthor{\bsnm{Vondra{\v{c}}ek},~\bfnm{Z.}\binits{Z.}}
(\byear{2012}).
\bhowpublished{Boundary Harnack principle for $\Delta+ \Delta
^{\alpha/2}$.
\textit{Trans. Amer. Math. Soc.} To appear}.
\end{bmisc}
%
\endbibitem

\bibitem{CK}
%
\begin{barticle}[mr]
\bauthor{\bsnm{Chen},~\bfnm{Zhen-Qing}\binits{Z.-Q.}} \AND
\bauthor{\bsnm{Kumagai},~\bfnm{Takashi}\binits{T.}}
(\byear{2003}).
\btitle{Heat kernel estimates for stable-like processes on {$d$}-sets}.
\bjournal{Stochastic Process. Appl.}
\bvolume{108}
\bpages{27--62}.
\bid{doi={10.1016/S0304-4149(03)00105-4}, issn={0304-4149}, mr={2008600}}
\end{barticle}
%
\endbibitem

\bibitem{CK2}
%
\begin{barticle}[mr]
\bauthor{\bsnm{Chen},~\bfnm{Zhen-Qing}\binits{Z.-Q.}} \AND
\bauthor{\bsnm{Kumagai},~\bfnm{Takashi}\binits{T.}}
(\byear{2008}).
\btitle{Heat kernel estimates for jump processes of mixed types on metric
measure spaces}.
\bjournal{Probab. Theory Related Fields}
\bvolume{140}
\bpages{277--317}.
\bid{doi={10.1007/s00440-007-0070-5}, issn={0178-8051}, mr={2357678}}
\end{barticle}
%
\endbibitem

\bibitem{CS1}
%
\begin{barticle}[mr]
\bauthor{\bsnm{Chen},~\bfnm{Zhen-Qing}\binits{Z.-Q.}} \AND
\bauthor{\bsnm{Song},~\bfnm{Renming}\binits{R.}}
(\byear{1998}).
\btitle{Estimates on {G}reen functions and {P}oisson kernels for symmetric
stable processes}.
\bjournal{Math. Ann.}
\bvolume{312}
\bpages{465--501}.
\bid{doi={10.1007/s002080050232}, issn={0025-5831}, mr={1654824}}
\end{barticle}
%
\endbibitem

\bibitem{CR}
%
\begin{barticle}[mr]
\bauthor{\bsnm{Chung},~\bfnm{K.~L.}\binits{K.~L.}} \AND
\bauthor{\bsnm{Rao},~\bfnm{K.~Murali}\binits{K.~M.}}
(\byear{1980}).
\btitle{A new setting for potential theory. {I}}.
\bjournal{Ann. Inst. Fourier (Grenoble)}
\bvolume{30}
\bpages{167--198}.
\bid{issn={0373-0956}, mr={0597022}}
\end{barticle}
%
\endbibitem

\bibitem{CW}
%
\begin{bbook}[mr]
\bauthor{\bsnm{Chung},~\bfnm{Kai~Lai}\binits{K.~L.}} \AND
\bauthor{\bsnm{Walsh},~\bfnm{John~B.}\binits{J.~B.}}
(\byear{2005}).
\btitle{Markov Processes, {B}rownian Motion, and Time Symmetry},
\bedition{2nd} ed.
\bseries{Grundlehren der Mathematischen Wissenschaften [Fundamental Principles
of Mathematical Sciences]}
\bvolume{249}.
\bpublisher{Springer}, \baddress{New York}.
\bid{mr={2152573}}
\end{bbook}
%
\endbibitem

\bibitem{CZ}
%
\begin{bbook}[mr]
\bauthor{\bsnm{Chung},~\bfnm{Kai~Lai}\binits{K.~L.}} \AND
\bauthor{\bsnm{Zhao},~\bfnm{Zhong~Xin}\binits{Z.~X.}}
(\byear{1995}).
\btitle{From {B}rownian Motion to {S}chr\"odinger's Equation}.
\bseries{Grundlehren der Mathematischen Wissenschaften [Fundamental Principles
of Mathematical Sciences]}
\bvolume{312}.
\bpublisher{Springer}, \baddress{Berlin}.
\bid{mr={1329992}}
\end{bbook}
%
\endbibitem

\bibitem{CrZ}
%
\begin{barticle}[mr]
\bauthor{\bsnm{Cranston},~\bfnm{M.}\binits{M.}} \AND
\bauthor{\bsnm{Zhao},~\bfnm{Z.}\binits{Z.}}
(\byear{1987}).
\btitle{Conditional transformation of drift formula and potential
theory for
{${1\over2}\Delta+b(\cdot)\cdot\nabla$}}.
\bjournal{Comm. Math. Phys.}
\bvolume{112}
\bpages{613--625}.
\bid{issn={0010-3616}, mr={0910581}}
\end{barticle}
%
\endbibitem

\bibitem{Du}
%
\begin{bbook}[auto:STB|2011-03-03|12:04:44]
\bauthor{\bsnm{Durrett},~\bfnm{R.}\binits{R.}}
(\byear{2005}).
\btitle{Probability: Theory and Examples},
\bedition{3rd} ed.
\bpublisher{Thomson Learning}, \baddress{Belmont, CA}.
\end{bbook}
%
\endbibitem

\bibitem{EK}
%
\begin{bbook}[mr]
\bauthor{\bsnm{Ethier},~\bfnm{Stewart~N.}\binits{S.~N.}} \AND
\bauthor{\bsnm{Kurtz},~\bfnm{Thomas~G.}\binits{T.~G.}}
(\byear{1986}).
\btitle{Markov Processes: Characterization and Convergence}.
\bpublisher{Wiley}, \baddress{New York}.
\bid{doi={10.1002/9780470316658}, mr={0838085}}
\end{bbook}
%
\endbibitem

\bibitem{G}
%
\begin{barticle}[mr]
\bauthor{\bsnm{Getoor},~\bfnm{R.~K.}\binits{R.~K.}}
(\byear{1971}).
\btitle{Duality of {L}\'evy systems}.
\bjournal{Z. Wahrsch. Verw. Gebiete}
\bvolume{19}
\bpages{257--270}.
\bid{mr={0301805}}
\end{barticle}
%
\endbibitem

\bibitem{Gr}
%
\begin{barticle}[mr]
\bauthor{\bsnm{Grigelionis},~\bfnm{B.}\binits{B.}}
(\byear{1973}).
\btitle{The relative compactness of sets of probability measures in
$D[0,\infty)$}.
\bjournal{Lith. Math. J.}
\bvolume{13}
\bpages{576--586}.
\end{barticle}
%
\endbibitem

\bibitem{J}
%
\begin{barticle}[mr]
\bauthor{\bsnm{Jakubowski},~\bfnm{Tomasz}\binits{T.}}
(\byear{2002}).
\btitle{The estimates for the {G}reen function in {L}ipschitz domains
for the
symmetric stable processes}.
\bjournal{Probab. Math. Statist.}
\bvolume{22}
\bpages{419--441}.
\bid{issn={0208-4147}, mr={1991120}}
\end{barticle}
%
\endbibitem

\bibitem{KS}
%
\begin{barticle}[mr]
\bauthor{\bsnm{Kim},~\bfnm{Panki}\binits{P.}} \AND
\bauthor{\bsnm{Song},~\bfnm{Renming}\binits{R.}}
(\byear{2006}).
\btitle{Two-sided estimates on the density of {B}rownian motion with singular
drift}.
\bjournal{Illinois J. Math.}
\bvolume{50}
\bpages{635--688 (electronic)}.
\bid{issn={0019-2082}, mr={2247841}}
\end{barticle}
%
\endbibitem

\bibitem{KS1}
%
\begin{barticle}[mr]
\bauthor{\bsnm{Kim},~\bfnm{Panki}\binits{P.}} \AND
\bauthor{\bsnm{Song},~\bfnm{Renming}\binits{R.}}
(\byear{2007}).
\btitle{Boundary {H}arnack principle for {B}rownian motions with measure-valued
drifts in bounded {L}ipschitz domains}.
\bjournal{Math. Ann.}
\bvolume{339}
\bpages{135--174}.
\bid{doi={10.1007/s00208-007-0110-6}, issn={0025-5831}, mr={2317765}}
\end{barticle}
%
\endbibitem

\bibitem{KS2}
%
\begin{barticle}[mr]
\bauthor{\bsnm{Kim},~\bfnm{Panki}\binits{P.}} \AND
\bauthor{\bsnm{Song},~\bfnm{Renming}\binits{R.}}
(\byear{2008}).
\btitle{On dual processes of non-symmetric diffusions with measure-valued
drifts}.
\bjournal{Stochastic Process. Appl.}
\bvolume{118}
\bpages{790--817}.
\bid{doi={10.1016/j.spa.2007.06.007}, issn={0304-4149}, mr={2411521}}
\end{barticle}
%
\endbibitem

\bibitem{Ku2}
%
\begin{barticle}[mr]
\bauthor{\bsnm{Kulczycki},~\bfnm{Tadeusz}\binits{T.}}
(\byear{1998}).
\btitle{Intrinsic ultracontractivity for symmetric stable processes}.
\bjournal{Bull. Pol. Acad. Sci. Math.}
\bvolume{46}
\bpages{325--334}.
\bid{issn={0239-7269}, mr={1643611}}
\end{barticle}
%
\endbibitem

\bibitem{L1}
%
\begin{bmisc}[mr]
\bauthor{\bsnm{Liao},~\bfnm{Ming}\binits{M.}}
(\byear{1984}).
\bhowpublished{Riesz representation and duality of Markov processes.
Ph.D. dissertation, Dept. Mathematics, Stanford Univ.}
\end{bmisc}
%
\endbibitem

\bibitem{L}
%
\begin{bincollection}[mr]
\bauthor{\bsnm{Liao},~\bfnm{Ming}\binits{M.}}
(\byear{1985}).
\btitle{Riesz representation and duality of {M}arkov processes}.
In \bbooktitle{S\'eminaire de Probabilit\'es, {XIX}, 1983/84}.
\bseries{Lecture Notes in Math.}
\bvolume{1123}
\bpages{366--396}.
\bpublisher{Springer}, \baddress{Berlin}.
\bid{doi={10.1007/BFb0075866}, mr={0889495}}
\end{bincollection}
%
\endbibitem

\bibitem{LZ}
%
\begin{barticle}[mr]
\bauthor{\bsnm{Liu},~\bfnm{Lu~Qin}\binits{L.~Q.}} \AND
\bauthor{\bsnm{Zhang},~\bfnm{Yi~Ping}\binits{Y.~P.}}
(\byear{1990}).
\btitle{Representation of conditional {M}arkov processes}.
\bjournal{J.~Math. (Wuhan)}
\bvolume{10}
\bpages{1--12}.
\bid{issn={0255-7797}, mr={1075304}}
\end{barticle}
%
\endbibitem

\bibitem{P}
%
\begin{barticle}[mr]
\bauthor{\bsnm{Pop-Stojanovi{\'c}},~\bfnm{Z.~R.}\binits{Z.~R.}}
(\byear{1988}).
\btitle{Continuity of excessive harmonic functions for certain diffusions}.
\bjournal{Proc. Amer. Math. Soc.}
\bvolume{103}
\bpages{607--611}.
\bid{doi={10.2307/2047186}, issn={0002-9939}, mr={0943091}}
\end{barticle}
%
\endbibitem

\bibitem{Sc}
%
\begin{bbook}[mr]
\bauthor{\bsnm{Schaefer},~\bfnm{Helmut~H.}\binits{H.~H.}}
(\byear{1974}).
\btitle{Banach Lattices and Positive Operators}.
\bseries{Die Grundlehren der Mathematischen Wissenschaften}
\bvolume{215}
\bpublisher{Springer}, \baddress{New York}.
\bid{mr={0423039}}
\end{bbook}
%
\endbibitem

\bibitem{Sh}
%
\begin{bbook}[mr]
\bauthor{\bsnm{Sharpe},~\bfnm{Michael}\binits{M.}}
(\byear{1988}).
\btitle{General Theory of {M}arkov Processes}.
\bseries{Pure and Applied Mathematics}
\bvolume{133}.
\bpublisher{Academic Press}, \baddress{Boston, MA}.
\bid{mr={0958914}}
\end{bbook}
%
\endbibitem

\bibitem{So}
%
\begin{barticle}[mr]
\bauthor{\bsnm{Song},~\bfnm{Renming}\binits{R.}}
(\byear{2004}).
\btitle{Estimates on the {D}irichlet heat kernel of domains above the
graphs of
bounded {$C\sp{1,1}$} functions}.
\bjournal{Glas. Mat. Ser. III}
\bvolume{39}
\bpages{273--286}.
\bid{doi={10.3336/gm.39.2.09}, issn={0017-095X}, mr={2109269}}
\end{barticle}
%
\endbibitem

\bibitem{SW}
%
\begin{barticle}[mr]
\bauthor{\bsnm{Song},~\bfnm{Renming}\binits{R.}} \AND
\bauthor{\bsnm{Wu},~\bfnm{Jang-Mei}\binits{J.-M.}}
(\byear{1999}).
\btitle{Boundary {H}arnack principle for symmetric stable processes}.
\bjournal{J. Funct. Anal.}
\bvolume{168}
\bpages{403--427}.
\bid{doi={10.1006/jfan.1999.3470}, issn={0022-1236}, mr={1719233}}
\end{barticle}
%
\endbibitem

\bibitem{So1}
%
\begin{barticle}[mr]
\bauthor{\bsnm{Song},~\bfnm{Ren~Ming}\binits{R.~M.}}
(\byear{1995}).
\btitle{Feynman--{K}ac semigroup with discontinuous additive functionals}.
\bjournal{J.~Theoret. Probab.}
\bvolume{8}
\bpages{727--762}.
\bid{doi={10.1007/BF02410109}, issn={0894-9840}, mr={1353551}}
\end{barticle}
%
\endbibitem

\bibitem{S}
%
\begin{barticle}[mr]
\bauthor{\bsnm{Shur},~\bfnm{M.~G.}\binits{M.~G.}}
(\byear{1977}).
\btitle{Dual {M}arkov processes}.
\bjournal{Theory Probab. Appl.}
\bvolume{22}
\bpages{264--278}.
\bid{issn={0040-361X}, mr={0443108}}
\end{barticle}
%
\endbibitem

\end{thebibliography}
\end{document}